\documentclass{article}

\usepackage{vmargin,epsfig}
\usepackage[francais]{babel}
\usepackage[latin1]{inputenc}
\usepackage[tbtags]{amsmath}
\usepackage{amsthm}
\usepackage{amssymb}
\usepackage{pifont}
\usepackage[all]{xy}
\usepackage{wrapfig}
\usepackage{calc}

\setmarginsrb{2.2cm}{1.5cm}{2.2cm}{1.5cm}{1cm}{0cm}{1cm}{1.5cm}

\newtheorem{theo}{Théorème}[subsection]
\newtheorem{lemme}[theo]{Lemme}
\newtheorem{prop}[theo]{Proposition}
\newtheorem{cor}[theo]{Corollaire}
\theoremstyle{definition}
\newtheorem{deftn}[theo]{Définition}

\def\leq{\leqslant}
\def\geq{\geqslant}

\def\Ax#1{\textrm{\bf (Ax{#1})}}

\def\N{\mathbb{N}}
\def\Z{\mathbb{Z}}
\def\Q{\mathbb{Q}}

\def\F{\mathbb{F}}
\def\O{\mathcal{O}}

\def\pa#1{\left(#1\right)}

\def\brac#1{\left< #1 \right>}

\def\epsilon{\varepsilon}

\def\calA{\mathcal{A}}
\def\calB{\mathcal{B}}
\def\calC{\mathcal{C}}
\def\calF{\mathcal{F}}

\def\calM{\mathcal{M}}

\def\calN{\mathcal{N}}
\def\calR{\mathcal{R}}

\def\calX{\mathcal{X}}
\def\frakM{\mathfrak{M}}

\def\calK{\mathcal{K}}

\def\si{\mathfrak{S}}

\def\id{\text{\rm id}}
\def\im{\text{\rm im}\:}
\def\coim{\text{\rm coim}\:}
\def\ker{\text{\rm ker}\:}
\def\coker{\text{\rm coker}\:}
\def\hom{\text{\rm Hom}}
\def\End{\text{\rm End}}

\def\Fil{\text{\rm Fil}}

\def\sep{\text{\rm sep}}

\def\pr{\text{\rm pr}}

\def\st{\text{\rm st}}
\def\cris{\text{\rm cris}}

\def\qst{\text{\rm qst}}

\def\Tst{T_\st}
\def\Tkst{T_{K\text{-}\st}}
\def\Tlst{T_{L\text{-}\st}}
\def\Tsthat{\hat T_\st}
\def\TstE{T_{\st,E}}
\def\Mst{M_\st}

\def\Mlst{M_{L\text{-}\st}}

\def\Tqst{T_\qst}

\def\Tsi{T_{\tilde \si}}

\def\Rep{\text{\rm Rep}}

\def\dd{\text{\rm dd}}
\def\pModpN#1{\text{\rm 'Mod}^{\phi,N}_{/#1}}
\def\ModpN#1{\text{\rm Mod}^{\phi,N}_{/#1}}
\def\ModpNdd#1{\text{\rm Mod}^{\phi,N,\dd}_{/#1}}
\def\Modst#1{\text{\rm Mod}^{\st}_{/#1}}
\def\Modcris#1{\text{\rm Mod}^{\cris}_{/#1}}
\def\pFilpN#1{\text{\rm 'Fil}^{\phi,N}_{/#1}}
\def\pGenpN#1{\text{\rm 'Gen}^{\phi,N}_{/#1}}
\def\pUnipN#1{\text{\rm 'Uni}^{\phi,N}_{/#1}}
\def\pRedpN#1{\text{\rm 'R\'ed}^{\phi,N}_{/#1}}
\def\RedpN#1{\text{\rm R\'ed}^{\phi,N}_{/#1}}
\def\MaxpN#1{\text{\rm Max}^{\phi,N}_{/#1}}
\def\MinpN#1{\text{\rm Min}^{\phi,N}_{/#1}}

\def\pModp#1{\text{\rm 'Mod}^{\phi}_{/#1}}
\def\Modp#1{\text{\rm Mod}^{\phi}_{/#1}}
\def\pFilp#1{\text{\rm 'Fil}^{\phi}_{/#1}}
\def\pGenp#1{\text{\rm 'Gen}^{\phi}_{/#1}}
\def\pUnip#1{\text{\rm 'Uni}^{\phi}_{/#1}}
\def\pRedp#1{\text{\rm 'R\'ed}^{\phi}_{/#1}}
\def\Redp#1{\text{\rm R\'ed}^{\phi}_{/#1}}

\def\Ind{\text{\rm Ind}}
\def\Gen{\text{\rm Gen}}
\def\Red{\text{\rm Réd}_\qst}
\def\RedN{\text{\rm Réd}_\st}
\def\Nil{\text{\rm Nil}_\qst}
\def\NilN{\text{\rm Nil}_\st}
\def\Mod{\text{\rm Mod}_\qst}
\def\ModN{\text{\rm Mod}_\st}

\def\Max{\text{\rm Max}}
\def\Maxst{\Max^\st}
\def\max{\text{\rm max}}

\def\Min{\text{\rm Min}}
\def\Minst{\Min^\st}
\def\min{\text{\rm min}}

\def\toinj{\hookrightarrow}
\def\tosurj{\twoheadrightarrow}

\def\Ob{}
\def\Comp{\text{\rm Comp}}
\def\Fin{\text{\rm Fin}}
\def\Oub{\text{\rm Oub}}
\def\Gal{\text{\rm Gal}}

\def\nr{\text{\rm nr}}
\def\ppcm{\text{\sc Ppcm}}

\newenvironment{axiome}[1]{%
  \begin{center}
  \hspace{1em}\begin{minipage}{\linewidth-2cm}%
  \hspace{-1.5em}\Ax{#1}\hspace{0.1em}\it}{%
  \end{minipage}\end{center}}

\def\tomonth{\ifcase\month\or
  Janvier\or Février\or Mars\or Avril\or Mai\or Juin\or
  Juillet\or Août\or Septembre\or Octobre\or Novembre\or Décembre\fi
  \space \number\year}

\title{$\F_p$-répresentations semi-stables}
\author{Xavier Caruso}
\date\tomonth

\begin{document}

\maketitle

\tableofcontents

\bigskip

\noindent \hrulefill

\bigskip

Soit $p$ un nombre premier. Soient $k$ un corps parfait de 
caractéristique $p$, $W = W(k)$ l'anneau des vecteurs de Witt à 
coefficients dans $k$, $K_0$ son corps des fractions et $K$ une 
extension finie de $K_0$ totalement ramifiée de degré $e$. Notons $G_K$ 
le groupe de Galois absolu de $K$ et fixons un entier $r \in \{0, 
\ldots, p-2\}$. À partir de ces données, Breuil a défini dans 
\cite{breuil-ens} et \cite{breuil-group} une certaine catégorie de 
modules de torsion, notée $\ModpN {\tilde S}$ dans cet article (et dont 
la définition est rappelée en \ref{subsec:defbreuil}). Celle-ci permet 
\emph{via} un foncteur $\Tst$ de construire certaines 
$\F_p$-représentations du groupe $G_K$. Les représentations ainsi 
obtenues sont intéressantes pour au moins deux raisons : d'une part, 
elles contiennent un grand nombre de représentations de nature 
géométrique (données typiquement par la cohomologie étale des variétés), 
et d'autre part elles regroupent tous les quotients annulés par $p$ de 
deux réseaux à l'intérieur d'une même représentation semi-stable à poids 
de Hodge-Tate compris entre $0$ et $r$. Ainsi la compréhension de cette 
catégorie et du foncteur associé permet-elle d'obtenir diverses 
informations générales pouvant trouver des applications variées (voir 
par exemple \cite{caruso-crelle}, \cite{hattori}, \cite{caruso-savitt}).

\medskip

Lorsque $er < p-1$, la situation est plutôt bien comprise : on sait, par
les résultats de \cite{caruso-crelle}, que la catégorie $\ModpN {\tilde
S}$ est abélienne et que le foncteur $\Tst$ est exact et pleinement
fidèle. Ainsi, en un certain sens, on ramène l'étude de ces objets
compliqués que sont les représentations galoisiennes à des questions
d'algèbre (semi-)linéaire d'apparence plus simple. Cependant, lorsque
$er \geq p-1$, les deux résultats essentiels cités précédemment sont
facilement mis en défaut. Le but de cet article est de dégager la
structure générale de la catégorie $\ModpN {\tilde S}$ et du foncteur
$\Tst$ : on prouve essentiellement que $\ModpN {\tilde S}$ admet une
sous-catégorie pleine\footnote{Lorsque $er <p-1$, on a $\MaxpN {tilde S}
= \ModpN {\tilde S}$ (autrement dit, tout objet est maximal, et on
retrouve la situation établie dans \cite{caruso-crelle}.} $\MaxpN
{\tilde S}$ (dont les objets sont qualifiés de \emph{maximaux}) qui est
abélienne et en restriction à laquelle le foncteur $\Tst$ est exact et
pleinement fidèle. De plus, on construit un foncteur $\Max : \ModpN
{\tilde S} \to \MaxpN {\tilde S}$ qui permet de réaliser $\MaxpN {\tilde
S}$ également comme un quotient de $\ModpN {\tilde S}$. On a en outre
$\Tst \circ \Max = \Tst$, ce qui assure en particulier que la catégorie
$\MaxpN {\tilde S}$ est suffisamment grosse pour être intéressante ; en
tout cas, elle capture autant de représentations galoisiennes que ne le
fait $\ModpN {\tilde S}$

\medskip

Afin de présenter les résultats obtenus de façon quelque peu
systématique, nous avons choisi d'isoler dans une première section toute
une axiomatique dont l'aboutissement est la notion de \emph{pylonet} qui
sera centrale dans la suite du texte, puisque c'est elle qui décrit avec
précision la structure de $\ModpN {\tilde S}$ et $\Tst$. Cette première
section est donc tout à fait générale et abstraite : on s'y borne
essentiellement à mener certains développements sur les catégories
fibrées.

Avec la deuxième section, on entre dans le vif du sujet : on donne les 
définitions des catégories de modules et les foncteurs évoquées 
précédemment puis on montre que $\Tst$ définit un \emph{pylonet} 
(théorème \ref{theo:tstpylonet}). La construction de la sous-catégorie 
$\MaxpN {\tilde S}$ et du foncteur $\Max$ dont il a été question 
auparavant découle alors de l'étude générale de la première partie. La 
section se termine par la preuve de la pleine fidélité de $\Tst$ en 
restriction à $\MaxpN {\tilde S}$ (théorème \ref{theo:pleinfid}). On 
notera que les méthodes de démonstration sont radicalement différentes 
de celles utilisées dans \cite{caruso-crelle} ; elles sont, selon nous, 
beaucoup plus conceptuelles, et semblent également avoir une portée bien 
plus importante.

Dans la troisième section, nous nous efforçons de rendre plus concrètes
les constructions faites dans la section \ref{sec:application},
notamment en ce qui concerne les noyaux, les conoyaux et le foncteur
$\Max$. Pour cela, on est amené à introduire la notion d'objets
$\Tst$-réduits. On montre que leur catégorie, notée $\RedpN {\tilde S}$,
est équivalente à $\ModpN {\tilde S}$, puis on explique comment de
nombreuses constructions se réalisent dans $\RedpN {\tilde S}$. On
démontre également une formule de réciprocité totalement explicite (même
si elle reste un peu compliquée à exprimer) qui permet de retrouver
l'unique objet de $\MaxpN {\tilde S}$ correspondant à une réprésentation
donnée. Combiné au résultat de \cite{caruso-invent}, cela donne en
particulier une recette pour calculer la cohomologie log-cristalline de
la fibre spéciale d'une variété $\calX$ à réduction semi-stable sur
$\O_K$ en fonction de la cohomologie étale $p$-adique de $\calX_{\bar
K}$. (Dans la référence précédente, on obtenait simplement une formule
pour aller dans l'autre sens.)

Dans la dernière section, nous poursuivons notre investigation
principalement en nous intéressant à certaines variantes de la catégorie
$\ModpN {\tilde S}$ obtenues en introduisant des coefficients ou des
données de descente. Dans les deux cas, on montre que l'on obtient
encore des pylonets et que la restriction de $\Tst$ aux objets maximaux
correspondant est à nouveau exacte et pleinement fidèle. Soulignons que
ces variantes interviennent de façon cruciale dans \cite{gee} pour
étudier certains problèmes de modularité de représentations galoisiennes
liés à la généralisation par Buzzard, Diamond et Jarvis de la conjecture
de modularité de Serre (voir \cite{bdj} pour l'énoncé de cette
généralisation). Bien que n'étant pas logiquement nécessaire, il nous
semble que le cadre théorique fourni par cet article éclaire de façon
spectaculaire les calculs de \cite{gee}, \S 3.4 (voir aussi
\cite{caruso-gtgee} à ce sujet).

Nous étudions ensuite une troisième variante, qui est celle que l'on
obtient lorsque l'on ne considère que les objets qui s'écrivent comme
quotients de deux modules fortement divisibles, et donc qui
correspondent à des quotients de deux réseaux dans une représentation
semi-stable. Encore une fois, on obtient des résultats analogues : le
foncteur $\Tst$ définit un pylonet et sa restriction à la sous-catégorie
des objets maximaux est exacte et pleinement fidèle. Finalement, on
donne une description complète des objets simples de $\MaxpN {\tilde S}$
(et de sa variante avec coefficients) lorsque le corps résiduel $k$ est
algébriquement clos.

\section{Notion de pylonet}
\label{sec:axiomes}

Ce premier chapitre est très général et très formel : on y développe une
certaine axiomatique de ce que l'on appelle des pylonets et qui sont des
\og catégories fibrées en sup-semi-treillis satisfaisant la condition de
chaîne croissante \fg. Les exemples et applications seront donnés dans
les chapitres ultérieurs, lorsque l'on s'intéressa plus précisément aux
représentations $p$-adiques.

\medskip

Dans la suite si $\calC$ est une catégorie, on notera parfois $C \in 
\calC$ pour dire que $C$ est un objet de $\calC$. La donnée de départ de 
notre travail est celle de deux catégories $\calC$ et $\calA$ et d'un 
foncteur $T : \calC \to \calA$ (sur lequel pour l'instant on ne fait 
aucune hypothèse) que l'on s'efforcera de considérer comme une 
fibration.

\subsection{Catégories fibres}
\label{subsec:fibres}

Soit $T : \calC \to \calA$ un foncteur covariant\footnote{Dans les 
applications, le foncteur $T$ sera en réalité plutôt contravariant. 
Cependant, quitte à remplacer $\calA$ par sa catégorie opposée, cela ne 
modifie en rien la théorie. Nous préférons donc, pour ce premier 
chapitre, ne pas introduire cette complication inutile.}. Fixons $A$ un 
objet de $\calA$. Il y a deux définitions naturelles pour la fibre de 
$T$ au-dessus de $A$ qui sont : 
\begin{itemize} 
\item la catégorie $F_A$ dont les objets sont les objets $C$ de $\calC$ 
tels que $T(C) = A$ (ceci est une vraie égalité !) et dont les 
morphismes sont les flèches de $\calC$ qui s'envoient sur l'identité de 
$A$ par le foncteur $T$ ; 
\item la catégorie $\calF_A$ dont les objets sont les couples $(C,f)$ 
où $C \in \Ob \calC$ et $f : T(C) \to A$ est un \emph{isomorphisme}, un 
morphisme de $(C,f)$ dans $(C',f')$ étant la donnée de $g \in 
\hom_\calC (C,C')$ vérifiant $f = f' \circ T(g)$.
\end{itemize}

\medskip

L'objet $A$ étant toujours fixé, les deux catégories $F_A$ et $\calF_A$ 
ne sont en général pas équivalentes. Précisément, on dispose d'un 
foncteur pleinement fidèle $F_A \to \calF_A$ défini par $C \mapsto (C, 
\id)$. L'essentielle surjectivité n'est par contre pas automatique, mais
équivaut par définition à l'axiome suivant :
\begin{axiome}{0}
Pour tout $C \in \Ob \calC$ et tout \emph{isomorphisme} (dans $\calA$)
$f : T(C) \to A'$, il existe un \emph{isomorphisme} (dans $\calC$) $g : 
C \to C'$ tel que $T(g) = f$ (et donc $T(C') = A'$).
\end{axiome}
On remarquera que cet axiome est une version très affaiblie de l'axiome
usuel de changement de base qui apparaît par exemple dans la théorie des
champs (algébriques). Nous ne pouvons nous permettre dans cet article de
supposer l'axiome usuel de changement de base, car il sera très loin 
d'être satisfait dans les exemples que nous souhaitons traiter.

\medskip

Il faut remarquer que \Ax 0 n'est pas du tout contraignant.
En effet, s'il n'est pas vérifié, il est toujours possible de 
remplacer $\calC$ par une catégorie \emph{équivalente} $\Comp(\calC,T)$ 
pour laquelle l'axiome est satisfait\footnote{Après ce remplacement, 
les fibres $\calF_A$ restent inchangées contrairement aux $F_A$. La 
première notion de fibre que nous évoquions n'est donc pas robuste, dans 
le sens où elle dépend de $\calC$ à l'intérieur même d'une classe 
d'équivalence de catégories (ou plus exactement de fibrations). Plus 
précisément, les $F_A$ obtenus dans une de ces classes admettent, en un 
certain sens, un élément maximal qui n'est autre que $\calF_A$.}. Cette 
catégorie $\Comp(\calC,T)$ est obtenue comme suit :
\begin{itemize}
\item ses objets sont les triplets $(C, A, f)$ où $C \in \Ob \calC$, $A 
\in \Ob \calA$ et $f : T(C) \to A$ est un isomorphisme ;
\item un morphisme de $(C, A, f)$ dans $(C', A', f')$ est la donnée de
deux morphismes $g \in \hom_\calC (C, C')$ et $h \in \hom_\calA (A, A')$ 
faisant commuter le diagramme suivant :
$$\xymatrix @R=15pt @C=40pt { 
T(C) \ar[r]^-{f}_-{\sim} \ar[d]_-{T(g)} & A \ar[d]^-{h} \\
T(C') \ar[r]^-{f'}_-{\sim} & A' }$$
\end{itemize}
Il est facile de vérifier que le foncteur $\calC \to \Comp(\calC,T)$, $C 
\mapsto (C, T(C), \id)$ est une équivalence de catégories. De plus, $T 
: \calC \to \calA$ se factorise par $\Comp(\calC,T)$ grâce au foncteur
précédent et au foncteur, que nous notons encore $T$, $\Comp(\calC,T) \to 
\calA$, $(C, A, f) \mapsto A$. Ce dernier vérifie l'axiome \Ax 0.

\medskip

Dans la suite de cette section, pour simplifier les écritures (par
exemple ceci nous permettra de travailler avec les fibres $F_A$ au lieu
de $\calF_A$), nous supposerons très fréquemment l'axiome \Ax 0. Malgré
tout, le lecteur doit garder à l'esprit que ce n'est pas du tout
essentiel, et que \emph{tous les résultats obtenus ne faisant pas
intervenir de véritables égalités entre objets (mais seulement des
isomorphismes) demeurent vrais sans aucune modification si l'hypothèse
\Ax 0 est relachée}. En réalité, dans les applications que l'on va
développer dans les sections suivantes, l'axiome \Ax 0 ne sera que très
rarement satisfait.

\subsection{Le foncteur $\Max$}
\label{subsec:max}

En supplément de \Ax 0, nous introduisons les trois axiomes suivants :
\begin{axiome}{1}
Le foncteur $T$ est fidèle.
\end{axiome}
\begin{axiome}{2}
Les catégories $\calC$ et $\calA$ admettent des sommes amalgamées, et 
le foncteur $T$ y est compatible.
\end{axiome}
\begin{axiome}{3}
Pour tout $A \in \Ob \calA$, soit la catégorie $\calF_A$ est vide, soit 
elle admet un objet final.
\end{axiome}
Dans ce paragraphe, nous supposons simplement \Ax 0, \Ax 2 et \Ax 3.
Nous avons préféré introduire \Ax 1 dès à présent car, comme nous allons
le voir, il joue déjà un rôle particulier dans la situation que nous
allons présenter.

\medskip

Nous construisons un foncteur $\Max : \calC \to \calC$ comme suit. Pour
tout objet $A$ de $\calA$ dont la fibre est non vide, choisissons un 
objet final $\Fin(A)$ de la catégorie $F_A$. Sur les objets, le foncteur 
$\Max$ est défini par $\Max(C) = \Fin(T(C))$. Le fait que $C$ soit un 
objet de la fibre $F_{T(C)}$ fournit un morphisme canonique 
$\iota_{\max}^{C} : C \to \Max(C)$ dans la catégorie $\calC$ vérifiant 
$T(\iota_{\max}^{C}) = \id_{T(C)}$. Il reste à définir $\Max(f)$ lorsque 
$f : C_1 \to C_2$ est un morphisme dans $\calC$. Considérons pour cela 
$C'_2 = C_2 \oplus_{C_1} \Max(C_1)$ la somme amalgamée du diagramme :
$$\xymatrix @C=40pt {
\Max(C_1) \\
C_1 \ar[u]^-{\iota_{\max}^{C_1}} \ar[r]^-{f} & C_2 } $$
Notons $\iota' : C_2 \to C'_2$ et $f' : \Max(C_1) \to C'_2$ les 
morphismes correspondants. Comme $T$ est compatible aux sommes 
amalgamées, quitte à modifier $C'_2$ par un objet isomorphe (ce 
que l'on peut faire par l'axiome \Ax 0) $T(C'_2) = T(C_2)$, $T(\iota') = 
\id_{T(C_2)}$ et $T(f') = T(f)$. On en déduit que $C'_2$ est un objet de 
$F_{T(C_2)}$, d'où on obtient le morphisme $\iota_{\max}^{C'_2} : C'_2 \to 
\Max(C_2)$. Le morphisme $\Max(f)$ recherché s'obtient alors comme la 
composée $\iota_{\max}^{C'_2} \circ f'$. Il vérifie $T(\Max(f)) = T(f)$.

\begin{lemme}
\label{lem:maxmorph}
Soit 
$$\xymatrix @C=40pt @R=20pt {
C'_1 \ar[r]^-{g} & C'_2 \\
C_1 \ar[u]^-{h_1} \ar[r]^-{f} & C_2 \ar[u]_-{h_2} }$$
un diagramme commutatif dans $\calC$ tel que $T(C_1) = T(C'_1)$, $T(h_1) 
= \id$ et $T(C_2) = T(C'_2)$, $T(h_2) = \id$. Alors $\Max(f) = \Max(g)$.
\end{lemme}

\begin{proof}
Remarquons tout d'abord que l'hypothèse assure que $\Max(C_1) = \Max(C'_1)$ 
et $\Max(C_2) = \Max(C'_2)$ de sorte que les morphismes $\Max(f)$ et
$\Max(g)$ ont bien même source et même but. Considérons le diagramme 
commutatif
$$\xymatrix @C=50pt {
 & \Max(C_2) \\
\Max(C'_1) \ar[r]^-{g'} & C'_2 \oplus_{C'_1} \Max(C'_1) \ar[u]_-{\iota'} \\
\Max(C_1) \ar[r]^-{f'} \ar@{=}[u] & C_2 \oplus_{C_1} \Max(C_1) 
\ar[u]^-{h_2 \oplus \id} \ar@/_2cm/[uu]_-{\iota} }$$
où $f'$ et $g'$ sont définis comme précédemment et où $\iota$ et 
$\iota'$ sont les morphismes canoniques d'un objet dans son $\Max$. 
Ainsi par définition, $\Max(f) = \iota \circ f'$ et $\Max(g) = \iota' 
\circ g'$. Par ailleurs, comme il y a par définition un unique morphisme 
dans un objet final et que $(h_2 \oplus \id)$, $\iota$ et $\iota'$ sont 
des flèches dans la catégorie $F_{T(C_2)}$, on a nécessairement $\iota = 
\iota' \circ (h_2 \oplus \id)$. Il s'ensuit
$$\Max(f) = \iota \circ f' = \iota' \circ (h_2 \oplus \id) \circ f' = 
\iota' \circ g' = \Max(g)$$
comme annoncé.
\end{proof}

\noindent
{\it Remarque.}
Sous \Ax 1, on remarque que $\Max(f)$ est l'unique morphisme tel que 
$T(\Max(f)) = T(f)$, ce qui permet de simplifier la preuve du lemme
précédent dans ce cas.

\begin{cor}
\label{cor:maxfoncteur}
La construction $\Max$ définit un foncteur $\calC \to \calC$,
et la collection des morphismes $(\iota_{\max}^{C})_{C \in \Ob \calC}$ 
définit une transformation naturelle $\iota_\max$ entre le foncteur 
identité et $\Max$.
\end{cor}

\begin{proof}
Le seul point qu'il reste à prouver est la compatibilité de $\Max$ à la 
composition des morphismes. Considérons pour cela deux morphismes 
composables $f$ et $g$. Par le lemme \ref{lem:maxmorph}, on a 
$\Max(f \circ g) = \Max(\Max(f) \circ \Max(g))$. Or $\Max(f) \circ
\Max(g)$ est un morphisme entre objets de l'image de $\Max$, et on
vérifie immédiatement sur la définition que $\Max$ ne modifie pas un
tel morphisme. Le corollaire s'ensuit.
\end{proof}

\subsubsection*{Digression sur les problèmes de logique}

Dans la construction précédente, on a eu besoin de choisir, pour tout $A 
\in \Ob \calA$, un objet final dans $F_A$. Étant donné que $\Ob \calA$ 
n'est \emph{a priori} pas un ensemble, on peut se demander dans quelle 
mesure, ce choix est légitime. Dans le cas général, il semble délicat de 
justifier cette opération sans introduire la théorie des univers de 
Grothendieck ou des considérations analogues.

\medskip

Toutefois, il est deux situations dans lesquelles on peut envisager des 
palliatifs satisfaisants. La première est bien entendu celle où $\calC$
est une petite catégorie, auquel cas il est suffisant d'invoquer 
l'axiome du choix. En pratique, si les catégories que l'on considère
ne sont pas petites, il sera néanmoins presque toujours possible de 
modifier les définitions pour les remplacer par des petites catégories. 
Ainsi, les problèmes de logique sous-jacents ne sont pas de véritables 
obstacles lorsque l'on envisage les applications.

\medskip

Malgré tout, si l'on suppose \Ax 1, il est possible de mener les
constructions précédentes sans même avoir recours à l'axiome du choix,
quitte à remplacer $\calC$ par une catégorie équivalente. Soit $\bar
\calC$ la catégorie dont les objets sont l'union disjointe
\begin{itemize}
\item des objets de $\calC$ qui ne sont pas des objets maximaux dans
leur fibre et
\item des objets de $\calA$ qui sont dans l'image de $T$.
\end{itemize}
La définition des morphismes est un peu plus délicate, et utilise
\Ax 1 (du moins si l'on souhaite se passer de l'axiome du choix). 
Soient $\bar C_1$ et $\bar C_2$ deux objets de $\bar \calC$. Si $\bar C_1$
est dans $\calC$, on pose $C_1 = \bar C_1$ ; sinon, on désigne par
$C_1$ un objet final (quelconque) de $F_{\bar C_1}$. On définit de 
même $C_2$. On pose :
$$\hom_{\bar \calC} (\bar C_1, \bar C_2) = \text{image} \big[\hom_\calC 
(C_1,C_2) \to \hom_\calA (T(C_1), T(C_2))\big].$$
Il s'agit de montrer que la quantité du membre de droite ne dépend pas
des choix de $C_1$ et $C_2$, lorsqu'il y a effectivement plusieurs choix
pour ces objets, c'est-à-dire lorsque $\bar C_1$ ou $\bar C_2$ est un
objet de $\calA$. On ne traite que le cas de $C_1$, celui de $C_2$ étant
analogue. Supposons donc que $C_1$ et $C'_1$ soient deux objets finaux
de la même fibre $F_A$. Alors, il existe un (unique) morphisme $f : 
C_1 \to C'_1$ tel que $T(f) = \id$. On a alors le diagramme commutatif 
suivant :
$$\xymatrix @C=40pt {
\hom_\calC (C_1,C_2) \ar[r] & \hom_\calA (T(C_1), T(C_2)) 
\ar@{=}[d] \\
\hom_\calC (C'_1,C_2) \ar[r] \ar[u]^-{f^\star} & \hom_\calA 
(T(C'_1), T(C_2)) }$$
qui permet de conclure.

La catégorie $\bar \calC$ est reliée aux données précédentes, notamment
grâce à un foncteur $F : \calC \to \bar \calC$ défini comme suit.
À un objet non final, il associe le même objet, alors qu'à un objet
final, il associe son image sous $T$. Sur les morphismes, il est 
donné par la corestriction du morphisme $\hom_\calC (C_1,C_2) \to 
\hom_\calA (T(C_1), T(C_2))$ déduit de $T$. Comme $T$ est fidèle, ce
dernier est par définition injectif, et donc la corestriction considérée 
est une bijection. Ceci assure que $F$ est pleinement fidèle. Par 
ailleurs, on vérifie facilement qu'il est aussi essentiellement 
surjectif. Ainsi $F$ est une équivalence de catégories.

La catégorie $\bar \calC$ permet aussi de factoriser $T$ : la définition 
même de la relation d'équivalence sur les objets de $\calC$ montre que 
la factorisation existe bien au niveau des objets, alors qu'au niveau 
des morphismes, cela découle de la pleine fidélité de $F$. (Notez que 
l'on ne peut pas dire plus simplement que cette factorisation est 
obtenue en considérant un quasi-inverse de $F$, puisqu'une telle 
construction utilise l'axiome du choix, ce qui est précisément ce que 
l'on souhaite éviter.) Notons $\bar T : \bar \calC \to \calA$ le 
foncteur obtenu. Il est facile de vérifier que la fibration $\bar T$ 
vérifie encore les axiomes \Ax 0, \Ax 1, \Ax 2 et \Ax 3, et que, par 
construction, l'objet maximal de chaque fibre non vide $F_A$ est 
uniquement déterminé (à rien près). Il n'y a donc plus besoin de 
l'axiome du choix pour définir $\Fin(A)$, ni donc le foncteur $\Max$.

\subsubsection*{La catégorie $\Max(\calC)$}

\begin{deftn}
Un objet $C \in \calC$ est dit \emph{maximal} si le morphisme 
$\iota_{\max}^{C} : C \to \Max(C)$ est un isomorphisme.
\end{deftn}

\begin{prop}
\label{prop:projection}
Le foncteur $\Max$ est idempotent, \emph{i.e.} $\Max \circ \Max = \Max$
(sur les objets et sur les morphismes).

L'image essentielle de $\Max$ est la sous-catégorie pleine de $\calC$
formée des objets maximaux. On la note $\Max(\calC)$.
\end{prop}

\begin{proof}
Clair d'après les définitions.
\end{proof}

\noindent
{\it Remarque.} L'image de $\Max$ peut être strictement plus petite que 
$\Max(\calC)$. Cela ne se produit toutefois pas si $\Ax 0$ est vérifié.

\bigskip

Nous prouvons à présent plusieurs propriétés de la catégorie 
$\Max(\calC)$ qui découlent toutes presque directement des définitions. 
Nous commençons pour cela par un lemme important.

\begin{lemme}
\label{lem:maxisom}
Soit $f : C \to C'$ un morphisme dans $\calC$. Alors $T(f)$ est un 
isomorphisme si, et seulement si $\Max(f)$ en est un.
\end{lemme}

\begin{proof}
Si $T(f)$ est un isomorphisme, quitte à remplacer $C'$ par un 
objet isomorphe, l'axiome \Ax 0 nous autorise à supposer que $T(f) = 
\id$. Alors $C$ et $C'$ sont deux objets d'une même fibre, et par 
définition $\Max(C) = \Max(C')$ et $\Max(f)$ n'est autre que l'identité 
entre ces deux objets.

Réciproquement si $\Max(f)$ est un isomorphisme, $T(\Max(f)) = T(f)$ en
est un aussi.
\end{proof}

\begin{prop}
\label{prop:adjonction}
Le foncteur $\Max : \calC \to \Max(\calC)$ est un adjoint à gauche du
foncteur d'inclusion $\Max(\calC) \to \calC$.
\end{prop}

\begin{proof}
Soient $C \in \Ob \calC$ et $M \in \Ob \Max(\calC)$. Nous voulons
exhiber une identification canonique entre $\hom_\calC (C,M)$ et
$\hom_\calC (\Max(C), M)$. Or, on dispose d'une application $\hom_\calC
(C,M) \to \hom_\calC (\Max(C), M)$ donnée par le foncteur $\Max$
(puisque par la proposition \ref{prop:projection}, $\Max(M)$ est
canoniquement isomorphe à $M$ \emph{via} $\iota_{\max}^{M}$) et d'une application
$\hom_\calC (\Max(C),M) \to \hom_\calC (C, M)$ obtenue en composant par
le morphisme canonique $\iota_{\max}^{C} : C \to \Max(C)$. On vérifie facilement
en utilisant $\Max(\iota_{\max}^{C}) = \id_{\Max(C)}$ que les deux applications
précédentes sont inverses l'une de l'autre.
\end{proof}

\begin{prop}
\label{prop:localisation}
Le foncteur $\Max : \calC \to \Max(\calC)$ réalise la localisation de
la catégorie $\calC$ par rapport aux morphismes $f$ tels que $T(f)$ est 
un isomorphisme.
\end{prop}

\begin{proof}
Supposons donné un foncteur $F$ de $\calC$ dans une catégorie $\calX$
tel que $F(f)$ est un isomorphisme dès que $T(f)$ en est un. Soit $G$ la
composée $\Max(\calC) \to \calC \to \calX$ où le premier foncteur est
l'inclusion canonique et le second est $F$. Pour tout $C \in \Ob \calC$,
$T(\iota_{\max}^{C}) = \id_{T(C)}$ est inversible, et donc par hypothèse
il en est de même de $F(\iota_{\max}^{C})$. La famille des
$F(\iota_{\max}^{C})$ définit donc une transformation naturelle
\emph{inversible} entre les foncteurs $F$ et $G \circ \Max$. Ceci
termine la preuve.
\end{proof}

En vertu de cette proposition, le foncteur $T : \calC \to \calA$ se 
factorise par $\Max(\calC)$ par l'intermédiaire d'un foncteur $T_\max : 
\Max(\calC) \to \calA$, qui n'est autre (d'après la preuve que l'on 
vient de donner) que la restriction de $T$ à $\Max(\calC)$.

\begin{prop}
\label{prop:conservatif}
Le foncteur $T_\max$ est conservatif, en ce sens qu'il vérifie : $f$ est 
un isomorphisme si, et seulement si $T_\max(f)$ en est un.

Si le foncteur $T$ est fidèle (resp. plein, resp. essentiellement 
surjectif), alors il en est de même de $T_\max$.
\end{prop}

\begin{proof}
La première partie de la proposition est une conséquence directe du
lemme \ref{lem:maxisom}. La seconde assertion à propos des propriétés de 
fidélite et de plénitude résulte de ce que $T_\max$ est obtenu comme une 
restriction du foncteur $T$. Pour la propriété d'essentielle 
surjectivité, elle résulte de l'égalité $T_\max \circ \Max = T$.
\end{proof}

\begin{prop}
\label{prop:fibtmax}
La fibration $T_\max : \Max(\calC) \to \calA$ vérifie encore les axiomes 
\Ax 0, \Ax 2 et \Ax 3 (où, bien entendu, les fibres $\calF_A$ sont 
calculées à partir du foncteur $T_\max$).

En outre, pour tout $A$, il existe un groupe $G_A$ tel que la catégorie 
$F_A$ soit équivalente à la catégorie ayant un unique objet $\bullet$ 
vérifiant $\End(\bullet) = G_A$. Si de plus \Ax 1 est vérifié, tous les 
groupes $G_A$ sont réduits à l'identité.
\end{prop}

\begin{proof}
Le seul point non trivial réside dans la vérification de \Ax 2. Mais, si 
$M \to M'$ et $M \to M''$ sont des morphismes dans la catégorie 
$\Max(\calC)$, on vérifie en utilisant la proposition 
\ref{prop:adjonction} que $\Max(M' \oplus_M M'')$ (où $M' \oplus_M M''$ 
désigne la somme amalgamée dans $\calC$) satisfait la propriété 
universelle de la somme amalgamée dans $\Max(\calC)$.
\end{proof}

\noindent
{\it Remarque.} Comme $T_\max$ vérifie encore les axiomes \Ax 0, \Ax 2
et \Ax 3, on peut répéter la construction $\Max$ et obtenir ainsi un
foncteur $\Max : \Max(\calC) \to \Max(\calC)$. Il est facile de voir à 
partir de ce qui précède que celui-ci est (isomorphe à) l'identité ; en 
particulier, $\Max(\Max(\calC)) = \Max(\calC)$.

\subsection{Dualité}
\label{subsec:dualite}

Introduisons les axiomes duaux de \Ax 2 et \Ax 3 à savoir 
respectivement :
\begin{axiome}{2*}
Les catégories $\calC$ et $\calA$ admettent des produits fibrés, et 
le foncteur $T$ y est compatible.
\end{axiome}
\begin{axiome}{3*}
Pour tout $A \in \Ob \calA$, la catégorie $\calF_A$ admet un objet initial.
\end{axiome}
Bien entendu, si ceux-ci sont satisfaits en plus de \Ax 0, on définit
par une construction analogue à la précédente un foncteur $\Min : \calC
\to \calC$ muni de morphismes naturels $\iota_{\min}^{C} : \Min(C) \to 
C$ pour tout $C \in \Ob \calC$. On dit qu'un objet est \emph{minimal} si
$\iota_{\min}^{C}$ est un isomorphisme, et on note $\Min(\calC)$ la
sous-catégorie pleine des objets minimaux. Toutes ces structures
vérifient évidemment des propriétes semblables à celles listées
précédemment pour le foncteur $\Max$ (que nous laissons au lecteur le 
soin d'écrire complètement). En particulier la fibration $T$ fournit
par restriction (ou, au choix, par passage au quotient) une fibration
$T_\min : \Min(\calC) \to \calA$.

Si $\calX$ est une catégorie, on définit une \emph{dualité} sur $\calX$ 
comme la donnée d'un foncteur contravariant $\calX \to \calX$, $X 
\mapsto X^\star$ et d'une identification fonctorielle entre 
$(X^\star)^\star$ et $X$. Considérons l'axiome suivant :
\begin{axiome}{4}
Il existe des dualités sur $\calC$ et sur $\calA$ compatibles au 
foncteur $T$ (c'est-à-dire qu'il existe une identification naturelle
entre $T(C^\star)$ et $T(C)^\star$).
\end{axiome}
S'il est vérifié, la dualité sur $\calC$ induit pour tout $A \in \Ob 
\calA$ une anti-équivalence de catégories entre $F_A$ et $F_{A^\star}$. 
On en déduit que, sous \Ax 4, les conditions \Ax 2 et \Ax {2*} (resp. 
\Ax 3 et \Ax {3*}) sont équivalentes.

\medskip

On suppose désormais \Ax 0, \Ax 2, \Ax 3, \Ax {2*} et \Ax{3*}. On
souhaite comparer les deux foncteurs $\Min : \calC \to \calC$ et $\Max 
: \calC \to \calC$, ainsi que les catégories $\Min(\calC)$ et 
$\Max(\calC)$ associées. On commence pour cela par un lemme.

\begin{lemme}
\label{lem:minomax}
On a $\Min \circ \Max = \Min$ et $\Max \circ \Min = \Max$ (sur les
objets et sur les morphismes).
\end{lemme}

\begin{proof}
Pour les objets, c'est immédiat au vu des définitions. Pour les 
morphismes, c'est une consé\-quence du lemme \ref{lem:maxmorph}.
\end{proof}

\begin{cor}
\label{cor:minmaxequiv}
Les restrictions $\Min : \Max(\calC) \to \Min(\calC)$ et $\Max : 
\Min(\calC) \to \Max(\calC)$ définissent des équivalences de catégories
inverses l'une de l'autre.
\end{cor}

\begin{proof}
D'après le lemme \ref{lem:minomax}, il suffit de montrer que le foncteur
$\Min$ (resp. $\Max$) est isomorphe au foncteur identité sur la catégorie
$\Min(\calC)$ (resp. $\Max(\calC)$), ce qui est immédiat par définition de
cette catégorie.
\end{proof}

\noindent
{\it Remarque.} Puisque les deux catégories $\Min(\calC)$ et 
$\Max(\calC)$ s'obtiennent comme localisation de $\calC$ par rapport
au même ensemble de morphismes (proposition \ref{prop:localisation}), on 
savait déjà qu'elles étaient équivalentes. Le corollaire précédent 
précise cela en donnant les foncteurs réalisant cette équivalence.

\begin{cor}
\label{cor:tmaxmin}
Les fibrations $T_\max : \Max(\calC) \to \calA$ et $T_\min : \Min(\calC) 
\to \calA$ satisfont toutes les deux les axiomes \Ax 0, \Ax 2, \Ax 3, 
\Ax {2*} et \Ax{3*}.
\end{cor}

\begin{proof}
Par la proposition \ref{prop:fibtmax}, $\Max(\calC)$ satisfait déjà \Ax
0, \Ax 2 et \Ax 3. En dualisant cette proposition, on obtient que
$\Min(\calC)$ satisfait \Ax {2*} et \Ax {3*}. Maintenant, les foncteurs
$\Min$ et $\Max$ commutant à $T$, le corollaire
\ref{cor:minmaxequiv} entraîne que les fibrations $T_\max$ et $T_\min$
sont isomorphes. La conclusion s'ensuit.
\end{proof}

\noindent
{\it Remarque.} Les sommes amalgamées et produits fibrés dans
$\Min(\calC)$ (resp. $\Max(\calC)$) s'obtiennent en appliquant le
foncteur $\Min$ (resp. $\Max$) aux constructions correspondantes dans
la catégorie $\calC$.

\bigskip

À partir de maintenant, on suppose en plus \Ax 4.

\begin{prop}
\label{prop:permute}
Pour $C \in \Ob \calC$, on a $\Max(C^\star) \simeq \Min(C)^\star$ et
$\Min(C^\star) \simeq \Max(C)^\star$. En particulier, la dualité de 
$\calC$ permute les catégories $\Min(\calC)$ et $\Max(\calC)$.

Le foncteur $C \mapsto \Max(C^\star)$ (resp. $C \mapsto \Min(C^\star)$)
définit une dualité de $\Max(\calC)$ (resp. $\Min(\calC)$) compatible au 
foncteur $T_\max$ (resp. $T_\min$).
\end{prop}

\begin{proof}
La première partie de la proposition résulte de ce que la dualité de 
$\calC$ induit une anti-équivalence de catégories entre $F_{T(C)}$ et 
$F_{T(C^\star)}$.

Si $D(C) = \Max(C^\star)$, on a, pour $C \in \Max(\calC)$ :
$$D(D(C)) = \Max(\Max(C^\star)^\star) \simeq \Max(\Min(C^{\star \star}))
\simeq \Max(\Min(C)) = \Max(C) \simeq C$$
dans l'ordre d'après la première partie de la proposition, la 
définition d'une dualité, le lemme \ref{lem:minomax}, et finalement le 
fait que $C$ soit maximal. Ce calcul assure que $D$ est une dualité.
La compatibilité à $T_\max$ est immédiate. Finalement, le même argument
fonctionne pour $C \mapsto \Min(C^\star)$.
\end{proof}

\subsection{Catégories fibrées en (semi-)treillis, pylonets}
\label{subsec:pylonets}

Dans les applications que l'on a en vue, on ne vérifiera jamais \Ax 3 
directement, mais on empruntera un chemin légèrement détourné que l'on
explique ci-dessous. Tout au long de ce paragraphe, on suppose \Ax 1.

\begin{lemme}
Soient $C$ et $C'$ deux objets d'une fibre $F_A$. Alors $\hom_{F_A} (C, 
C')$ a au plus un élément.
\end{lemme}

\begin{proof}
Tout $f \in \hom_{F_A} (C, C')$ vérifie par définition $T(f) = \id_A$. 
Le lemme résulte alors de la fidélité de $T$.
\end{proof}

On rappelle qu'une catégorie vérifiant la condition du lemme correspond 
simplement à un préordre sur l'\og ensemble \fg\ de ses objets : un
objet $C$ est plus petit que $C'$ s'il existe effectivement un morphisme
de $C$ dans $C'$. On rappelle également que les constructions usuelles
sur les ensembles (pré)ordonnés ont en général des équivalents simples 
en langage des catégories : par exemple, pour ne citer que celles qui 
vont nous intéresser dans la suite, une borne supérieure est une somme 
directe, et un élément maximal est un objet final\footnote{Ce que 
justifie la notation $\Max$ pour le foncteur construit précédemment.}.
Sachant cela, on démontre facilement (supposant toujours \Ax 1) que \Ax 
3 est impliqué par les deux axiomes suivants :
\begin{axiome}{3a}
Les catégories $\calF_A$ admettent des sommes directes (finies).
\end{axiome}
\begin{axiome}{3b}
Les catégories $\calF_A$ satisfont la condition de chaîne croissante 
(c.c.c) : pour tout suite infinie de morphismes
$$\xymatrix {
C_1 \ar[r]^-{f_1} & C_2 \ar[r]^-{f_2} & C_3 \ar[r]^-{f_3} & \cdots
\ar[r]^-{f_{n-1}} & C_n \ar[r]^-{f_n} & \cdots } $$
il existe un entier $N$ tel que $f_n$ soit un isomorphisme pour tout $n 
\geq N$.
\end{axiome}
Bien entendu, il existe des versions duales de ces axiomes à savoir :
\begin{axiome}{3a*}
Les catégories $\calF_A$ admettent des produits (finis).
\end{axiome}
\begin{axiome}{3b*}
Les catégories $\calF_A$ satisfont la condition de chaîne décroissante 
(c.c.d) : pour tout suite infinie de morphismes
$$\xymatrix {
C_1 & \ar[l]_-{f_1} C_2 & \ar[l]_-{f_2} C_3 & \ar[l]_-{f_3} \cdots &
\ar[l]_-{f_{n-1}} C_n & \ar[l]_-{f_n} \cdots } $$
il existe un entier $N$ tel que $f_n$ soit un isomorphisme pour tout $n 
\geq N$.
\end{axiome}
Sous l'hypothèse \Ax 1, ils impliquent \Ax {3*}. Par ailleurs, si l'on
suppose \Ax 4, les énoncés \Ax{3a} et \Ax{3a*} d'une part, et \Ax{3b} et 
\Ax{3b*} d'autre part sont équivalents.

\bigskip

Terminons ce paragraphe par quelques remarques et un peu de
terminologie. En théorie des ordres, un ensemble ordonné satisfaisant
les (équivalents des) axiomes \Ax {3a} et \Ax {3a*} est ce que l'on
appelle un \emph{treillis}. De même, les conditions qui apparaissent
dans \Ax{3b} et \Ax{3b*} sont ainsi nommées car les propriétés
correspondantes sur les ensembles ordonnées portent ces noms.
Tout ceci conduit à poser la définition suivante.

\begin{deftn}
\label{def:catfib}
Une fibration $T : \calC \to \calA$ vérifiant les axiomes \Ax 1, \Ax 2,
\Ax{3a} (resp. \Ax{3a*}) est appelée une \emph{catégorie fibrée en
sup-semi-treillis} (resp. une \emph{catégorie fibrée en
inf-semi-treillis}). Si les deux axiomes \Ax{3a} et \Ax{3a*} sont
vérifiés, on parlera simplement de \emph{catégorie fibrée en treillis}.

On dit que $T$ vérifie c.c.c (resp. c.c.d) si l'axiome \Ax{3b} (resp.
\Ax{3b*}) est satisfait.
On dit qu'elle est \emph{autoduale} si l'axiome \Ax 4 est satisfait.
\end{deftn}

\noindent
{\it Remarque.} 
On pourra s'étonner de ne pas voir apparaître \Ax 0 dans la définition 
précédente, alors que toute la théorie que nous avons développée semble
reposer sur cet axiome. Toutefois, comme nous l'avons expliqué en 
\ref{subsec:fibres}, on peut toujours remplacer $\calC$ par une 
catégorie équivalente pour laquelle \Ax 0 est satisfait. Nous préférons 
ne pas inclure \Ax 0 dans la définition précédente, car il ne sera en 
fait pas vérifié dans les exemples que nous allons manipuler par la 
suite.

\bigskip

Pour simplifier la terminologie, nous introduisons la définition 
suivante.

\begin{deftn}
\label{def:pylonet}
Une catégorie fibrée en sup-semi-treillis satisfaisant c.c.c est 
appelée un \emph{pylonet}.
\end{deftn}

\begin{wrapfigure}{r}{4.2cm}
\vspace{-1em}
\raggedleft
\includegraphics[width=4cm]{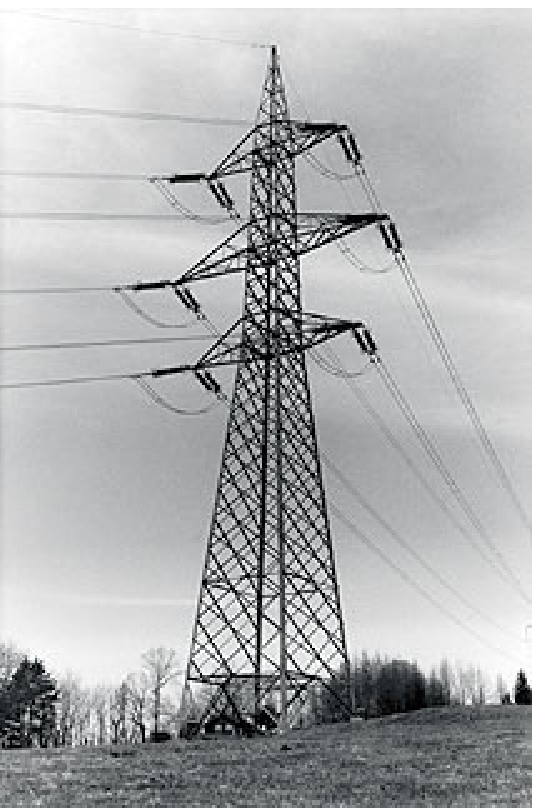}
\end{wrapfigure}

\noindent
{\it Remarques.} Cette terminologie est basée sur la concaténation des
deux mots \emph{pylone} et \emph{net}. Le premier d'entre eux se
rapporte aux fibres de $T$ qui, en un sens imagée, ressemblent à des
pylones électriques (voir photo ci-contre), la structure métallique de
ceux-ci pouvant évoquer un ordre admettant des bornes supérieures finies
et satisfaisant c.c.c (voire la condition plus forte \Ax{3c} donnée
plus bas). Le mot \emph{net}, quant à lui, est un anglicisme à prendre
dans le sens de réseau : il faut imaginer que ces pylones sont rélies
par tout un tissu de câbles (électriques) qui correspondent aux
morphismes de la catégorie $\calC$ dont l'image par $T$ n'est pas
l'identité. La propriété fondamentale des pylonets est que tout pylone
(\emph{i.e.} toute fibre) admet un sommet (\emph{i.e.} un élément
maximal) et qu'à tout cable reliant deux pylones (\emph{i.e.} tout
morphisme de $\calC$), il est associé un unique cable reliant les
sommets des pylones correspondants.

Si le foncteur $T$ est contravariant, on dira que le pylonet est
lui-même \emph{contravariant}.

\bigskip

Finalement, il est possible d'imaginer une version forte des axiomes 
\Ax{3b} et \Ax{3b*} qui est :
\begin{axiome}{3c}
Les catégories $\calF_A$ sont de \emph{hauteur finie}, dans le sens où il
existe un entier $N$ (qui dépend de $A$) telle que toute suite de $N$ 
morphismes
$$\xymatrix {
C_1 \ar[r]^-{f_1} & C_2 \ar[r]^-{f_2} & C_3 \ar[r]^-{f_3} & \cdots
\ar[r]^-{f_{N}} & C_{N+1} } $$
contient au moins un isomorphisme.
\end{axiome}
Contrairement à \Ax{3b} et \Ax{3b*}, l'axiome \Ax{3c} est autodual, et 
comme nous l'avons dit (ou du moins sous-entendu) précédemment, il 
implique à lui seul les deux énoncés \Ax{3b} et \Ax{3b*}. Encore une
fois, signalons que la terminologie \og de hauteur finie \fg\ est 
recopiée de celle couramment utilisée pour les treillis.

\subsection{Le cas additif}

Nous étudions à présent le cas particulier décrit par l'axiome suivant.
\begin{axiome}{5}
La catégorie $\calC$ est additive, la catégorie $\calA$ est abélienne
et le foncteur $T$ est additif.
\end{axiome}
On dira dans ce cas que la fibration $T$ est \emph{additive}. En
particulier, on pourra parler de \emph{catégories fibrées en 
(semi-)treillis additives}, et même de \emph{pylonets additifs}.

Remarquons que, sous \Ax 5, la condition \Ax 2 (resp. \Ax{2*}) est
équivalente à l'existence de conoyaux (resp. de noyaux) dans $\calC$ et
au fait que $T$ commute à la formation de ceux-ci. Supposons à partir de
maintenant, en plus de \Ax 5, les axiomes \Ax 0, \Ax 2, \Ax 3. D'après la
discussion menée en \ref{subsec:max}, on dispose d'un foncteur 
$\Max : \calC \to \calC$ et d'une sous-catégorie pleine $\Max(\calC)$ 
vérifiant un certain nombre de propriétés sympathiques.

\begin{lemme}
\label{lem:maxadditif}
Le foncteur $\Max$ est additif.
\end{lemme}

\begin{proof}
Il suffit de montrer que $\Max(C \oplus C')$ est naturellement isomorphe
à $\Max(C) \oplus \Max(C')$. Or, les inclusions canoniques $C \to C \oplus
C'$ et $C' \to C \oplus C'$ permettent de construire un morphisme
$\alpha : \Max(C) \oplus \Max(C') \to \Max(C \oplus C')$, tandis que les
projections $C \oplus C' \to C$ et $C \oplus C' \to C'$ donnent un
morphisme $\beta : \Max(C \oplus C') \to \Max(C) \oplus \Max(C')$. Il
est formel de vérifier que $\beta \circ \alpha$ est l'identité. Par
ailleurs, du fait que $T$ est additif, on déduit que $T(\alpha \circ 
\beta) = \id_{T(C \oplus C')}$. Ainsi $\alpha \circ \beta$ est un
endomorphisme de l'objet final de $F_{T(C \oplus C')}$ ; il ne peut
donc être que l'identité et le lemme en découle.
\end{proof}

Supposons maintenant en supplément de ce qui précède les axiomes duaux
\Ax{2*} et \Ax{3*}. Ils permettent à leur tour de construire un foncteur
additif $\Min : \calC \to \calC$ et une sous-catégorie $\Min(\calC)$.

\begin{prop}
\label{prop:abelien}
Dans la situation précédente,
les catégories $\Max(\calC)$ et $\Min(\calC)$ sont abéliennes et la
restriction du foncteur $T$ à ces catégories est exact.
\end{prop}

\begin{proof}
Nous ne donnons la preuve que pour $\Max(\calC)$, le cas de
$\Min(\calC)$ se traitant pareillement. Par le lemme
\ref{lem:maxadditif}, $\Max(\calC)$ contient l'objet nul et est stable
par somme directe ; c'est donc déjà une catégorie additive. Sachant
cela, le corollaire \ref{cor:tmaxmin} entraîne l'existence de noyaux et
de conoyaux dans $\Max(\calC)$. Pour conclure, il suffit de montrer que
si $f$ est une flèche dans $\Max(\calC)$, le morphisme induit $\bar f :
\coim f \to \im f$ est un isomorphisme. Or, comme $T$ commute à la
formation des noyaux et conoyaux dans $\Max(\calC)$ (corollaire
\ref{cor:tmaxmin}), $T(\bar f)$ s'identifie au morphisme $\coim T(f) \to
\im T(f)$ induit par $T(f)$. Comme $\calA$ est une catégorie abélienne,
$T(\bar f)$ est un isomorphisme, et donc, par le lemme
\ref{lem:maxisom}, $\Max(\bar f)$ également. Finalement, étant donné que  
par construction $\bar f$ est un morphisme entre deux objets maximaux, 
il s'identifie à $\Max(\bar f)$ et est par suite lui aussi un 
isomorphisme.

Il reste à montrer que la restriction de $T$ à $\Max(\calC)$ est exact, 
mais ceci découle directement de la commutation de ce foncteur à la
formation des noyaux et des conoyaux.
\end{proof}

Terminons ce paragraphe en soulignant qu'il est possible d'obtenir un
substitut à la proposition précédente dans une situation légèrement
différente. Précisément, on remplace \Ax{3*} par la nouvelle hypothèse
\Ax 1. Ce cas paraît de prime abord un peu batard car il ne permet pas
de définir le foncteur $\Min$. Malgré tout, on dispose de la proposition
suivante :

\begin{prop}
Dans la situation précédente, la catégorie $\Max(\calC)$ est abélienne
et la restriction du foncteur $T$ à cette catégorie est exacte.
\end{prop}

\begin{proof}
Prouvons tout d'abord que $\Max$ est un foncteur fidèle. Soit $f$ un 
morphisme de $\calC$ tel que $\Max(f) = 0$. En appliquant $T$ à cette
dernière égalité, on obtient $T(f) = 0$, puis $f = 0$ par fidélité de
$T$. Ceci démontre notre assertion.

Par le lemme \ref{lem:maxadditif}, $\Max(\calC)$ est une
catégorie additive. Par la première partie de la proposition
\ref{prop:fibtmax}, $\Max(\calC)$ admet des conoyaux et la formation
de ceux-ci commute au foncteur $T$. Il reste à montrer qu'il en est
de même pour les noyaux. En effet, après, on pourra appliquer le 
même raisonnement que dans la preuve de la proposition \ref{prop:abelien}
pour obtenir l'isomorphisme entre image et coimage.

Nous montrons en fait le résultat plus général suivant : si $f : C \to
C'$ est un morphisme dans $\calC$ et si $K$ est son noyau, alors
$\Max(K)$ est le noyau dans $\Max(\calC)$ du morphisme $\Max(f)$. Soit
$X$ un objet de $\Max(\calC)$ muni d'un morphisme $G : X \to \Max(C)$
tel que $\Max(f) \circ G = 0$. Nous voulons montrer que $G$ se factorise
de façon unique par $\Max(K)$. Notons $X' = X \times_{\Max(C)} C$ le
produit fibré de $X$ et $C$ au-dessus de $\Max(C)$. Quitte à remplacer
$X'$ par un objet isomorphe grâce à l'axiome \Ax 0, on a $T(X') = T(X)$, 
d'où il suit $\Max(X') \simeq X$. De plus, le morphisme canonique $g : X'
\to C$ vérifie $\Max(g) = G$, d'où on déduit $\Max(f \circ g) = 0$ puis
$f \circ g = 0$ par fidélité de $\Max$. Puisque $K$ est le noyau de $f$,
il suit que $g$ se factorise par $K$, et donc par fonctorialité, $G = 
\Max(g)$ se factorise par $\Max(K)$. L'unicité de cette factorisation
résulte à nouveau de la fidélité de $\Max$.
\end{proof}

\subsection{Avant-goût des applications}

Au fil des chapitres suivants, nous verrons que les pylonets sont des
structures qui apparaissent naturellement en géométrie algébrique, et
plus précisément en théorie de Hodge $p$-adique.  L'exemple fondamental
duquel tout ce travail est inspiré est le suivant : $\calC$ est la
catégorie des schémas en groupes commutatifs finis et plats annulés par
$p$ (ou une puissance de $p$) sur l'anneau des entiers d'un corps
$p$-adique $K$, alors que le foncteur $T$ est celui qui à un tel groupe
associe la représentation galoisienne donnée par ses $\bar K$-points. Il
résulte des travaux de Raynaud (voir \cite{raynaud}, \S 2.2) que cette
fibration $T$ est un pylonet autodual\footnote{Les dualités sont d'une
part la dualité de Cartier sur les schémas en groupes, et d'autre part
la dualité usuelle twistée (par le twist de Tate) sur les
représentations galoisiennes.} (définitions \ref{def:catfib} et
\ref{def:pylonet}) et additif (\emph{i.e.} satisfaisant \Ax 5). Le but
de cet article est de montrer que cette situation n'est pas isolée, mais
au contraire se généralise à de nombreuses autres fibrations rencontrées
en théorie de Hodge $p$-adique. Ce papier fait en réalité suite à un
travail antérieur \cite{carliu}, dans lequel il est montré (sans le dire
explicitement) que certaines catégories de modules définissent des
pylonets additifs et, dans certains cas favorables, autoduaux.

\medskip

Nous nous intéresserons donc par la suite à d'autres exemples de 
fibrations : $\calA$ sera la catégorie des $\F_p$-représentations (ou 
de $E$-représentations pour une extension finie $E$ de $\F_p$) du 
groupe de Galois absolu d'un corps $p$-adique, $\calC$ s'instanciera en 
certaines catégories de \og modules de Breuil \fg, et $T$ sera le 
foncteur de réalisation galoisienne correspondant. On rappelle (pour 
l'instant très brièvement) que $T$ admet généralement une version 
covariante et une contravariante. La version contravariante sera plus 
adaptée aux cas qui nous intéressent et c'est donc celle que nous 
manipulerons tout au long de cet article.

Nous allons montrer que ces données fournissent des pylonets 
(contravariants) autoduaux et additifs, c'est-à-dire, d'après les 
définitions, qu'elles obéissent à \Ax 1, \Ax 2, \Ax {3a}, \Ax {3b}, \Ax 
4 et \Ax 5.  De façon générale, la vérification de \Ax 5 sera toujours 
immédiate, alors que celle de \Ax 4, \Ax 1 et \Ax {3b} résultera 
directement de travaux antérieurs (\cite{carliu}). Ainsi, l'essentiel du 
travail consistera en l'établissement des énoncés \Ax 2 et \Ax {3a}. 
Après cela, on pourra déduire toute une liste de propriétés agréables 
sur la fibration $T$. Afin de faciliter la tâche du lecteur (et bien que 
cela fasse certainement redite), nous avons choisi de les regrouper dans 
le théorème suivant :

\begin{theo}
\label{theo:conseqpylonet}
Soit $T : \calC \to \calA$ un pylonet contravariant additif et autodual.
Alors :
\begin{itemize}
\item[$\bullet$] \emph{(cf \S \ref{subsec:pylonets})}
Pour tout $C \in \Ob \calC$, il existe un unique (à isomorphisme unique
près) couple $(\Max(C), \iota_{\max}^{C} : C \to \Max(C))$ (resp
$(\Min(C), \iota_{\min}^{C} : \Min(C) \to C)$) satisfaisant la propriété
universelle suivante :
\begin{itemize}
\item le $\calA$-morphisme $T(\iota_{\max}^{C})$ (resp.
$T(\iota_{\min}^{C})$) est un isomorphisme ;
\item pour tout $C' \in \Ob \calC$ muni d'une flèche $f : C \to C'$
(resp. $f : C' \to C$) telle que $T(f)$ est un isomorphisme, il existe
un unique $g : C' \to \Max(C)$ (resp. $g : \Min(C) \to C'$) tel que
$g \circ f = \iota_{\max}^{C}$ (resp. $f \circ g = \iota_{\min}^{ C}$).
\end{itemize}

\item[$\bullet$] \emph{(cf \S \ref{subsec:max})}
Ceci conduit à un foncteur \og idempotent \fg\ $\Max : \calC \to \calC$
(resp. $\Min : \calC \to \calC$).
\end{itemize}

\smallskip

\noindent
Si l'on note $\Max(\calC)$ (resp. $\Min(\calC)$) l'image essentielle de 
$\Max$ (resp. $\Min$), on a :
\begin{itemize}
\item[$\bullet$] \emph{(cf proposition \ref{prop:adjonction})}
La corestriction $\Max : \calC \to \Max(\calC)$ (resp. $\Min : \calC \to
\Min(\calC)$) est un adjoint à gauche (resp. à droite) du morphisme
d'inclusion.

\item[$\bullet$] \emph{(cf proposition \ref{prop:localisation} et
corollaire \ref{cor:minmaxequiv})}
Les foncteurs $\Max : \calC \to \Max(\calC)$ et $\Min : \calC \to
\Min(\calC)$ réalisent tous les deux la localisation de la catégorie
$\calC$ par rapport aux morphismes $f$ tels que $T(f)$ est un
isomorphisme. En particulier, les catégories $\Max(\calC)$ et
$\Min(\calC)$ sont équivalentes et, concrètement, cette équivalence se
réalise \emph{via} les foncteurs $\Min$ et $\Max$.

\item[$\bullet$] \emph{(cf proposition \ref{prop:conservatif})}
La restriction du foncteur $T$ à $\Max(\calC)$ d'une part, et à
$\Min(\calC)$ d'autre part est fidèle et conservative.

\item[$\bullet$] \emph{(cf proposition \ref{prop:permute})}
La dualité sur $\calC$ permute les catégories $\Max(\calC)$ et
$\Min(\calC)$. La composition de celle-ci avec le foncteur $\Max$ (resp.
$\Min$) induit une dualité sur $\Max(\calC)$ (resp. sur $\Min(\calC)$)
qui commute au foncteur $T$.

\item[$\bullet$] \emph{(cf proposition \ref{prop:abelien})}
La catégorie $\Max(\calC)$ (resp. $\Min(\calC)$) est abélienne, les
noyaux et conoyaux s'obtenant en appliquant le foncteur $\Max$ (resp.
$\Min$) aux constructions correspondantes dans $\calC$. La restriction
du foncteur $T$ à cette sous-catégorie est exacte.
\end{itemize}
\end{theo}

\section{Application à la théorie de Hodge $p$-adique}
\label{sec:application}

Nous reprenons à partir de maintenant les notations de l'introduction : 
$p$ est un nombre premier, $k$ un corps parfait de caractéristique $p$, 
$W$ l'anneau des vecteurs de Witt à coefficients dans $k$, $K_0$ son 
corps des fractions, $K$ une extension totalement ramifiée de $K_0$ de 
degré $e$. Fixons en outre $\bar K$ une clôture algébrique de $K$ et 
notons $G_K = \Gal(\bar K / K)$ le groupe de Galois absolu de $K$. 
Appelons $\O_K$ (resp. $\O_{\bar K}$) l'anneau des entiers de $K$ (resp. 
de $\bar K$). Soient également $\pi$ une uniformisante de $K$ et 
$(\pi_n)$ (resp. $(p_n)$) un système compatible de racines $p^n$-ièmes 
de $\pi$ (resp. de $p$). Soit $G_1 \subset G_K$ le groupe de Galois 
absolu de l'extension $K(\pi_1)$.

Dans tout le reste de cet article, on fixe un entier $r \in \{1, \ldots,
p-2\}$. Nous préférons éviter dès à présent le cas $r = 0$ car, bien
que fondamentalement plus simple, il conduit souvent à des discussions
assez peu intéressantes, et dans tous les cas, il vérifie certainement
l'inégalité $er < p-1$ et donc relève de l'étude menée dans
\cite{caruso-crelle}. Remarquons qu'ainsi, on peut supposer $p > 2$
(sinon aucun $r$ ne convient).

\subsection{Rappel sommaire de la théorie de Breuil}
\label{subsec:defbreuil}

On se borne dans cette sous-section à présenter les aspects \og annulés 
par $p$ \fg\ de la théorie de Breuil (développée de façon générale dans 
\cite{breuil-ens}, \cite{breuil-invent} et \cite{caruso-crelle}). 
Certains définitions (ou constantes) que nous allons introduire sont 
motivées par les aspects entiers de cette théorie (qui n'apparaitront 
que superficiellement dans cet article en \ref{subsec:reseaux}) et 
pourront de fait paraître étrange au lecteur qui n'est pas familier. 
Pour palier ce manque, nous renvoyons aux articles précédemment cités.

\subsubsection*{Les catégories de modules}

Posons $\tilde S = k[u]/u^{ep}$. Soit $c \in \tilde S$ la réduction
modulo $p$ du coefficient constant de polynôme minimal sur $K_0$ de
l'uniformisante $\pi$ fixée. On définit plusieurs catégories de
modules sur $\tilde S$. Tout d'abord, une grosse catégorie $\pFilpN
{\tilde S}$ (la notation deviendra claire en \ref{subsec:adjonctions})
dont les objets sont la donnée de :
\begin{enumerate}
\item un $\tilde S$-module $\calM$ ;
\item un sous-module $\Fil^r \calM \subset \calM$ contenant $u^{er}
\calM$ ;
\item un opérateur (dit de Frobenius) $\phi_r : \Fil^r \calM \to \calM$
semi-linéaire par rapport au Frobenius (c'est-à-dire l'élévation à la
puissance $p$) sur $\tilde S$ ;
\item un opérateur (dit de monodromie) $N : \calM \to \calM$ vérifiant :
\begin{itemize}
\item (condition de Leibniz) $N(ux) = u N(x) - ux $ pour tout $x \in
\calM$ ;
\item (transversalité de Griffith) $u^e N(\Fil^r \calM) \subset \Fil^r
\calM$ ;
\item le diagramme suivant est commutatif :
$$\xymatrix @C=50pt {
\Fil^r \calM \ar[r]^-{\phi_r} \ar[d]_-{u^e N} & \calM \ar[d]^-{c N} \\
\Fil^r \calM \ar[r]^-{\phi_r} & \calM }$$
\end{itemize}
\end{enumerate}
Les morphismes dans $\pFilpN {\tilde S}$ sont sans surprise les
applications $\tilde S$-linéaires qui commutent à toutes les structures
supplémentaires. Pour tout entier $t \leq r$, l'anneau $\tilde S$
lui-même muni de $\Fil^r \tilde S = u^{et} \tilde S$, de $\phi_r$
défini par $\phi_r(u^{et}) = c^t$ et de l'opérateur $N$ tel que $N(1) =
0$ est un exemple d'objet de $\pFilpN {\tilde S}$.
Avant de passer à la définition des autres catégories, signalons que
l'on dispose d'une notion de suite exacte dans $\pFilpN {\tilde S}$ : 
une suite d'objets de cette catégorie est dite exacte, si elle est 
exacte en tant que suite de $\tilde S$-modules, et si elle induit une 
suite exacte de $\tilde S$-modules au niveau des $\Fil^r$.

\medskip

Soit $\pModpN {\tilde S}$ la sous-catégorie pleine de $\pFilpN {\tilde
S}$ regroupant les objets $\calM$ pour lesquels $\phi_r(\Fil^r \calM)$
engendre $\calM$ comme $\tilde S$-module. La catégorie essentielle dont
nous voulons mener l'étude est encore une sous-catégorie de $\pModpN
{\tilde S}$ ; c'est celle qui regroupe les objets $\calM \in \pModpN
{\tilde S}$ qui sont des $\tilde S$-modules libres de type fini. On la
note $\ModpN {\tilde S}$ (sans l'apostrophe donc).

\subsubsection*{Foncteur vers Galois}

La catégorie $\pFilpN {\tilde S}$ est munie d'un foncteur $\Tst$ vers la 
catégorie $\Rep_{\F_p}(G_K)$ des $\F_p$-représenta\-tions du groupe 
$G_K$. Pour le définir, nous avons besoin d'introduire des anneaux de 
périodes : comme nous restons toujours dans le cas des représentations 
annulées par $p$, ces anneaux sont exceptionnellement faciles à décrire. 
Le premier d'entre eux est, en tant que $k$-algebre, $\hat A_0 = k 
\otimes_{(\phi),k} \O_{\bar K}/p$. Il est muni de l'action de $G_K$ 
naturelle. De plus, pour $0 \leq i \leq r$, on définit $\Fil^i \hat A_0$ 
comme l'idéal principal engendré par $1 \otimes p_1^i$ (où l'on rappelle 
que $p_1$ est une racine $p$-ième de $p$ fixée) : cela forme une 
filtration (finie) décroissante. On définit aussi pour les mêmes entiers 
$i$ une application $\phi_i : \Fil^i \hat A_0 \to \hat A_0$ en envoyant 
$1 \otimes p_1^i x$ sur la réduction modulo $p$ de $(-1)^i \otimes \hat 
x^p$ où $\hat x \in \O_{\bar K}$ est un relevé quelconque de $x$. (On 
montre que le résultat ne dépend que de $p_1^i x$, et pas de $x$ ni de 
son relevé $\hat x$.)

L'anneau qui nous intéressera le plus est $\hat A$, défini comme suit. 
En tant que $k$-algèbre, il vaut $\hat A_0 \brac X$ où la notation 
$\brac \cdot$ fait référence à l'algèbre polynomiale à puissance 
divisées. Il est muni d'un idéal $\Fil^r \hat A$ engendré par les 
produits $\Fil^{r-i} \hat A_0 \cdot X^i$ pour $0 \leq i \leq r$ et par 
les $\gamma_i(X) = \frac {X^i}{i!}$ pour $i > r$. On dispose également 
d'un morphisme $\phi_r : \Fil^r \hat A \to \hat A$ ; c'est celui qui 
envoie les éléments $\gamma_i(X)$ ($i > r$) sur $0$ et l'élément $a X^i$ 
($0 \leq i \leq r$, $a \in \Fil^{r-i} \hat A_0$) sur $\phi_{r-i}(a) Y^i$ 
avec $Y = \frac{(1+X)^p -1} p$, le calcul de cette dernière fraction se 
faisant bien entendu dans $\Z_p[X]$ avant d'être réduit dans $\hat A$. 
L'action de $G_K$ se prolonge à $\hat A$ grâce à la formule $g(X) = 
\frac{g(\pi_1)}{\pi_1} (1+X) - 1$. Enfin $\hat A$ apparaît comme une 
$\tilde S$-algèbre grâce au morphisme $\tilde S \to \hat A$, $u \mapsto 
\frac{\pi_1}{1+X}$. Tout cela fait de $\hat A$ un objet de $\pFilpN 
{\tilde S}$, et on peut poser pour tout objet $\calM \in \pFilpN {\tilde 
S}$ 
$$\Tst(\calM) = \hom_{\pFilpN S} (\calM, \hat A).$$ 
On définit comme ceci un foncteur contravariant $\Tst : \pFilpN {\tilde 
S} \to \Rep_{\F_p}(G_K)$.

\subsubsection*{Dualité}

La catégorie $\ModpN {\tilde S}$ est munie d'une dualité introduite dans
le chapitre V de \cite{caruso-thesis}. Rappelons que si $\calM$ en est
un objet, son dual $\calM^\star$ est défini comme suit :
\begin{enumerate}
\item $\calM^\star = \hom_{\tilde S\text{-mod}} (\calM, S)$ ;
\item $\Fil^r \calM^\star = \{ f \in \calM^\star \, / \, f(\Fil^r \calM)
\subset u^{er} S\}$ ;
\item pour $f \in \Fil^r \calM^\star$, $\phi_r^\star(f)$ est l'unique
application vérifiant $\phi_r^\star(f)(\phi_r(x)) = \phi_r \circ f (x)$
pour tout $x \in \Fil^r \calM$ où $\phi_r : u^{er} \tilde S \to \tilde
S$ est l'unique application semi-linéaire envoyant $u^{er}$ sur $c^r$ ;
\item pour $f \in \calM^\star$, $N^\star(f) = N \circ f - f \circ N$.
\end{enumerate}
L'association $\calM \mapsto \calM^\star$ définit une dualité dans le
sens du paragraphe \ref{subsec:dualite}. De plus, par le théorème
V.4.3.1 de \emph{loc. cit.} et la remarque qui le suit :
\begin{equation}
\label{eq:tstax4}
\Tst(\calM^\star) = \Tst(\calM)^\vee (r)
\end{equation}
où par définition \og $(r)$ \fg\ désigne le twist de Tate et où
$T^\vee$ est la représentation contragrédiente $\hom_{\F_p\text{-mod}} 
(T, \F_p)$. Autrement dit, si l'on muni la catégorie $\Rep_{\F_p}(G_K)$ 
de la dualité $T \mapsto T^\star = T^\vee(r)$ le foncteur $\Tst$ vérifie
l'axiome \Ax 4.

\subsubsection*{Sans l'opérateur de monodromie}

Il sera important dans la suite de considérer un analogue des objets
précédents dans lequel l'opérateur de monodromie est omis. Ceci amène à
définir tout d'abord la catégorie $\pFilp {\tilde S}$ dont les objets
sont les $\tilde S$-modules $\calM$ munis d'un $\Fil^r \calM$ et d'un
Frobenius $\phi_r : \Fil^r \calM \to \calM$ (mais \emph{pas} d'un
opérateur $N$) vérifiant les mêmes axiomes que précédemment. On isole
ensuite deux sous-catégories, à savoir $\pModp {\tilde S}$, $\Modp
{\tilde S}$, les définitions de celles-ci étant identiques à celles de
leurs analogues.

Le morphisme $\tilde S \to \hat A_0$, $u \mapsto \pi_1$ fait de $\hat
A_0$ une $\tilde S$-algèbre et permet de voir $\hat A_0$ comme un objet
de $\pFilp {\tilde S}$. On définit alors
$$\Tqst(\calM) = \hom_{\pFilp {\tilde S}} (\calM, \hat A_0)$$
pour $\calM \in \pModp {\tilde S}$. Il faut toutefois faire attention au
point suivant : le module $\Tqst(\calM)$ n'est pas une représentation de
$G_K$, mais seulement du sous-groupe $G_1$ étant donné que le morphisme
structural $\tilde S \to \hat A_0$ n'est pas $G_K$-équivariant (mais
seulement $G_1$-équivariant). On a malgré tout un lemme important qui
permet de comparer les foncteurs $\Tst$ et $\Tqst$.

\begin{lemme}
\label{lem:tsttqst}
Soit $\calM$ un objet de $\pModp {\tilde S}$. La projection
$\hat A_0$-linéaire $\hat A \to \hat A_0$, $\gamma_n(X) \mapsto 0$ ($n
\geq 1$) induit un isomorphisme $G_1$-équivariant $\Tst(\calM) \to
\Tqst(\calM)$.
\end{lemme}

\begin{proof}
La preuve est une version simplifiée de celle du lemme 2.3.1.1 de
\cite{breuil-invent} que l'on ne recopie pas. On notera par contre que 
celle-ci donne une formule explicite pour l'inverse $\Tqst(\calM) \to 
\Tst(\calM)$ : à $f_0 \in \Tqst(\calM)$, on associe l'application $f$ 
définie par
\begin{equation}
\label{eq:inverse}
f(x) = \sum_{i=0}^\infty f_0(N^i(x)) \gamma_i(\log(1+X))
\end{equation}
où la somme converge pour la topologie \og $\Fil$-adique \fg.
\end{proof}

Ces catégories \og sans $N$ \fg\ sont intéressantes car elles admettent
une description alternative plus simple. Soit $\tilde \si = k[[u]]$ que
l'on munit d'un opérateur de Frobenius $\phi : \tilde \si \to \tilde
\si$ qui agit comme l'élévation à la puissance $p$. Lorsque $\frakM$
est un module sur $\tilde \si$, on notera $\phi^\star \frakM = \tilde
\si \otimes_{(\phi),\tilde \si} \frakM$. Introduisons $\pModp \si$ la
catégorie dont les objets sont la donnée de :
\begin{enumerate}
\item un $\tilde \si$-module $\frakM$ ;
\item un opérateur $\phi$-semi-linéaire $\phi : \frakM \to \frakM$ tel
que le conoyau de $\id \otimes \phi : \phi^\star \frakM \to \frakM$
soit annulé par $u^{er}$.
\end{enumerate}
Comme dans les cas précédents, définissons $\Modp {\tilde \si}$ la
sous-catégorie pleine de $\pModp {\tilde \si}$ formée des objets libres
de type fini sur $\tilde \si$. On peut alors construire une équivalence
de catégories\footnote{Comme les définitions précises de tous les
objets qui interviennent ne nous seront pas vraiment utiles ici, nous ne
nous attardons par plus sur le sujet et nous contentons de renvoyer par
exemple à \cite{carliu} pour une présentation succinte de la théorie.} 
$M_{\tilde \si} : \Modp {\tilde \si} \to \Modp {\tilde S}$ qui jouït de 
propriétés intéressantes. En particulier, la composée $\Tsi = \Tqst 
\circ M_{\si_\infty}$ a une expression simple : pour tout $\frakM \in 
\Modp {\tilde \si}$, on a un isomorphisme $G_\infty$-équivariant, 
canonique et fonctoriel
\begin{equation}
\label{eq:Tsi}
\Tsi(\frakM) = \hom_{\tilde \si, \phi} (\frakM, k((u))^\sep)
\end{equation}
où $k((u))^\sep$ désigne une clôture séparable de $k((u))$ et est muni
du Frobenius usuel (l'élévation à la puissance $p$). Dans cette dernière
formule, l'action de $G_\infty$ sur le $\hom$ se fait par
l'intermédiaire d'une action sur $k((u))^\sep$ qui provient de la
théorie du corps des nombres. Il existe d'autres résultats concernant le 
foncteur $M_{\tilde \si}$ qui nous intéresserons particulièrement dans 
la suite. Nous les regroupons dans le théorème suivant.

\begin{theo}
\label{theo:sitilde}
Le foncteur $M_{\tilde \si}$ est une équivalence de catégories exacte.
Tout quasi-inverse est également exact et, de plus, respecte les
injections et les surjections.

Le foncteur $\Tsi : \Modp {\tilde \si} \to \Rep_{\F_p}(G_1)$ est un
pylonet additif et autodual et sa restriction à $\Max(\Modp {\tilde
\si})$ est pleinement fidèle.
\end{theo}

\begin{proof}
La première assertion est une généralisation directe (déjà utilisée
par ailleurs dans la littéra\-ture) d'un résultat de Breuil (théorème
4.1.1 de \cite{breuil-compo}). La phrase suivante concernant les
injections et les surjections est prouvée dans la proposition 2.3.2 de
\cite{carliu}.
Le second alinéa est, quant à lui, une version faible du résultat
principal\footnote{Le travail de \cite{carliu} n'est pas valable
seulement pour les objets annulés par $p$, pour une catégorie plus
grosse d'objets annulés par une puissance de $p$.} de \emph{loc. cit.},
même s'il n'est à aucun endroit écrit sous une forme aussi concise.
\end{proof}

\subsection{Opérateur de monodromie et prolongement de l'action de
Galois}
\label{subsec:tong}

Dans ce paragraphe, nous démontrons un résultat essentiel (proposition
\ref{prop:stabN}) qui précise les liens entre la donnée supplémentaire
d'un opérateur de monodromie $N$ et le prolongement de l'action de
Galois de $G_1$ à $G_K$. Les méthodes de démonstration (ainsi que
les énoncés d'ailleurs) sont très largement inspirées de celles
développées par Liu dans \cite{liu-compo}.

\medskip

On rappelle que $\pi_1 \in \O_K$ est une racine $p$-ième fixée de $\pi$.
Pour tout $\sigma \in G_K$, on définit $\epsilon(\sigma)$ comme l'image
dans $\hat A_0$ du quotient $\frac{\sigma \pi_1}{\pi_1}$ ; c'est une
racine $p$-ième de l'unité, qui vaut $1$ si $\sigma \in G_1$. On pose
en outre
$$t(\sigma) = \sum_{i=1}^{p-2} \frac{(1-\epsilon(\sigma))^i} i.$$
(On remarque que $(1-\epsilon(\sigma))^{p-1}$ s'annule dans $\hat A_0$,
ce qui est en accord avec le fait que l'on arrête la somme à $p-2$.) Il
est clair que si $\sigma \in G_1$, alors $t(\sigma) = 0$. Sinon,
$t(\sigma)$ est un élément de valuation $\frac 1 {p-1}$ et vérifie donc
en particulier $t(\sigma)^{p-1} = 0$. En outre, $t$ définit un cocycle,
\emph{i.e.} il est soumis à la relation $t(\sigma \sigma') = t(\sigma) +
\sigma t(\sigma')$, valable pour $\sigma$ et $\sigma'$ dans $G_K$.

Dans la suite, lorsque $\calM$ est un objet de $\Modp {\tilde S}$, on
sera amené à considérer le produit tensoriel $\calM \otimes_{\tilde S}
\hat A_0$ : il est naturellement muni d'un $\Fil^r$ (défini par $\Fil^r
\calM \otimes_{\tilde S} \hat A_0$), d'un $\phi_r$ (qui provient de
l'application $\phi_r : \Fil^r \calM \to \calM$) et d'une action de
$G_1$ (obtenue par son action naturellement sur le second facteur).
Lorsque de surcroît $\calM \in \ModpN {\tilde S}$, on prolonge l'action
de $G_1$ à $G_K$ tout entier en utilisant l'opérateur de monodromie
grâce à la formule
\begin{equation}
\label{eq:action}
\sigma(x \otimes a) = \sum_{i=0}^{p-2} N^i(x) \otimes \sigma(a) \: \frac
{t(\sigma)^i} {i!}
\end{equation}
avec $\sigma \in G_K$, $x \in \calM$ et $a \in \hat A_0$. En utilisant
$t(\sigma)^i t(\sigma')^j = 0$ pour $i + j \geq p-1$, on vérifie
aisément que l'égalité (\ref{eq:action}) définit bien une action. De
plus, on a la relation
\begin{equation}
\label{eq:tN}
N(x) \otimes t(\sigma) = \sum_{i=1}^{p-2} (-1)^{i-1} \frac {(\sigma -
1)^i} {i} (x \otimes 1).
\end{equation}
pour tout $\sigma \in G_K$.
Fixons $\tau$ un élément de $G_K$ qui n'appartient pas à $G_1$ ; avec
$G_1$, il engendre $G_K$ tout entier (puisque $G_1$ est d'indice premier
dans $G_K$). Notons $\epsilon = \epsilon(\tau)$ et $t = t(\tau)$ et
posons pour finir $\hat A_0^\star = \hat A_0$ que l'on munit de $\Fil^r
\hat A_0^\star = \hat A_0^\star$ et $\phi_r = \phi$. Comme dans la
démonstration du théorème 4.3.4 de \cite{liu-compo}, on construit (de
façon fonctorielle) des morphismes canoniques
$$\iota_\calM : \calM \otimes_{\tilde S} \hat A_0 \to \Tqst(\calM)^\star
\otimes_{\F_p} \hat A_0
\quad \text{et} \quad
\iota_\calM^\star : \Tqst(\calM)^\star \otimes_{\F_p} \hat A_0^\star \to
\calM \otimes_{\tilde S} \hat A_0$$
compatibles à $\Fil^r$, $\phi_r$ et l'action de $G_1$ (resp. $G_K$). En
outre, ils sont soumis à la relation $\iota_\calM \circ
\iota_\calM^\star = \id \otimes t^r$.

Toute cette artillerie permet de démontrer les deux propositions 
suivantes.

\begin{prop}
\label{prop:commN}
Soient $\calM$, $\calM'$ des objets de $\pModpN {\tilde S}$ et $f :
\calM \to \calM'$ un morphisme dans $\pModp {\tilde S}$. On suppose
que $\Tqst(f) : \Tst(\calM') \to \Tst(\calM)$ est $G_K$-équivariant.
Alors $f$ commute à $N$ (\emph{i.e.} $f$ est un morphisme dans
$\pModpN {\tilde S}$).
\end{prop}

\begin{proof}
La formule \eqref{eq:tN} implique que $f$ commute à $tN$ agissant 
sur les produits tensoriels $\calM \otimes_{\tilde S} \hat A_0$ et
$\calM' \otimes_{\tilde S} \hat A_0$. Dans $\hat A_0$, écrivons $t 
= q u^e$ où $q$ est un élément de valuation $\frac 1 {p(p-1)}$.
L'application $g = f \circ (u^e N) - (u^e N) \circ f$ prend alors
ses valeurs dans $\calM' \otimes_{\tilde S} m_q$ où $m_q$ est le
noyau de la multiplication par $q$ sur $\hat A_0$, c'est-à-dire l'idéal 
des éléments de valuation supérieure ou égale à $1 - \frac 1 {p(p-1)}$. 
Maintenant, pour $x \in \Fil^r \calM$, on a
$$\frac 1 c \: g \circ \phi_r(x) = \phi_r \circ (f \circ N - N \circ f) 
(x) \in \phi_r(\calM \otimes_{\tilde S} m_q) \subset \calM 
\otimes_{\tilde S} \phi_r(m_q) = 0$$
la dernière égalité provenant d'une simple calcul de valuation. On en
déduit, comme souhaité, que $f$ et $N$ commutent.
\end{proof}

\begin{prop}
\label{prop:stabN}
Soient $\calM \in \pModpN {\tilde S}$, $\calM' \in \pModp {\tilde S}$ et
$f : \calM \to \calM'$ un morphisme surjectif dans $\pModp {\tilde S}$.
On suppose que $\Tqst(\calM')$ (identifié grâce à $\Tqst(f)$ à une
sous-$G_1$-représentation de $\Tst(\calM)$) est stable par $G_K$. Alors,
il existe sur $\calM'$ un unique opérateur de monodromie pour lequel
$f$ est un morphisme dans $\pModpN {\tilde S}$.
\end{prop}

On commence par démontrer deux lemmes.

\begin{lemme}
\label{lem:keriota}
Pour tout $\calM \in \Modp {\tilde S}$, on a $\ker \iota_\calM
\subset t (\Fil^r \calM \otimes_{\tilde S} \hat A_0)$.
\end{lemme}

\begin{proof}
Notons $d$ le rang de $\calM$ sur $\tilde S$.
Posons $\calA = \Fil^r \calM \otimes_{\tilde S} \hat A_0$ et $\calB =
\iota_\calM(\calA)$. En s'appuyant sur le fait que $\hat A_0$ est un
anneau de Bézout et sur l'inclusion $\Fil^r \calM \subset \Fil^r S
{\cdot} \calM$, on montre que $\calA / t\calA$ est libre de rang $d$
sur $k \otimes_{(\phi),k} \O_{\bar K}/t$.
Par ailleurs, la compatibilité de $\iota_\calM^\star$ à $\Fil^r$ montre
que l'image de ce morphisme est incluse dans $\calA$. De la relation
$\iota_\calM \circ \iota_\calM^\star = \id \otimes t^r$, on déduit $t^r 
(\Tqst(\calM)^\star \otimes_{\F_p} \hat A_0) \subset \calB$, d'où il 
suit, comme précédemment, que $\calB / t\calB$ est aussi libre de rang 
$d$ sur $k \otimes_{(\phi),k} \O_{\bar K}/t$.
L'application $\iota_\calM$ induit une surjection linéaire $\calA /
t\calA \to \calB / t\calB$. Comme les espaces de départ et d'arrivée
sont des modules libres de même rang, c'est un isomorphisme et le lemme
en découle.
\end{proof}

\begin{lemme}
\label{lem:noyau}
Soit $f : \calM \to \calM'$ un morphisme surjectif dans $\Modp {\tilde
S}$ (resp. $\ModpN {\tilde S}$). Alors $\calK = \ker f$ (avec les
structures induites) est aussi dans $\Modp {\tilde S}$ (resp. $\ModpN
{\tilde S}$).
\end{lemme}

\begin{proof}
Il suffit de traiter le cas de $\Modp {\tilde S}$, l'opérateur de
monodromie ne posant pas de problèmes. Étant donné ce que nous avons vu,
le plus simple est de passer par l'équivalence avec $\Modp {\tilde
\si}$. D'après le théorème \ref{theo:sitilde}, $f$ provient d'un
morphisme surjectif $g : \frakM \to \frakM'$ de $\Modp {\tilde \si}$.
D'après la définition des objets de cette catégorie, il est clair que
$\ker g$ en est un objet. L'exactitude de $M_{\tilde \si}$ montre alors
que $M_{\tilde \si} (\ker g)$ s'identifie à $\calK$, d'où résulte le
lemme.
\end{proof}

\begin{proof}[Démonstration de la proposition \ref{prop:stabN}]
Soit $\calK$ le noyau de $f$ ; d'après le lemme \ref{lem:noyau}, c'est
un objet de $\Modp {\tilde S}$. Pour conclure, il suffit de montrer
qu'il est stable par $N$.
L'hypothèse assure que $\Tqst(\calK)^\star$ est stable par $G_K$ dans
$\Tst(\calM)^\star$. On en déduit, en utilisant l'exactitude de $\Tqst$
et l'égalité (\ref{eq:tN}), que $t f \circ N(\calK) \subset \ker
\iota_{\calM'}$. Avec le lemme \ref{lem:keriota}, on récupère $t f \circ
N(\calK) \subset t (\Fil^r \calM' \otimes_{\tilde S} \hat A_0)$. On 
suit alors la méthode de démonstration utilisée pour la proposition 
\ref{prop:commN} : dans $\hat A_0$, on peut écrire $t = q u^e$ où $q$ 
est un élément de valuation $\frac 1{p(p-1)}$. En \og divisant \fg\ la 
dernière inclusion par $q$, on obtient
$$f \circ (u^e N) (\calK) \subset (\calM' \otimes_{\tilde S} \ker m_q)
+ (u^e \Fil^r \calM' \otimes_{\tilde S} \hat A_0)$$
où $m_q$ désigne la multiplication par $q$ sur $\hat A_0$. On remarque
que $\ker m_q$ (resp. $\Fil^r S {\cdot} \hat A_0$) est formé des 
éléments de valuation supérieure ou égale à $1 - \frac 1 {p(p-1)}$ 
(resp. $\frac r p$). On en déduit $\ker m_q \subset \Fil^r S {\cdot} 
\hat A_0$, puis $\calM' \otimes_{\tilde S} \ker m_q \subset u^e \Fil^r 
\calM' \otimes_{\tilde S} \hat A_0$. Ainsi $f \circ (u^e N) (\calK) 
\subset u^e \Fil^r \calM' \otimes_{\tilde S} \hat A_0$.
Soit maintenant $x \in \Fil^r \calK$. Posons $y = \phi_r(x)$ et $z =
N(y)$. Par ce qui précède :
$$c f(z) = f \circ (cN) \circ \phi_r(x) = f \circ \phi_r \circ
(u^e N) (x) = \phi_r \circ f \circ (u^e N) (x) \in
\phi_r (u^e \Fil^r \calM' \otimes_{\tilde S} \hat A_0) = 0$$
d'où $f(z) = 0$, \emph{i.e.} $z \in \calK$ (car la flèche $\calK \to
\calK \otimes_{\tilde S} \hat A_0$ est injective). Puisque $\phi_r 
(\Fil^r
\calK)$ engendre $\calK$, on en déduit que $\calK$ est stable par $N$
comme voulu.
\end{proof}

\subsubsection*{Application : découpage par une sous-représentation}

Si $\calM$ est un objet de $\Modp {\tilde S}$ (resp. $\ModpN {\tilde 
S}$), tout quotient de $\calM$ (dans cette catégorie) détermine une 
sous-représentation de $T = \Tqst(\calM)$ (resp. $T = \Tst(\calM)$). 
Nous donnons ici une construction dans l'autre sens : à partir d'une 
sous-représentation de $T$, on retrouve un quotient (en fait, l'unique 
quotient) de $\calM$ qui lui correspond.

\begin{prop}
\label{prop:decoup}
Soient $\calM$ un objet de $\Modp {\tilde S}$ (resp. $\ModpN {\tilde
S}$) et $T'$ une sous-$G_1$-représentation de $T = \Tqst(\calM)$ (resp.
une sous-$G_K$-représentation de $\Tst(\calM)$). Alors, il existe un
unique quotient $\calM'$ de $\calM$ qui est un objet de $\Modp {\tilde 
S}$ (resp. $\ModpN {\tilde S}$) et pour lequel, en notant $f$ la 
projection canonique $\calM \to \calM'$, $\Tqst(f)$ (resp. $\Tst(f)$) 
s'identifie à l'inclusion $T' \toinj T$.
\end{prop}

\begin{proof}
On commence par traiter le cas des objets de $\Modp {\tilde S}$ 
(\emph{i.e.} sans monodromie). En utilisant l'équivalence avec $\Modp 
{\tilde \si}$, il revient au même de travailler dans cette dernière 
catégorie. Notons donc $\frakM$ l'objet de $\Modp {\tilde \si}$ 
correspondant à $\calM$. On rappelle que $G_\infty$ est le sous-groupe 
de $G_1$ correspondant à l'extension $K_\infty = \bigcup_{n \in \N} 
K(\pi_n)$. 
Nous allons utiliser la classification usuelle des représentations de 
$G_\infty$ à coefficients dans $\F_p$ telle que développée dans 
\cite{fontaine-fest}, \S A.1 (pour une présentation bien plus succinte, 
on pourra se reporter à \cite{carliu}, \S 3.1). Soit $M'$ le 
$\phi$-module sur $k((u))$ associé à $T'_{|G_\infty}$. La donnée de 
l'inclusion $T' \toinj T$ fait apparaître $M'$ comme un quotient de $M = 
\frakM[1/u]$. On note $\frakM'$ l'image de $\frakM$ dans $M'$. C'est un 
objet de $\Modp {\tilde \si}$ dont la $G_1$-représentation associée 
s'identifie à $T'$, au moins en tant que $G_\infty$-représentation. 
Toutefois, comme l'inclusion $T' \toinj T$ est par hypothèse 
$G_1$-équivariante, l'isomorphisme $\Tsi(\frakM') \simeq T'$ doit lui 
aussi être $G_1$-équivariant et l'existence est démontrée. L'unicité 
résulte de l'égalité $\ker f = \bigcup_{h \in T'} \ker h$, elle-même 
conséquence du lemme 2.1.5 de \cite{carliu}.

Le cas \og avec monodromie \fg\ s'obtient directement en combinant ce
que l'on vient de démontrer avec la proposition \ref{prop:stabN}.
\end{proof}

\begin{cor}
\label{cor:sousobquo}
Les images essentielles de $\Tqst$ défini sur $\Modp {\tilde S}$ et de 
$\Tst$ défini sur $\ModpN {\tilde S}$ sont stables par sous-objets et 
quotients.
\end{cor}

\begin{proof}
La stabilité par sous-objets est immédiate après la proposition
\ref{prop:decoup}. La stabilité par quotients s'obtient par dualité.
\end{proof}

Énonçons pour conclure ce paragraphe un corollaire de la proposition
\ref{prop:decoup} qui nous sera utile à plusieurs reprises dans la
suite.

\begin{cor}
\label{cor:decoup}
Soient $\calM$ un objet de $\Modp {\tilde S}$ (resp. $\ModpN {\tilde
S}$) et $T'$ une sous-$G_1$-représentation de $T = \Tqst(\calM)$ (resp.
une sous-$G_K$-représentation de $T = \Tst(\calM)$). On suppose $\big(
\bigcap_{h \in T'} \ker h \big) \subset u \calM$. Alors $T = T'$.
\end{cor}

\begin{proof}
La proposition \ref{prop:decoup} montre que $T' \toinj T$ provient d'un
morphisme surjectif $f : \calM \to \calM'$ dans $\Modp {\tilde S}$
(resp. $\ModpN {\tilde S}$), alors que le lemme \ref{lem:noyau} assure
que $\calK = \ker f$ est un objet de $\Modp {\tilde S}$ (resp. $\ModpN
{\tilde S}$). Par ailleurs, on a clairement $\calK \subset \bigcap_{h
\in T'} \ker h$, d'où $\calK \subset u \calM$. Par liberté de $\calK$ 
sur $\tilde S$, ceci ne peut se produire que si $\calK = 0$. Ainsi $f$ 
est un isomorphisme et $T = T'$.
\end{proof}

\subsection{Vérification des axiomes}

On est à présent en mesure de montrer que le foncteur $\Tst : \ModpN 
{\tilde S} \to \Rep_{\F_p}(G_K)$ satisfait certains axiomes de la 
section \ref{sec:axiomes}. Plus précisément, nous allons montrer que ce 
foncteur est un pylonet additif et autodual en établissant \Ax 1, \Ax 2, 
\Ax {3a}, \Ax {3b}, \Ax 4 et \Ax 5. En fait, \Ax 1 est déjà connu 
(corollaire 2.3.3 de \cite{carliu}), de même que \Ax 4 que nous avons 
rappelé brièvement en \ref{subsec:defbreuil}. L'axiome \Ax 5, quant à 
lui, est immédiat, tandis que \Ax{3b} résulte de la véracité de la 
proposition correspondant pour la fibration $\Tqst$. Il ne reste donc
qu'à prouver \Ax 2 et \Ax{3a}. C'est l'objet des deux propositions qui
suivent.

\begin{prop}
\label{prop:axiome2}
La fibration $\Tst : \ModpN {\tilde S} \to \Rep_{\F_p} (G_K)$ vérifie
l'axiome \Ax 2.
\end{prop}

\begin{proof}
Il faut prendre garde au fait que $\Tst$ est un foncteur contravariant.
On rappelle que notre convention à ce propos est de le considérer comme 
un foncteur covariant de $\ModpN {\tilde S}$ dans la catégorie opposée 
de $\Rep_{\F_p} (G_K)$. En particulier, \Ax 2 signifie que $\ModpN 
{\tilde S}$ admet des \emph{conoyaux} et que $\Tst$ transforme ceux-ci 
en noyaux. C'est ce que nous allons démontrer.

Soit $f : \calM \to \calM'$ un morphisme dans $\ModpN {\tilde S}$. On 
note $T = \Tst(\calM)$ et $T' = \Tst(\calM')$. Soit $\calC$ le quotient
de $\calM'$ associé à $K = \ker \Tst(f) \subset T'$ par la 
correspondance de la proposition \ref{prop:decoup}. Par définition, on a 
$\Tst(\calC) = \ker \Tst(f)$ et il suffit donc pour conclure de montrer 
que $\calC$ est un conoyau de $f$ dans $\ModpN {\tilde S}$. On considère 
pour cela $\calX \in \ModpN {\tilde S}$ muni d'un morphisme $g : \calM' 
\to \calX$ tel que $g \circ f = 0$. Notons $\pr : \calM' \to \calC$ la 
projection canonique. Soit $M_{\tilde S}$ un quasi-inverse de 
l'équivalence de catégories $M_{\tilde \si} : \Modp {\tilde \si} \to 
\Modp {\tilde S}$. \emph{Via} l'identification $T' = \Tsi \circ 
M_{\tilde S} (\calM')$ (voir formule \eqref{eq:Tsi}), on peut voir les 
éléments de $T'$ comme des morphismes de $M_{\tilde S} (\calM')$ dans 
$k((u))^\sep$, et c'est ce que nous ferons. Le lemme 2.1.5 de 
\cite{carliu} donne alors 
$$\ker M_{\tilde S}(\pr) = \bigcap_{h \in K} \ker h
\quad \text{et} \quad
\ker M_{\tilde S}(g) = \bigcap_{h \in L} \ker h$$
où $L$ est l'image de $\Tst(g) : \Tst(\calX) \to T'$. Du fait que $g 
\circ f = 0$, on déduit $L \subset K$ et donc, par les formules 
précédentes, que $\ker M_{\tilde S}(\pr) \subset \ker M_{\tilde S}(g)$. 
On en déduit que $M_{\tilde S}(g)$ se factorise par $M_{\tilde S}(\pr)$. 
En appliquant $M_{\tilde \si}$, on obtient un morphisme $g' : \calC \to 
\calX$ dans $\Modp {\tilde S}$ tel que $g' \circ \pr = g$. On remarque
alors que $\Tqst(g')$ n'est rien d'autre que la corestriction à $T'$ de 
$\Tqst(g) = \Tst(g)$. Ainsi $\Tqst(g)$ est $G_K$-équivariant et par la 
proposition \ref{prop:commN}, $g'$ commute à $N$, \emph{i.e.} $g'$ est 
un morphisme dans la catégorie $\ModpN {\tilde S}$. Pour montrer la
propriété universelle du conoyau, il ne reste plus qu'à justifier 
l'unicité de $g'$ mais elle est claire une fois que l'on a remarqué que 
$\pr$ est surjectif.
\end{proof}

\begin{prop}
\label{prop:axiome3a}
La fibration $\Tst : \ModpN {\tilde S} \to \Rep_{\F_p} (G_K)$ vérifie
l'axiome \Ax {3a}.
\end{prop}

\begin{proof}
Soient $T$ une $\F_p$-représentation de $G_K$, et $\calM_1$ et $\calM_2$ 
deux objets de $\ModpN {\tilde S}$ munis d'une identification 
$\Tst(\calM_1) \simeq \Tst(\calM_2) \simeq T$. On pose $\calM = \calM_1 
\oplus \calM_2$ et on définit $\calM'$ comme le quotient de $\calM$ 
attachée à la représentation diagonale de $\Tst(\calM) \simeq T \oplus 
T$ \emph{via} la correspondance de la proposition \ref{prop:decoup}. 
Montrons que $\calM'$ est la somme directe de $\calM_1$ et $\calM_2$ 
dans la fibre au-dessus de $T$.

Par construction, $\calM'$ est muni de morphismes $f_1 : \calM_1 \to
\calM'$ et $f_2 : \calM_2 \to \calM'$ (obtenus en plongeant d'abord
$\calM_1$ et $\calM_2$ dans $\calM$) qui induisent des isomorphismes
après application de $\Tst$. Pour conclure, il suffit de montrer que si
$\calN$ est un objet de $\ModpN {\tilde S}$ munis de morphismes $g_1 :
\calM_1 \to \calN$ et $g_2 : \calM_2 \to \calN$ induisant des
isomorphismes \emph{via} $\Tst$, alors il existe un unique morphisme $h
: \calN \to \calM'$ tel que $h \circ f = g$ où $f = f_1 \oplus (-f_2)$
et $g = g_1 \oplus (-g_2)$. Cela se fait de même que dans la preuve de 
la proposition \ref{prop:axiome2}.
\end{proof}

Pour récapituler, on a prouvé le théorème suivant :

\begin{theo}
\label{theo:tstpylonet}
La fibration $\Tst : \ModpN {\tilde S} \to \Rep_{\F_p} (G_K)$ est un
pylonet (contravariant) additif et autodual. En particulier, tous les
résultats du théorème \ref{theo:conseqpylonet} s'appliquent.
\end{theo}

\begin{deftn}
On pose $\MaxpN {\tilde S} = \Max(\ModpN {\tilde S})$ et $\MinpN {\tilde
S} = \Min(\ModpN {\tilde S})$.
\end{deftn}

Terminons par quelques remarques en revenant tout d'abord un instant sur 
le cas de $\Modp {\tilde S}$ (\emph{i.e.} sans $N$). \emph{Via} 
l'équivalence avec $\Modp {\tilde \si}$, les travaux de \cite{carliu} 
montrent que le foncteur $\Max$ s'interprète alors simplement comme 
l'extension des scalaires de $\tilde \si$ à $\tilde \si[1/u]$. Le 
théorème \ref{theo:tstpylonet} fournit donc en un certain sens un 
subsitut à cette extension des scalaires (qui n'est en effet par 
réalisable directement ici étant donné que $u$ est inversible). On peut 
se demander si emprunter ce chemin détourné est vraiment nécessaire, ou 
si au contraire, il n'existe pas une catégorie équivalente à $\ModpN 
{\tilde S}$ pour laquelle l'opération $\Max$ se réaliserait par un 
simple produit tensoriel. Le cas échéant, il serait intéressant de se 
demander en outre si ces nouveaux objets ont une interprétation 
cohomologique.

\medskip

Une des conséquences du théorème \ref{theo:tstpylonet} est le fait que 
la catégorie $\MaxpN {\tilde S}$ est abélienne. Un des intérets que cela 
peut présenter est l'utilisation des méthodes (co)homologiques en lien 
avec cette catégorie. Hélas, cela ne peut se faire directement car 
$\MaxpN {\tilde S}$ ne possède pas assez d'injectifs. Il s'agit par 
contre d'une catégorie dont tous les objets sont de longueur finie à 
laquelle on peut appliquer les méthodes de \cite{wildeshaus} : $\MaxpN 
{\tilde S}$ se plonge de façon pleinement fidèle dans la catégorie des 
ind-objets $\Ind(\MaxpN {\tilde S})$ dans laquelle on peut calculer les 
foncteurs dérivés de façon classique en utilisant des résolutions 
injectives.

\subsection{Un résultat de pleine fidélité}

\begin{theo}
\label{theo:pleinfid}
La restriction du foncteur $\Tst$ à $\MaxpN {\tilde S}$ est pleinement 
fidèle.
\end{theo}

\begin{proof}
Soient $\calM$ et $\calM'$ deux objets de $\MaxpN {\tilde S}$. On pose 
$T = \Tst(\calM)$, $T' = \Tst(\calM')$, et on suppose donné un morphisme 
$G_K$-équivariant $g : T' \to T$. En factorisant $g$ par $T' \tosurj \im 
g \toinj T$ et en se rappelant que l'image essentielle de $\Tst$ est 
stable par sous-objets (corollaire \ref{cor:sousobquo}), on se ramène à 
supposer successivement que $g$ est injectif puis surjectif.
Si $g$ est injectif, la proposition \ref{prop:decoup} montre l'existence
d'un morphisme $f : \calM \to \calM''$ dans $\ModpN {\tilde S}$ tel que
$\Tst(f) = g$. Le morphisme $\Max(f)$ relève alors $g$ dans la catégorie
$\MaxpN {\tilde S}$. La cas \og $g$ surjectif \fg\ s'obtient par 
dualité.
\end{proof}

\section{Quelques formules explicites}
La méthode que nous avons utilisée dans la section \ref{sec:application} 
pour démontrer le théorème \ref{theo:tstpylonet} a l'avantage d'être 
efficace mais, en contrepartie, elle donne une présentation des objets
construits (conoyaux, bornes supérieures dans une fibre) en termes de
représentations galoisiennes. En un sens, ceci n'est pas satisfaisant
car un des objectifs recherchés par l'introduction de la catégorie
$\ModpN {\tilde S}$ est de pouvoir faire des calculs entièrement du
côté \og algèbre linéaire \fg\ sans jamais avoir affaire aux 
représentations galoisiennes.

Cette section a pour but de remédier à ce problème. Pour cela, après 
avoir fait quelques développements sur le calcul des noyaux et conoyaux 
dans $\pModpN {\tilde S}$ en \ref{subsec:adjonctions}, nous construisons 
une nouvelle catégorie, notée $\RedpN {\tilde S}$, d'objets que nous 
qualifions de $\Tst$-réduits (sous-section \ref{subsec:red}). Nous 
montrons ensuite que cette catégorie est équivalente à $\ModpN {\tilde 
S}$ (sous-section \ref{subsec:equivs}) et nous explicitons enfin les 
constructions qui nous intéressent au niveau de $\RedpN {\tilde S}$ 
(sous-section \ref{subsec:explicite}).

Finalement, dans une dernière partie, nous donnons encore une formule 
explicite qui permet de retrouver à partir d'une représentation 
galoisienne $T$ appartenant à l'image essentielle de $\Tst$, l'objet 
maximal de $\RedpN {\tilde S}$ (ou $\ModpN {\tilde S}$) qui lui est 
associé.

\subsection{Deux adjonctions}
\label{subsec:adjonctions}

On commence par introduire de nouvelles catégories encore plus vastes 
que les précédentes. La première d'entre elle est $\pUnipN {\tilde S}$ 
($\text{Uni}$ pour \og univers \fg). Elle regroupe les objets qui sont 
la donnée des points suivants :
\begin{enumerate}
\item un $\tilde S$-module $\calM$ ;
\item un $\tilde S$-module $\Fil^r \calM$ muni d'un morphisme (pas
nécessairement injectif) $\iota : \Fil^r \calM \to \calM$ dont
l'image contient $u^{er} \calM$ ;
\item un morphisme $\phi$-semi-linéaire $\phi_r : \Fil^r \calM \to
\calM$ ;
\item des morphismes $N : \calM \to \calM$ et $N_{\Fil} : \Fil^r \calM
\to \Fil^r \calM$ vérifiant :
\begin{itemize}
\item (condition de Leibniz) $N(ux) = u N(x) - ux$ pour tout $x \in
\calM$ et $N_\Fil(ux) = u N_\Fil(x) - ux$ pour tout $x \in \Fil^r \calM$
\item les deux diagrammes suivant sont commutatifs :
\begin{equation}
\label{eq:diagNfil}
\raisebox{0.5\depth}{\xymatrix @C=50pt {
\Fil^r \calM \ar[r]^-{\iota} \ar[d]_-{N_\Fil} & \calM \ar[d]^-{u^e N} \\
\Fil^r \calM \ar[r]^-{\iota} & \calM }
\qquad \qquad
\xymatrix @C=50pt {
\Fil^r \calM \ar[r]^-{\phi_r} \ar[d]_-{N_\Fil} & \calM \ar[d]^-{c N} \\
\Fil^r \calM \ar[r]^-{\phi_r} & \calM } }
\end{equation}
\end{itemize}
\end{enumerate}
Les morphismes dans $\pUnipN {\tilde S}$ sont les paires $(f : \calM \to
\calM', f_\Fil : \Fil^r \calM \to \Fil^r \calM')$ qui sont compatibles à
toutes les structures additionnelles. On isole la sous-catégorie pleine 
$\pGenpN {\tilde S}$ de $\pUnipN {\tilde S}$ qui regroupe les objets 
$\calM$ tels que $\phi_r(\Fil^r \calM)$ engendre $\calM$ comme $\tilde 
S$-module. On dispose du diagramme suivant :
$$\xymatrix @R=7pt {
& \pFilpN {\tilde S} \ar@{^(->}[rd] \\
\pModpN {\tilde S} \ar@{^(->}[ru] \ar@{^(->}[rd] & & \pUnipN {\tilde
S} \\
& \pGenpN {\tilde S} \ar@{^(->}[ru] }$$
où les flèches $\toinj$ symbolisent des foncteurs pleinement fidèles.
En outre, l'image de $\pModpN {\tilde S}$ dans $\pGenpN {\tilde S}$
(resp. de $\pFilpN {\tilde S}$ dans $\pUnipN {\tilde S}$) est constituée
des objets pour lesquels le morphisme $\iota$ est injectif. Il est
finalement facile de voir que le parallélogramme précédent est
cartésien, c'est-à-dire que l'intersection (calculée dans $\pUnipN
{\tilde S}$) des catégories $\pFilpN {\tilde S}$ et $\pGenpN {\tilde S}$
n'est autre que $\pModpN {\tilde S}$.

Les notations $\Gen$ et $\Fil$ doivent maintenant être plus claires : on
utilise $\Gen$ (comme \og en\emph{gen}dre \fg\ ou \og \emph{gen}erate
\fg) pour désigner les objets sur lesquels l'image de $\phi_r$ engendre
tout, et $\Fil$ (comme \og \emph{fil}tration \fg) pour les objets pour
lesquels $\Fil^r \calM$ définit un véritable sous-module, c'est-à-dire
pour lesquels l'application $\Fil^r \calM \to \calM$ est injective.

\medskip

Le but de cette sous-section est de construire des adjoints aux quatre
foncteurs d'inclusion que nous venons d'introduire. Étant donné que nous
ne souhaitons pas nous limiter à une catégorie d'objets de type fini (en
particulier pour les constructions menées en \ref{subsec:pleinfid}), la
construction de ces adjonctions va reposer sur une induction transfinie.
Si le lecteur n'est pas familier avec ce type de manipulations, et qu'il
ne souhaite pas s'impliquer trop loin dans cette direction, nous
l'invitons à supposer que tous les objets $\calM$ sont de type fini sur
$\tilde S$, à remplacer systématiquement dans la suite le mot \og
ordinal \fg\ (resp. \og induction transfinie \fg) par la locution \og
entier naturel \fg\ (resp. \og récurrence \fg), et à ignorer tout ce qui
concerne les ordinaux limites. L'hypothèse de type finitude, combinée au
fait que $\tilde S$ soit un anneau artinien, entraîne que toutes les
constructions itératives que nous allons entreprendre se stabilisent au
bout d'un nombre fini (et pas transfini) d'étapes.

\subsubsection*{Le foncteur $\Gen$}

Soit $\calM$ un objet de $\pUnipN {\tilde S}$. On définit par induction
transfinie une suite décroissante $(\Gen_\alpha(\calM)$) (indexée par
les ordinaux $\alpha$) de sous-objets (dans $\pUnipN {\tilde S}$) de
$\calM$. Pour $\alpha = 0$, on pose simplement $\Gen_0(\calM) = \calM$.
Si $\alpha$ est un ordinal limite, on pose $\Gen_\alpha(\calM) =
\bigcap_{\beta < \alpha} \Gen_\beta(\calM)$ et $\Fil^r
\Gen_\alpha(\calM) = \bigcap_{\beta < \alpha} \Fil^r \Gen_\beta(\calM)$.
Finalement, si $\alpha = \beta + 1$ est un ordinal successeur,
$\Gen_\alpha(\calM)$ est le sous-$\tilde S$-module engendré par
$\phi_r(\Fil^r\Gen_\beta(\calM))$ (qui est bien un sous-module de
$\Gen_\beta(\calM)$ puisque, par construction, $\Gen_\beta(\calM)$ est
un objet de $\pUnipN {\tilde S}$). On le munit de $\Fil^r
\Gen_\alpha(\calM) = \iota^{-1} (\Gen_\alpha (\calM))$. Par construction
l'application $\phi_r$ envoie $\Fil^r \Gen_\beta(\calM)$ dans
$\Gen_\alpha(\calM)$ et donc induit bien par restriction un morphisme
$\phi_r : \Fil^r \Gen_\alpha(\calM) \to \Gen_\alpha(\calM)$. De même,
par définition de $\Fil^r \Gen_\alpha(\calM)$, le morphisme $\iota$
envoie $\Fil^r \Gen_\alpha(\calM)$ dans $\Gen_\alpha(\calM)$. Les
diagrammes (\ref{eq:diagNfil}) impliquent dans l'ordre que $N$ stabilise
$\Gen_\alpha(\calM)$ puis que $N_\Fil$ stabilise $\Fil^r
\Gen_\alpha(\calM)$. On a ainsi bien défini un objet
$\Gen_\alpha(\calM)$ de $\pUnipN {\tilde S}$, ce qui termine notre
induction transfinie.

Les propriétés usuelles des ordinaux impliquent que la suite 
$(\Gen_\alpha(\calM))$ est stationnaire. On appelle $\Gen(\calM)$ la 
valeur limite atteinte par cette suite. Il est alors clair que 
$\phi_r(\Fil^r \Gen(\calM))$ engendre $\Gen(\calM)$, c'est-à-dire que 
$\Gen(\calM)$ est un objet de $\pGenpN {\tilde S}$. Par induction 
transfinie, on montre qu'un morphisme $\calM \to \calM'$ induit par 
restriction des flèches $\Gen_\alpha(\calM) \to \Gen_\alpha (\calM')$ 
pour tout ordinal $\alpha$, et donc finalement un morphisme $\Gen(\calM) 
\to \Gen (\calM')$. On obtient ainsi des foncteurs $\Gen_\alpha : 
\pUnipN {\tilde S} \to \pUnipN {\tilde S}$ pour tout ordinal $\alpha$ et 
un foncteur limite $\Gen : \pUnipN {\tilde S} \to \pGenpN {\tilde S}$. 
Par ailleurs, il est facile de voir que la restriction de $\Gen$ à 
$\pFilpN {\tilde S}$ prend ses valeurs dans $\pModpN {\tilde S}$.

\begin{lemme}
\label{lem:genadjoint}
Le foncteur $\Gen : \pUnipN {\tilde S} \to \pGenpN {\tilde S}$ (resp.
$\Gen : \pFilpN {\tilde S} \to \pModpN {\tilde S}$) est un adjoint à
droite de l'inclusion canonique $\pGenpN {\tilde S} \toinj \pUnipN
{\tilde S}$ (resp. $\pModpN {\tilde S} \toinj \pFilpN {\tilde S}$).
\end{lemme}

\begin{proof}
On ne donne la preuve que pour les catégories $\pUnipN {\tilde S}$ et
$\pGenpN {\tilde S}$, l'autre cas étant absolument semblable. Soient
$\calM \in \pGenpN {\tilde S}$ et $\calM' \in \pUnipN {\tilde S}$. Il
suffit de montrer que tout morphisme $f : \calM \to \calM'$ se factorise
de façon unique par $\Gen(\calM')$. L'unicité résulte de ce que les
flèches $\Gen(\calM') \to \calM'$ et $\Fil^r \Gen(\calM') \to \Fil^r
\calM'$ sont injectives. Pour l'existence, il suffit de remarquer que
$\Gen(f)$ permet cette factorisation.
\end{proof}

\begin{cor}
La catégorie $\pModpN {\tilde S}$ admet des noyaux.
\end{cor}

\begin{proof}
Soit $f : \calM \to \calM'$ un morphisme dans $\pModpN {\tilde S}$. Le
noyau au sens usuel de $f$, disons $\calK$, hérite par restriction des
structures supplémentaires de $\calM$ : on pose $\Fil^r \calK = \calK
\cap \Fil^r \calM$, et on vérifie directement $\phi_r(\Fil^r \calK)
\subset \calK$ et $N(\calK) \subset \calK$. On obtient comme ci un
objet de $\pFilpN {\tilde S}$. Le lemme \ref{lem:genadjoint} assure
alors que $\Gen(\calK)$ est un noyau de $f$ dans la catégorie $\pModpN
{\tilde S}$.
\end{proof}

Il existe une version légèrement plus précise du corollaire précédent.
Elle dit que la catégorie $\pFilpN {\tilde S}$ admet des noyaux (ceux-ci
sont construits de la manière naïve) et que si $f : \calM \to \calM'$
est un morphisme de $\pFilpN {\tilde S}$ qui admet pour noyau $\calK$,
alors $\Gen(f)$ admet pour noyau $\Gen(\calK)$ dans la catégorie
$\pModpN {\tilde S}$.
On prendra garde par contre au fait que ceci n'implique aucune
exactitude (au sens des suites exactes dans $\pFilpN {\tilde S}$) pour
le foncteur $\Gen$. On a toutefois, à ce sujet, le résultat très partiel
suivant :

\begin{lemme}
\label{lem:gennul}
Soit $0 \to \calM'' \to \calM \to \calM' \to 0$ une suite exacte dans
$\pFilpN {\tilde S}$. On suppose qu'il existe un ordinal $\alpha$ tel 
que $\Gen_\alpha(\calM') = 0$. Alors pour tout ordinal $\beta$, on a
$\Gen_{\alpha + \beta} (\calM) \subset \Gen_\beta (\calM'')$.
\end{lemme}

\begin{proof}
Il suffit de prouver le lemme pour $\beta = 0$, les autres cas se
déduisant de celui-ci par une induction transfinie immédiate. Or
l'image de $\Gen_\alpha(\calM)$ dans $\calM'$ est contenue dans
$\Gen_\alpha(\calM')$ qui est nul par hypothèse. Ainsi
$\Gen_\alpha(\calM) \subset \calM'' = \Gen_0(\calM'')$, comme voulu.
\end{proof}

\begin{cor}
\label{cor:gennul}
Si $0 \to \calM'' \to \calM \to \calM' \to 0$ est une suite exacte et
si $\Gen(\calM') = \Gen(\calM'') = 0$, alors $\Gen(\calM) = 0$.
\end{cor}

\begin{proof}
L'hypothèse $\Gen(\calM') = \Gen(\calM'') = 0$ implique l'existence
d'ordinaux $\alpha$ et $\beta$ tels que $\Gen_\alpha(\calM') = 0$ et
$\Gen_\beta(\calM'') = 0$. Par le lemme précédent $\Gen_{\alpha + \beta}
(\calM) = 0$ et la corollaire en résulte.
\end{proof}

\subsubsection*{Le calcul de $\Gen(\hat A)$}

Si $\calM$ est un objet de $\pModpN {\tilde S}$, on a $\Gen(\calM) =
\calM$, ce qui entraîne $\Tst(\calM) = \hom_{\pModpN {\tilde S}} (\calM,
\Gen(\hat A))$. Il semble donc intéressant de calculer $\Gen(\hat A)$,
et c'est ce que nous nous proposons de faire ci-après comme premier
exercice de manipulation du foncteur $\Gen$. Cerise sur le gâteau, nous
allons constater qu'il a une structure très simple.

\medskip

Du fait que tout élément de $\O_{\bar K}$ possède certainement une
racine $p$-ième, on déduit que $\phi_i : \Fil^i (k \otimes_{(\phi),k}
\O_{\bar K} / p) \to k \otimes_{(\phi),k} \O_{\bar K} / p$ est
surjectif. De la description de l'action du Frobenius sur $\hat A$
donnée lors de la définition, il suit
$$\phi_r(\Fil^r \hat A) = \sum_{i=0}^r (k \otimes_{(\phi),k} \O_{\bar
K}/p) \cdot Y^i \subset \hat A.$$
Par définition $\Gen_1(\hat A)$ est le sous-$\tilde S$-module de $\hat
A$ engendré par $\phi_r(\Fil^r \hat A)$. Ainsi si l'on note $\hat S$
l'image du morphisme naturel $\tilde S \otimes_k \hat A_0 \to \hat A$,
on a $\Gen_1(\hat A) = \sum_{i=0}^r \hat S \cdot Y^i$.
Il s'agit maintenant de calculer les itérés suivants, mais par chance,
cela est assez simple. En effet, on remarque que l'élément $(1 \otimes
\pi_1^{r-i}) Y^i$ est simultanément dans $\Gen_1(\hat A)$ et dans
$\Fil^r \hat A$. De plus $\phi_r((1 \otimes \pi_1^{r-i}) Y^i) =
(-1)^{r-i} Y^i$. Ceci entraîne $\Gen_2(\hat A) \supset \Gen_1(\hat A)$
puis l'égalité, l'inclusion réciproque étant contenue dans la
définition. Ainsi la suite des $\Gen_\alpha(\hat A)$ est stationnaire à
partir de $\alpha = 1$. L'expression de $\Gen(\hat A)$ en résulte
directement
$$\Gen(\hat A) = \Gen_1(\hat A) = \sum_{i=0}^r \hat S \cdot Y^i.$$
Le lemme suivant termine de préciser la structure algébrique de
$\Gen(\hat A$).

\begin{lemme}
La famille des $Y^i$ ($0 \leq i \leq r$) est libre sur $\hat S$.
Ainsi, $\Gen(\hat A)$ est un $\hat S$-module libre de rang $r+1$
de base $(Y^i)_{0 \leq i \leq r}$.
\end{lemme}

\begin{proof}
Nous montrons un résultat légèrement plus fort, à savoir la liberté
sur un anneau plus gros $\hat B$. Cet anneau est défini comme la
sous-$(k \otimes_{(\phi),k} \O_{\bar K}/p)$-algèbre de $\hat A$
engendrée par $X$. Il est isomorphe à $(k \otimes_{(\phi),k} \O_{\bar
K}/p)[X]/X^p$ et fait de $\hat A$ un $\hat B$-module libre de base
$(\gamma_{pj}(X))_{j \geq 0}$. Par ailleurs, un calcul direct montre que
sur cette base $Y^i$ (pour $0 \leq i \leq r$) n'a de composantes non
nulles que pour $j \leq i$, et que la composante en $j = i$ est
inversible (c'est un élément non nul de $\F_p$). La conclusion en
découle.
\end{proof}

\noindent
{\it Remarque.} De ce qui précède, il résulte sans mal que l'image
de $\phi_r : \Fil^r \Gen(\hat A) \to \Gen(\hat A)$ est le sous
$(k \otimes_{(\phi),k} \O_{\bar K}/p)$-module (libre) engendré par les
$Y^i$ avec $0 \leq i \leq r$.

\subsubsection*{Le foncteur $\Fil$}

De façon assez semblable à ce qui vient d'être fait, on construit
maintenant un foncteur $\Fil : \pUnipN {\tilde S} \to \pFilpN {\tilde
S}$. Soit $\calM \in \pUnipN {\tilde S}$. On définit par induction
une suite transfinie $(\Fil_\alpha(\calM))$ de quotients successifs de
$\calM$. On pose tout d'abord $\Fil_0(\calM) = \calM$.

Supposons que $\alpha = \beta + 1$ soit un ordinal successeur. Notons
$K$ le noyau de $\iota : \Fil^r \Fil_\beta(\calM) \to \Fil_\beta(\calM)$ 
et
$Q$ le sous-$\tilde S$-module de $\Fil_\beta(\calM)$ engendré par
$\phi_r(K)$. Définissons :
$$\Fil_\alpha(\calM) = \frac{\Fil_\beta(\calM)} Q \quad \text{et} \quad
\Fil^r \Fil_\alpha(\calM) = \frac{\Fil^r \Fil_\beta(\calM)} K.$$
Par construction, $\phi_r : \Fil^r \Fil_\beta(\calM) \to 
\Fil_\beta(\calM)$
se factorise en un morphisme $\phi_r : \Fil^r \Fil_\alpha(\calM) \to
\Fil_\alpha(\calM)$. En outre, $\iota$ induit une application 
(injective)
$\Fil^r \Fil_\alpha(\calM) \to \Fil_\beta(\calM)$ que l'on peut
composer avec la projection canonique $\Fil_\beta(\calM) \to
\Fil_\alpha(\calM)$ pour obtenir un nouveau morphisme $\iota : \Fil^r
\Fil_\alpha(\calM) \to \Fil_\alpha(\calM)$. La commutation des
diagrammes (\ref{eq:diagNfil}) implique dans l'ordre $N_\Fil(K) \subset
K$ puis $N(Q) \subset Q$. On en déduit des opérateurs $N_\Fil$ et $N$
agissant respectivement sur les quotients $\Fil^r \Fil_\alpha(\calM)$ et
$\Fil_\alpha(\calM)$ dont il est facile de vérifier qu'ils font encore
commuter les diagrammes (\ref{eq:diagNfil}). Bref, on obtient comme cela
un objet $\Fil_\alpha(\calM) \in \pUnipN {\tilde S}$ muni d'un morphisme
surjectif $\Fil_\beta(\calM) \to \Fil_\alpha(\calM)$. Ceci montre que
$\Fil_\alpha(\calM)$ apparaît comme un quotient de $\Fil_\beta(\calM)$
et donc aussi de $\calM$.

Finalement, si $\alpha$ est un ordinal limite, on pose simplement
$\Fil_\alpha(\calM) = \varinjlim_{\beta < \alpha} \Fil_\beta(\calM)$.
La suite des $\Fil_\alpha(\calM)$ est stationnaire, et sa limite, notée
$\Fil(\calM)$, est nécessairement un objet de $\pUnipN {\tilde S}$ sur
lequel $\iota$ est injectif, c'est-à-dire un objet de $\pFilpN
{\tilde S}$. Par ailleurs, si $f : \calM \to \calM'$ est un morphisme
dans $\pUnipN {\tilde S}$, on vérifie par induction transfinie qu'il
induit pour tout ordinal $\alpha$ un morphisme $\Fil_\alpha(\calM) \to
\Fil_\alpha(\calM')$ et donc, par passage à la limite, une flèche
$\Fil(\calM) \to \Fil(\calM')$. Ainsi, obtient-on pour tout ordinal
$\alpha$ un foncteur $\Fil_\alpha : \pUnipN {\tilde S} \to \pUnipN
{\tilde S}$, ainsi qu'un foncteur $\Fil : \pUnipN {\tilde S} \to \pFilpN
{\tilde S}$. Le fait que $\Fil(\calM)$ apparaisse comme un quotient de
$\calM$ montre que $\Fil$ stabilise la catégorie $\pGenpN {\tilde S}$.
Ainsi, il induit par restriction un foncteur $\Fil : \pGenpN {\tilde S}
\to \pModpN {\tilde S}$.

\begin{lemme}
\label{lem:filadjoint}
Le foncteur $\Fil : \pUnipN {\tilde S} \to \pFilpN {\tilde S}$ (resp.
$\Fil : \pGenpN {\tilde S} \to \pModpN {\tilde S}$) est un
adjoint à gauche du foncteur d'inclusion $\pFilpN {\tilde S} \toinj
\pUnipN {\tilde S}$ (resp. $\pModpN {\tilde S} \toinj \pGenpN {\tilde
S}$).
\end{lemme}

\begin{proof}
Soient $\calM \in \pUnipN {\tilde S}$ et $\calM' \in \pFilpN {\tilde
S}$. Il suffit de montrer que tout morphisme $f : \calM \to \calM'$ se
factorise de façon unique par $\Fil(\calM)$. L'unicité résulte de ce que
la flèche $\calM \to \Fil(\calM)$ est surjective (sur les modules
sous-jacents et sur les $\Fil^r$). Pour l'existence, on remarque que
$\Fil(\calM') = \calM'$ puis que le morphisme $\Fil(f)$ convient. On
raisonne de même avec les catégories $\pUnipN {\tilde S}$ et $\pFilpN
{\tilde S}$.
\end{proof}

\begin{cor}
\label{cor:pcoker}
La catégorie $\pModpN {\tilde S}$ admet des conoyaux.
\end{cor}

\begin{proof}
En vertu du lemme \ref{lem:filadjoint}, il suffit de montrer que
$\pGenpN {\tilde S}$ admet des conoyaux. Soit donc $f : \calM \to
\calM'$ un morphisme dans $\pGenpN {\tilde S}$. On note $\calC$ (resp.
$\Fil^r \calC$) le conoyau du morphisme sous-jacent à $f$ (resp. du
morphisme donné par $f$ sur les $\Fil^r$). Il est alors aisé de vérifier
que les structures supplémentaires sur $\calM'$ passent au quotient pour
faire de $\calC$ un objet de $\pUnipN {\tilde S}$ qui, est en fait dans
$\pGenpN {\tilde S}$.
\end{proof}

\noindent
{\it Remarque.} Si $f : \calM \to \calM'$ est strictement compatible
à $\Fil^r$ dans le sens où $f(\Fil^r \calM) = f(\calM) \cap \Fil^r
\calM'$, alors le conoyau de $f$ calculé dans $\pGenpN S$ est déjà un
objet de $\pModpN {\tilde S}$.

\bigskip

Terminons par un dernier résultat important concernant le foncteur
$\Fil$.

\begin{lemme}
Pour tout objet $\calM \in \pGenpN {\tilde S}$, on a une identification
canonique et fonctorielle :
$$\Tst(\Fil(\calM)) = \hom_{\pUnipN {\tilde S}} (\calM, \hat A).$$
\end{lemme}

\begin{proof}
C'est une conséquence directe du lemme \ref{lem:filadjoint}.
\end{proof}

\begin{cor}
\label{cor:tstcoker}
Soit $f : \calM \to \calM'$ un morphisme dans $\pModpN {\tilde S}$. Si
$\calC$ désigne son conoyau (dans $\pModpN {\tilde S}$), alors
$\Tst(\calC)$ est le noyau de $\Tst(f) : \Tst(\calM') \to \Tst(\calM)$.
\end{cor}

\subsubsection*{Sans la monodromie}

Bien entendu, tout ce qui vient d'être fait peut se refaire sans
difficulté supplémentaire avec les objets \og sans $N$ \fg. On notera
$\pGenp {\tilde S}$, $\pFilp {\tilde S}$ et $\pUnip {\tilde S}$ les
catégories obtenues et encore $\Gen$ et $\Fil$ les foncteurs adjoints.
Ceci ne prête pas à confusion car on vérifie facilement que ces
foncteurs commutent aux foncteurs d'oubli.

\subsection{Éléments nilpotents et objets réduits}
\label{subsec:red}

On introduit ici la catégorie $\RedpN {\tilde S}$ qui va être amené à 
jouer un grand rôle dans la suite.

\begin{deftn}
Soit $\calM \in \pModpN {\tilde S}$. Un élément $x \in \calM$ est
$\Tst$-\emph{nilpotent} si $f(x) = 0$ pour tout $f \in \Tst(\calM)$.
L'ensemble des éléments $\Tst$-nilpotents de $\calM$ est noté
$\NilN(\calM)$.

Le module $\calM$ est dit $\Tst$-\emph{réduit} si $\NilN(\calM) = 0$.
\end{deftn}

\noindent
{\it Remarque.} Les définitions précédentes auraient un sens pour
une catégorie plus générale que $\pModpN {\tilde S}$, mais finalement
que peu d'intérêt pour les applications que nous souhaitons développer
ici. Pour simplifier un peu la présentation, nous nous restreignons
donc au cas de $\pModpN {\tilde S}$.

\bigskip

On vérifie sans mal que cette construction définit un foncteur $\NilN :
\pModpN {\tilde S} \to \pFilpN {\tilde S}$. Soit $\RedN(\calM)$ le
quotient de $\calM$ par $\NilN(\calM)$. Les structures supplémentaires
sur $\calM$ passent au quotient et font de $\RedN(\calM)$ un objet de
$\pModpN {\tilde S}$. De plus, l'application de passage au quotient
$\calM \to \RedN(\calM)$ induit un isomorphisme $\Tst(\RedN(\calM))
\simeq \Tst (\calM)$. Il en résulte que $\RedN(\calM)$ est un objet
$\Tst$-réduit, ou si l'on préfère que $\RedN \circ\RedN = \RedN$.
L'objet $\RedN(\calM)$ est appelé le $\Tst$-\emph{réduit} de $\calM$.

On note $\pRedpN {\tilde S}$ (resp. $\RedpN {\tilde S}$) la
sous-catégorie pleine de $\pModpN {\tilde S}$ qui regroupe les objets
$\Tst$-réduits (resp. $\Tst$-réduits et de type fini sur $S$). On a le
diagramme suivant :
\begin{equation}
\label{eq:diagred}
\raisebox{0.5\depth}{\xymatrix @C=50pt {
\ModpN {\tilde S} \ar@{^(->}[r] \ar[d]_-{\RedN} & \pModpN {\tilde S}
\ar@{->>}[d]^-{\RedN} \\
\RedpN {\tilde S} \ar@{^(->}[r] & \pRedpN {\tilde S}
\ar@/^1em/@{^(->}[u]}}
\end{equation}
où les flèches $\toinj$ symbolisent les inclusions et la flèche 
$\tosurj$
un foncteur essentiellement surjectif. Par un argument analogue à celui
utilisé dans la preuve du lemme \ref{lem:filadjoint}, on obtient :

\begin{lemme}
\label{lem:redadjoint}
Le foncteur $\RedN : \pModpN {\tilde S} \to \pRedpN {\tilde S}$ est un
adjoint à gauche du foncteur d'inclusion.
\end{lemme}

\begin{cor}
\label{cor:redcoker}
Les catégories $\pRedpN {\tilde S}$ et $\RedpN {\tilde S}$ admettent
des conoyaux. De plus, si $f$ est un morphisme dans une de ces deux
catégories, et si $\calC$ est le conoyau de $f$, alors $\Tst(\calC)$
s'identifie au noyau de $\Tst(f)$.
\end{cor}

\begin{proof}
Soit $f : \calM \to \calM'$ un morphisme dans $\pRedpN {\tilde S}$. Le
corollaire \ref{cor:pcoker} assure que $f$ admet un conoyau $\calC$ dans
$\pModpN {\tilde S}$. Le lemme \ref{lem:redadjoint} montre que
$\RedN(\calC)$ est un conoyau de $f$ dans $\pRedpN {\tilde S}$. Par
ailleurs, si $\calM'$ est de type fini, il en est de même de $\calC$
puis de $\RedN(\calC)$ puisque ce sont des quotients successifs de
$\calM'$. Donc, si $f$ est un morphisme dans $\RedpN {\tilde S}$, son
conoyau dans $\pRedpN {\tilde S}$ est un objet de $\RedpN{\tilde S}$.
Ainsi $\RedpN {\tilde S}$ admet, elle aussi, des conoyaux.

La propriété de compatibilité au foncteur $\Tst$ résulte du corollaire
\ref{cor:tstcoker} et de l'identification canonique $\Tst(\Red(\calC))
\simeq \Tst(\calC)$.
\end{proof}

\subsubsection*{Sans la monodromie}

Évidemment, il est possible de rejouer la chanson en omettant partout
l'opérateur $N$. Si $\calM$ est un objet de $\pModp {\tilde S}$, on dit
que $x \in \calM$ est $\Tqst$-\emph{nilpotent} si $f(x) = 0$ pour tout
$f \in \Tqst(\calM)$ ; on note $\Nil(\calM)$ l'ensemble des éléments
$\Tqst$-nilpotents et $\Red(\calM) = \calM / \Nil(\calM)$. La projection
canonique $\calM \to \Red(\calM)$ est envoyé sur un isomorphisme par le
foncteur $\Tqst$. Les équivalents du lemme \ref{lem:redadjoint} et du
corollaire \ref{cor:redcoker} sont encore vrais dans ce contexte et on
définit de façon analogue les catégories $\pRedp {\tilde S}$ et $\Redp
{\tilde S}$ ; elles apparaissent dans un diagramme analogue à
(\ref{eq:diagred}).

Question terminologie, un objet $\calM$ pour lequel $\Nil(\calM) = 0$
est dit $\Tqst$-\emph{réduit} et $\Red(\calM)$ est encore appelé le
$\Tqst$-\emph{réduit} de $\calM$. Malheureusement, si $N$ est un objet
de $\pModpN {\tilde S}$, les notions \og $\Tqst$-\emph{nilpotent} \fg\
et \og $\Tst$-\emph{nilpotent} \fg\ ne coïncident pas ; il est donc
nécessaire de faire la distinction dans l'écriture et la terminologie.
On a malgré tout le lemme suivant.

\begin{lemme}
\label{lem:nilnilN}
Soit $\calM \in \pModpN {\tilde S}$. Alors
$$\NilN(\calM) = \{ x \in \calM \, / \, \forall i \geq 0, \,
N^i(x) \in \Nil(\calM) \}.$$
En particulier $\NilN(\calM) \subset \Nil(\calM)$ et la projection
$\calM \to \Red(\calM)$ se factorise par $\RedN(\calM)$.

De plus $\Nil(\RedN(\calM))$ est l'image de $\Nil(\calM)$ dans
$\RedN(\calM)$, et $\Red \circ \RedN (\calM) = \Red(\calM)$.
\end{lemme}

\begin{proof}
La première partie du lemme est une conséquence directe du lemme
\ref{lem:tsttqst} et de la remarque faite dans sa démonstration.
De $\NilN(\calM) \subset \Nil(\calM)$, on déduit que la projection
$\calM \to \RedN(\calM)$ induit un isomorphisme après application de
$\Tqst$. La dernière partie du lemme résulte facilement de cette
remarque.
\end{proof}

\subsection{Des équivalences de catégories}
\label{subsec:equivs}

Le but de cette sous-section est de démontrer le théorème suivant.

\begin{theo}
\label{theo:equivs}
Les foncteurs $\Red : \Modp {\tilde S} \to \Redp {\tilde S}$ et $\RedN :
\ModpN {\tilde S} \to \RedpN {\tilde S}$ sont des équivalences de
catégories.
\end{theo}

\noindent
{\it Remarque.} Combiné à ce qui a été développé précédemment, ce
théorème permet de construire un objet de $\Modp {\tilde S}$ (resp.
$\ModpN {\tilde S}$) à partir de n'\emph{importe} quel objet de type 
fini
de $\pUnip {\tilde S}$ (resp.  $\pUnipN {\tilde S}$), simplement en lui
appliquant successivement les foncteurs $\Fil$, $\Gen$, $\Red$ (resp.
$\RedN$) puis $\Red^{-1}$ (resp. $\RedN^{-1}$). Ceci nous sera
particulièment utile dans la suite pour mener à bien un certain nombre
de constructions.

\subsubsection*{Définition des quasi-inverses}

Soit $\Mod : \Redp {\tilde S} \to \Modp {\tilde S}$ le foncteur défini
par la formule usuelle :
$$\Mod(\calM) = \tilde S \otimes_{(\phi), k[u]/u^e} \frac{\Fil^r \calM}
{u^e \Fil^r \calM}$$
les structures additionnelles s'obtenant comme suit. Remarquons tout
d'abord que l'on dispose d'un morphisme $\tilde S$-linéaire $\pr = \id
\otimes \phi_r : \Mod(\calM) \to \calM$. Celui-ci est surjectif puisque
par hypothèse $\phi_r(\Fil^r \calM)$ engendre $\calM$. On pose $\Fil^r
\Mod(\calM) = \pr^{-1} (\Fil^r \calM)$ et $\phi_r(x) = 1 \otimes 
\pr(x)$ pour $x \in \Fil^r\Mod(\calM)$. Ainsi, $\Mod(\calM)$ est un 
objet $\pModp {\tilde S}$, qui est bien entendu de type fini de $\tilde 
S$. Il reste à voir qu'il est bien dans $\Modp {\tilde S}$, c'est-à-dire 
qu'il est libre sur $\tilde S$. Cela résulte directement de la 
proposition \ref{prop:redlibre} (ci-dessous) dont la démonstration est 
basée sur le lemme suivant.

\begin{lemme}
\label{lem:nula0}
Soit $y \in \Fil^r \hat A_0$. On suppose $\phi_r(y) \neq 0$. Alors
$u^{e-1} y \neq 0$.
\end{lemme}

\begin{proof}
Comme $\phi_r(y)$ est non nul, $y$ est lui-même non nul et possède une
valuation $p$-adique $v$ bien définie (on rappelle que $\hat A_0$ est
isomorphe en tant qu'anneau à $\O_{\bar K}/p$). De $\phi_r(y) \neq 0$, 
il découle $pv - r < 1$. Il s'ensuit $v < \frac{r+1} p \leq 1 - \frac 1 
p$ puis $u^{e-1} y = \pi_1^{e-1} y \neq 0$.
\end{proof}

\begin{prop}
\label{prop:redlibre}
Soit $\calM \in \Redp {\tilde S}$. Alors $\Fil^r \calM/u^e \Fil^r
\calM$ est libre sur $k[u]/u^e$.
\end{prop}

\begin{proof}
Puisque $\calM$ est de type fini, il existe un $k[[u]]$-module libre de
rang fini $\hat \calM$ muni d'un morphisme surjectif $\pr : \hat \calM
\to \calM$. On peut en outre supposer que $\hat \calM$ est de rang
minimal, ce qui entraîne facilement \emph{via} la théorie des diviseurs
élémentaires que $\ker \pr \subset u \hat \calM$, c'est-à-dire que $\pr$
induit un isomorphisme $\hat \calM / u\hat \calM \simeq \calM / u
\calM$. Notons $d$ le rang de $\hat \calM$. Soit $\Fil^r \hat \calM$
l'image réciproque par $f$ de $\Fil^r \calM$. Il existe $\hat x_1,
\ldots, \hat x_d$ une base de $\Fil^r \hat \calM$ et des entiers $n_1,
\ldots, n_d$ tels que $u^{n_1} \hat x_1, \ldots, u^{n_d} \hat x_d$ soit
une base de $\ker \pr$. Soit $x_i \in \calM$ l'image de $\hat x_i$. Du
fait que $\phi_r(\Fil^r \calM)$ engendre $\calM$, on déduit que la
famille des $\phi_r(x_i) \mod u$ engendre $\calM / u \calM$. Comme elle
est de cardinal $d$, elle en est une base. En particulier, aucun des
$\phi_r(x_i)$ n'est nul. Comme $\calM$ est $\Tqst$-réduit, il existe
pour tout $i$, un élément $f_i \in \Tqst(\calM)$ tel que
$f_i(\phi_r(x_i)) \neq 0$. Fixons un indice $i$ et posons $y = f_i(x_i)
\in \Fil^r \hat A_0$. On a $\phi_r(y) \neq 0$, et donc par le lemme
\ref{lem:nula0}, $f_i(u^{e-1} x_i) = u^{e-1} y \neq 0$. On en déduit que
$u^{e-1} x_i$ est lui-même non nul, \emph{i.e} $u^{e-1} \hat x_i \not\in
\ker \pr$. Ainsi $n_i \geq e$, et comme ceci est vrai pour tout $i$, la
proposition est démontrée.
\end{proof}

On peut procéder de même lorsque l'opérateur de monodromie est présent :
à un objet $\calM$ de $\RedpN {\tilde S}$, on associe $\ModN(\calM)$
défini par les mêmes formules que $\Mod(\calM)$ auxquelles on ajoute
$N(s \otimes x) = c^{-1} \otimes u^e N(x) - u s' \otimes x$ ($s \in
\tilde S$, $x \in \Fil^r \calM/u^e \Fil^r \calM$) où $s'$ désigne la
dérivée de $s$ vu comme polynôme en $u$. La liberté (sur $k[u]/u^e$) du
quotient $\Fil^r \calM / u^e \Fil^r \calM$ s'obtient alors comme dans la
preuve de la proposition \ref{prop:redlibre} en remplaçant (\og $\Tqst$
\fg\ par \og $\Tst$ \fg\ et) la référence au lemme \ref{lem:nula0} par
une référence au lemme suivant.

\begin{lemme}
Soit $y \in \Fil^r \hat A$. On suppose que $\phi_r(y) \neq 0$. Alors
$u^{e-1} y \neq 0$.
\end{lemme}

\begin{proof}
Comme $u$ et $\pi_1$ diffèrent d'une unité, il suffit de montrer que
$\pi_1^{e-1} y$ ne s'annule pas. Par définition, $y$ s'écrit de façon
unique sous la forme :
$$y = a_0 + a_1 X + a_2 \frac{X^2} 2 + \cdots + a_n \frac {X^n}{n!}$$
pour un certain entier $n$, avec $a_i \in \Fil^{r-i} \hat A_0$ (où par
convention $\Fil^j \hat A_0 = \hat A_0$ pour $j \leq 0$). On a alors
$\phi_r(y) = \phi_r(a_0) + \phi_{r-1} (a_1) Y + \cdots + \phi_0(a_r)
\frac{Y^r}{r!}$ (où on rappelle que $Y = \frac {(1+X)^p-1} p$). Comme 
$\phi_r(y) \neq 0$, il existe un indice $i \in \{0, \ldots, r\}$ tel que
$\phi_{r-i}(a_i) \neq 0$. Le lemme \ref{lem:nula0} assure alors que
$\pi_1 ^{e-1} a_i$ est non nul, et donc qu'il en est de même de $\pi_1
^{e-1} y$.
\end{proof}

\subsubsection*{Calcul des composées}

Nous allons montrer que $\Mod$ (resp. $\ModN$) est un quasi-inverse de
$\Red$ (resp. $\RedN$) simplement en calculant les composées dans les
deux sens. Nous commençons par le calcul de $\Mod \circ \Red$ (resp.
$\ModN \circ \RedN$) largement basé sur le lemme suivant.

\begin{lemme}
\label{lem:modnilfil}
Soit $\calM$ un objet de $\Modp {\tilde S}$ (resp. $\ModpN {\tilde S}$).
Alors $\Nil(\calM)$ (resp. $\NilN(\calM)$) est inclus dans $u^e \Fil^r
\calM$.
\end{lemme}

\begin{proof}
Grâce à l'inclusion $\NilN(\calM) \subset \Nil(\calM)$ énoncée dans
le lemme \ref{lem:nilnilN}, il suffit de montrer le résultat lorsqu'il
n'y a pas d'opérateur de monodromie.

Notons $d$ le rang de $\calM$ comme $\tilde S$-module. Montrons tout
d'abord que $\Nil(\calM) \subset u \calM$. On remarque à cet effet que
$\Red(\calM)/u \Red(\calM)$ est naturellement un quotient de $\calM / u
\calM$. Soit $d'$ sa dimension (sur $k$). On a $d' \leq d$. Par
ailleurs, en relevant une famille génératrice de $\Red(\calM)/u
\Red(\calM)$, on montre que $\Red(\calM)$ est lui aussi engendré (sur
$\tilde S$) par $d'$ éléments. Autrement dit, il existe un morphisme
surjectif $\tilde S$-linéaire $f : \calN \to \Red(\calM)$ où $\calN$ est
un $\tilde S$-module libre de rang $d'$. Posons $\Fil^r \calN = f^{-1}
(\Fil^r \Red(\calM))$. En appliquant la théorie des diviseurs
élémentaires à l'inclusion $\Fil^r \calN \subset \calN$, on définit
facilement une application $\phi_r : \Fil^r \calN \to \calN$ qui fait de
$\calN$ un objet de $\Modp {\tilde S}$ et de $f$ un morphisme dans cette
catégorie. On en déduit une injection $\Tqst(\calM) \simeq
\Tqst(\Red(\calM)) \toinj \Tqst(\calN)$. L'espace de départ est un
$\F_p$-espace vectoriel de dimension $d$, alors que celui d'arrivée est
de dimension $d'$. Il en résulte $d \leq d'$, et puis $d = d'$. Ainsi
$\Red(\calM) / u \Red(\calM) = \calM / u \calM$, ce qui entraîne
$\Nil(\calM) \subset u \calM$ comme annoncé.

La fin de la preuve consiste à répéter l'argument de la preuve de la
proposition \ref{prop:redlibre} en appliquant la théorie des diviseurs
élémentaires non pas à l'inclusion $\ker \pr \subset \Fil^r \hat \calM$
mais aux deux sous-modules $\pr^{-1}(\Nil(\calM))$ et $\Fil^r \hat
\calM$, les entiers $n_i$ qui apparaissent étant alors \emph{a priori}
relatifs. Nous laissons au lecteur le soin de faire ces modifications
mineures.
\end{proof}

\begin{cor}
Pour tout objet $\calM$ de $\Modp {\tilde S}$ (resp. $\ModpN {\tilde
S}$), on a un isomorphisme canonique et fonctoriel :
$$\Mod \circ \Red(\calM) \simeq \calM \quad \text{(resp. }
\ModN \circ \RedN(\calM) \simeq \calM \text{).}$$
\end{cor}

\begin{proof}
On n'écrit la preuve que pour $\calM \in \Modp {\tilde S}$, l'autre
cas étant entièrement analogue. Par le lemme \ref{lem:modnilfil}, le
morphisme canonique $\calM \to \Red(\calM)$ induit un isomorphisme
$$\Fil^r \calM / u^e \Fil^r \calM \simeq \Fil^r \Red(\calM) / u^e \Fil^r
\Red(\calM).$$
Il s'ensuit que $\Mod \circ \Red(\calM)$ s'identifie à $\tilde S
\otimes_{(\phi), k[u]/u^e} \Fil^r \calM/u^e \Fil^r \calM$. Le morphisme
$\pr = \id \otimes \phi_r$ induit une application surjective (compatible
aux structures additionnelles) $\Mod \circ \Red(\calM) \to \calM$. Comme
les espaces de départ et d'arrivée sont des $\tilde S$-modules libres de
même rang, c'est un isomorphisme et le corollaire est démontré.
\end{proof}

On termine à présent la preuve du théorème \ref{theo:equivs} en faisant
le calcul de $\Red \circ \Mod$ (resp. $\RedN \circ \ModN$).

\begin{prop}
Pour tout objet $\calM$ de $\Redp {\tilde S}$ (resp. $\RedpN {\tilde
S}$), on a un isomorphisme canonique et fonctoriel :
$$\Red \circ \Mod(\calM) \simeq \calM \quad \text{(resp. }
\RedN \circ \ModN(\calM) \simeq \calM \text{).}$$
\end{prop}

\begin{proof}
Comme précédemment, on ne donne la preuve que pour $\Redp {\tilde S}$.
Posons $\calM' = \Mod(\calM)$ et notons $\pr = \id \otimes \phi_r :
\calM' \to \calM$ la projection canonique. Nous allons montrer que
$\Red(\pr)$ est un isomorphisme (ce qui permettra de conclure). Étant
donné que $\pr$ est surjectif, il est clair déjà que $\Red(\pr)$ l'est
aussi. D'après le lemme \ref{lem:redinjectif}, pour prouver qu'il est
injectif, il suffit de montrer que $\Tqst(\pr) : \Tqst(\calM')
\toinj \Tqst(\calM)$ est surjectif. Or, du fait que $\calM$ est 
$\Tqst$-réduit, on a
$$\big\{ x \in \calM' \, / \, \forall h \in \Tqst(\calM), \,
h(x) = 0 \big\} = \ker \pr \subset u \hat \calM'$$
la dernière inclusion se vérifiant aisément à la main (on pourra
remarquer que $\pr$ induit une application surjective --- et donc un
isomorphisme --- de $\calM' / u \calM'$ dans $\calM / u
\calM$). Le corollaire \ref{cor:decoup} s'applique alors et termine la
preuve.
\end{proof}

\subsection{La structure de pylonet en termes d'objets réduits}
\label{subsec:explicite}

Le théorème \ref{theo:tstpylonet} montre que la catégorie $\ModpN 
{\tilde S}$ admet une structure riche. Le but de cette partie est de
la comprendre en termes d'objets $\Tst$-réduits, c'est-à-dire sous
l'équivalence de catégories $\RedN$. On commence par un lemme très
simple qui nous sera utile à plusieurs reprises dans la suite.

\begin{lemme}
\label{lem:sousred}
Soit $f : \calM \to \calM'$ un morphisme dans la catégorie $\pModp
{\tilde S}$ (resp. $\ModpN {\tilde S}$). On suppose que $f$ est injectif
et que $\calM'$ est $\Tqst$-réduit (resp. $\Tst$-réduit). Alors $\calM$
l'est aussi.
\end{lemme}

\begin{proof}
De la fonctorialité de $\Nil$ (resp. $\NilN$), on déduit que $f$ envoie
$\Nil(\calM)$ sur $\Nil(\calM') = 0$ (resp. $\NilN(\calM)$ sur
$\NilN(\calM') = 0$). Le lemme résulte alors de l'injectivité de $f$.
\end{proof}

\subsubsection*{Noyaux, images et conoyaux}

Considérons $f : \calM \to \calM'$ un morphisme dans $\RedpN {\tilde 
S}$. D'après les axiomes \Ax 2 et \Ax {2*}, $f$ admet un noyau et un
conoyau dans $\RedpN {\tilde S}$ que l'on note respectivement $\calK$ et 
$\calC$. Il résulte facilement des diverses propriétés d'adjonction 
démontrées précédemment et du lemme \ref{lem:sousred} que $\calK$ et
$\calC$ se calculent explicitement comme suit :
\begin{equation}
\calK = \Gen(\ker f) \quad \text{et} \quad
\calC = \RedN \circ \Fil (\coker f)
\end{equation}
où $\ker f$ et $\coker f$ sont respectivement le noyau et le conoyau de 
$f$ au sens usuel. La formule pour le noyau est intéressante car elle ne 
fait plus intervenir à aucun moment les représentations galoisiennes ! 
Malheureusement, ce n'est pas le cas pour le conoyau puisque la formule
que l'on obtient fait apparaître le foncteur $\RedN$ dans la définition
duquel intervient de façon essentielle le foncteur $\Tst$. Malgré tout,
avec un peu de pratique, il ne semble pas très difficile d'avoir une
intuition du résultat final et de le démontrer \emph{a posteriori}. Dans 
tous les cas, si le calcul pose vraiment un problème, on a toujours 
comme recours l'utilisation de la dualité.

Il n'est sans doute finalement pas anodin de remarquer que l'on dispose 
également d'une formule --- qui plus est très simple --- pour le calcul 
de l'image (c'est-à-dire le noyau du conoyau) dans $\RedN {\tilde S}$ 
puisque celle-ci s'identifie à l'image usuelle.

\subsubsection*{La relation d'ordre}

Soit $T$ une $\F_p$-représentation galoisienne dans l'image essentielle 
de $\Tst$. D'après le théorème \ref{theo:tstpylonet}, la fibre au-dessus 
de $T$ (c'est-à-dire l'ensemble des $\calM \in \pModpN {\tilde S}$ dont 
l'image par $\Tst$ est $T$) a une structure de treillis. Nous allons
voir que celle-ci se comprend immédiatement en termes d'objets 
$\Tst$-réduits. Le lemme suivant (très facile) est la clé de cette
compréhension.

\begin{lemme}
\label{lem:redinjectif}
Soit $f : \calM \to \calM'$ un morphisme dans la catégorie $\pModp
{\tilde S}$ (resp. $\pModpN {\tilde S}$). On suppose $\Tqst(f)$
(resp. $\Tst(f)$) surjectif. Alors $\Red(f)$ (resp. $\RedN(f)$) est
injectif.
\end{lemme}

\begin{proof}
On peut bien sûr supposer que $\calM$ et $\calM'$ sont $\Tqst$-réduits
(resp $\Tst$-réduits), et on veut alors montrer que $f$, lui-même, est
injectif. Soit $x \in \ker f$. Par hypothèse, tout $g \in \Tqst(\calM)$
(resp. $g \in \Tst(\calM)$) se factorise par $f$ et donc s'annule sur
$x$. On en déduit que $x \in \Nil(\calM)$ (resp. $x \in \NilN(\calM)$)
et donc que $x = 0$ comme souhaité.
\end{proof}

Il résulte de ce lemme que tous les morphismes dans la catégorie fibre 
$\calF_T$ sont injectifs. Ainsi l'ordre que l'on cherche à décrire 
correspond simplement à l'inclusion naturelle sur les objets 
$\Tst$-réduits. On a en outre un petit rabiot qui montre en un certain 
sens qu'il n'y a pas de \og trous \fg.

\begin{prop}
\label{prop:redinterm}
Soient $\calM' \subset \calM''$ des objets de $\RedpN {\tilde S}$. On
note $\iota$ l'inclusion de $\calM'$ dans $\calM''$ et on suppose que
$\Tst(\iota)$ est un isomorphisme. Soit $\calM \in \pModpN {\tilde S}$
tel que $\calM' \subset \calM \subset \calM''$. Alors $\calM \in \RedpN
{\tilde S}$ et les deux flèches déduites par fonctorialité :
$$\Tst(\calM'') \to \Tst(\calM) \to \Tst(\calM')$$
sont des isomorphismes.
\end{prop}

\begin{proof}
Soit $d$ la dimension sur $\F_p$ de $T = \Tst(\calM')$ ; c'est aussi le
rang de $\Mod(\calM')$, d'où on déduit que $\calM'$ est engendré par au
plus $d$ éléments. Comme par hypothèse $\Tst(\calM') \simeq
\Tst(\calM'')$, la même conclusion vaut par $\calM''$. Étant donné que
$\tilde S$ est un anneau principal (non intègre), on en déduit que
$\calM$ est lui aussi engendré par au plus $d$ éléments. En particulier,
il est de type fini. D'autre part, le lemme \ref{lem:sousred} entraîne
que $\calM$ est $\Tst$-réduit. Ainsi, $\calM$ est bien un objet de $\RedpN
{\tilde S}$.

La composée des morphismes $\Tst(\calM'') \to \Tst(\calM) \to
\Tst(\calM')$ est bijective car elle s'identifie à $\Tst(\iota)$. Par
ailleurs, la dimension (sur $\F_p$) de $\Tst(\calM)$ est égale au rang
de $\Mod(\calM)$ et donc majorée par $d$. La proposition en résulte.
\end{proof}

Forts de ces résultats, il devient possible de décrire le calcul des 
bornes supérieures et inférieures dans la fibre $F_T$. Donnons-nous pour 
cela $\calM \in \RedpN {\tilde S}$ un objet maximal et posons $T = 
\Tst(\calM)$. Soient également $\calM_1$ et $\calM_2$ deux objets de 
$\RedpN {\tilde S}$ dont l'image par $\Tst$ s'identifie à $T$. D'après 
la définition des objets maximaux et le lemme \ref{lem:redinjectif},
$\calM_1$ et $\calM_2$ apparaissent comme des sous-objets de $\calM$. 
Dans ces conditions, la borne supérieure de $\calM_1$ et $\calM_2$ 
s'identifie à leur somme dans $\calM$ tandis que leur borne inférieure 
est $\Gen(\calM_1 \cap \calM_2)$. (La vérification est immédiate et 
laissée au lecteur.) On notera que ceci vaut encore pour n'importe 
quelle famille $(\calM_i)_{i \in I}$ par nécessairement finie.

\subsection{Une formule de réciprocité}
\label{subsec:pleinfid}

Nous avons vu que la restriction de $\Tst$ à $\MaxpN {\tilde S}$ est 
pleinement fidèle. Ainsi, étant donnée une représentation $T$ dans 
l'image essentielle de $\Tst$, il y a un unique objet de $\MaxpN {\tilde 
S}$ qui lui correspond. Nous montrons ci-après qu'il est possible de
retrouver cet objet par une formule explicite.

\subsubsection*{Le foncteur $\Mst$}

Soit $T$ une $\F_p$-représentation de $G_K$. Le $\tilde S$-module
$\hom_{\F_p[G_K]} (T, \hat A)$ hérite des structures supplémentaires de
$\hat A$, ce qui en fait un objet de la catégorie $\pFilpN {\tilde S}$.
L'association
$$T \mapsto \Mst(T) = \Gen\big(\hom_{\F_p[G_K]} (T, \hat A)\big)$$
définit un foncteur contravariant $\Mst : \Rep_{\F_p}(G_K) \to \pModpN
{\tilde S}$. Pour $\calM \in \pModpN {\tilde S}$ et $T \in
\Rep_{\F_p}(G_K)$, on dispose en outre d'applications de bidualité
$$\alpha_\st(\calM) : \calM \to \Mst \circ \Tst(\calM)
\quad \text{et} \quad
\beta_\st(T) : T \to \Tst \circ \Mst(T)$$
qui sont des morphismes respectivement dans $\pModpN {\tilde S}$ et
$\Rep_{\F_p}(G_K)$. De plus, en déroulant les définitions, on obtient
$\ker \alpha_\st(\calM) = \NilN(\calM)$, tandis qu'une vérification
simple montre que $\Tst(\alpha_\st(\calM)) \circ \beta_\st(\Tst(\calM))
= \id_{\Tst(\calM)}$ et donc que $\Tst(\calM)$ apparaît (\emph{via}
$\beta_\st(\Tst(\calM))$) comme un facteur direct de $\Tst \circ
\Mst \circ \Tst(\calM)$.

\begin{lemme}
\label{lem:bidual}
Soit $\calM$ un objet de $\ModpN {\tilde S}$. Posons $T = \Tst(\calM)$
et $\calM' = \Mst(T)$. Alors :
$$\big\{ x \in \calM' \, / \, \forall h \in T, \, \beta_\st(T)(h)(x) =
0\big\}$$
est réduit à $0$.
\end{lemme}

\begin{proof}
C'est évident après avoir remarqué que les $x \in \calM'$ sont des
morphismes ($G_K$-équivariants) de $T$ dans $\hat A$ et que
$\beta_\st(T)(h)(x)$ n'est rien d'autre que $x(h)$.
\end{proof}

\begin{cor}
\label{cor:msttst1}
Soit $\calM \in \ModpN {\tilde S}$. Alors $\Mst \circ \Tst(\calM)$ est
$\Tst$-réduit.
\end{cor}

\begin{prop}
\label{prop:msttst2}
Soit $T$ une $\F_p$-représentation de dimension finie de $G_K$. Alors
$\Mst(T)$ est un $\tilde S$-module de type fini.
\end{prop}

\begin{proof}
Soit $L$ une extension finie de $K$ dont le groupe de Galois absolu,
noté $G_L$, agit trivialement sur $T$. Quitte à agrandir $L$, on peut
supposer $\pi_1 \in L$. Les morphismes $G_K$-équivariants de $T$ dans
$\hat A$ prennent alors leurs valeurs dans $\hat A^{G_L}$, d'où on
déduit $\Mst(T) = \Gen(\hom_{\F_p[G_K]} (T, \hat A^{G_L}))$. Une
induction transfinie à partir de la définition de $\Gen$ montre
directement l'inclusion
$$\Mst(T) \subset \hom_{\F_p[G_K]} (T, \Gen(\hat A^{G_L})).$$
Ainsi, puisque $\tilde S$ est n\oe thérien, il suffit pour conclure de
montrer que $\Gen(\hat A^{G_L})$ est de type fini sur $\tilde S$. Nous
allons en fait montrer que $\Gen_1(\hat A^{G_L})$ est déjà de type fini
sur $\tilde S$. Soit $a = a_0 + a_1 X + \cdots + a_n \frac{X^n}{n!} \in
\hat A^{G_L} \cap \Fil^r \hat A$. Ici, donc, les $a_i$ sont \emph{a
priori} des éléments de $\Fil^{r-i} \hat A_0$. Étant donné que $G_L$
n'agit pas sur $X$ (on rappelle que l'on a supposé $\pi_1 \in L$), le
fait que $a \in \hat A^{G_L}$ implique que chacun des $a_i$ est lui-même
fixe par $G_L$. Soit $L_1$ une extension de $L$ obtenue en ajoutant une
racine $p$-ième d'une uniformisante de $L$, et soit $v_p$ la valuation
$p$-adique sur $\bar K$ normalisée par $v_p(p) = 1$. D'après les
résultats de \cite{leborgne}, on peut écrire $a_i = b_i + c_i$ avec $b_i
\in \O_{L_1}/p$ et $v_p(c_i) \geq 1 - \frac 1 {p(p-1)}$. La dernière
condition sur la valuation montre que $c_i$ (et donc aussi $b_i$) est
dans $\Fil^{r-i} \hat A_0$ et que $\phi_{r-i}(c_i) = 0$. Ainsi 
trouve-t-on
$$\phi_r(a) = \phi_r(b_0) + \phi_{r-1}(b_1) Y + \cdots + \phi_1(b_{r-1})
\frac{Y^{r-1}}{(r-1)!} + \phi(b_r) \frac{Y^r}{r!}$$
où on rappelle que $Y = \frac{(1+X)^p-1} p$. On en déduit que 
$\Gen_1(\hat A^{G_L})$ est inclus dans le $\tilde S$-module engendré par 
les éléments de la forme $\phi_{r-i} (b) Y^i$ pour $b \in \O_{L_1}/p$ et 
$0 \leq i \leq r$. Comme $\O_{L_1}/p$ est de type fini sur $\O_K/p$ 
(puisque $L_1$ est une extension finie de $K$), on a bien montré que 
$\Gen_1(\hat A^{G_L})$ est de type fini sur $\tilde S$.
\end{proof}

\subsubsection*{La composée $\Mst \circ \Tst$}

Le corollaire \ref{cor:msttst1} combiné à la proposition
\ref{prop:msttst2} montre que la composée $\Mst \circ \Tst$ definit
un foncteur de $\ModpN {\tilde S}$ dans $\RedpN {\tilde S}$.

\begin{theo}
\label{theo:reciprocite}
Pour tout $\calM \in \ModpN {\tilde S}$, le morphisme
$$\alpha_\st(\Max(\calM)) : \Max(\calM) \to \Mst \circ \Tst(\Max(\calM))
\simeq \Mst \circ \Tst(\calM)$$
est surjectif (et donc induit un isomorphisme entre $\Red \circ
\Max(\calM)$ et $\Mst \circ \Tst(\calM)$).
\end{theo}

\begin{proof}
Quitte à remplacer $\calM$ par $\Max(\calM)$, on peut bien sûr supposer
que $\calM$ est maximal. Notons $T = \Tst(\calM)$, $\calM' = \Mst(T)$,
$T' = \Tst(\calM')$ et $f : \Red(\calM) \toinj \calM'$ le morphisme
(injectif) induit par $\alpha_\st(\calM)$. Il s'agit de montrer que $f$
est un isomorphisme.

On a vu que $T$ apparaît \emph{via} $\beta_\st(T)$ comme une
sous-représentation (et même un facteur direct) de $T'$. Par ailleurs,
le lemme \ref{lem:bidual} donne :
$$\big\{ x \in \Mod(\calM') \, / \, \forall h \in T, \, h(x) = 0\big\}
= \NilN(\Mod(\calM')) \subset u \Mod(\calM')$$
la dernière inclusion provenant du lemme \ref{lem:modnilfil}. Le
corollaire \ref{cor:decoup} entraîne $T = T'$, \emph{i.e.} $\Tst(f)$ est
un isomorphisme. Étant donné que $\calM$ est maximal, ceci implique
l'existence d'un morphisme $g : \calM' \to \Red(\calM)$ tel que $g \circ
f = \id_{\Red(\calM)}$. Par le lemme \ref{lem:redinjectif}, $g$ est
injectif. Il s'ensuit, comme annoncé, que $f$ est un isomorphisme.
\end{proof}

\noindent
{\it Remarque.}
Le théorème donne une formule qui permet de retrouver $\Max(\calM)$ à 
partir de $\Tst(\calM)$. L'intérêt de disposer d'une telle formule est 
de pouvoir relever facilement au niveau de $\MaxpN {\tilde S}$ des 
applications vivant \emph{a priori} sur les représentations 
galoisiennes. Par exemple, voici comment on peut l'utiliser pour donner 
une seconde preuve de la pleine fidélité de $\Tst : \MaxpN {\tilde S} 
\to \Rep _{\F_p} (G_K)$. Soient $\calM$ et $\calM'$ dans $\MaxpN {\tilde 
S}$. Posons $T = \Tst (\calM)$ et $T' = \Tst(\calM')$. D'après les 
théorèmes \ref{theo:equivs} et \ref{theo:reciprocite}, la composée : 
$$\hom_{\ModpN {\tilde S}} (\calM, \calM') \stackrel{v}{\longrightarrow} 
\hom_{\F_p[G_K]} (T', T) \stackrel{w}{\longrightarrow} \hom_{\RedpN 
{\tilde S}} (\Mst(T), \Mst (T'))$$
est un isomorphisme. On veut montrer que $v$ est un isomorphisme, et 
pour cela il suffit de justifier que $w$ est injective. Or, tout $f \in 
\ker w$ s'insère dans le diagramme commutatif suivant :

$$\xymatrix @C=60pt {
T' \ar[r]^-{\beta_\st(T')} \ar[d]_-{f} & \Tst \circ \Mst(T')
\ar[d]^-{0} \\
T \ar[r]^-{\beta_\st(T)} & \Tst \circ \Mst(T) }$$
à partir duquel on déduit directement que $f = 0$ en utilisant
l'injectivité de $\beta_\st(T)$ (on rappelle que ce dernier morphisme
admet $\Tst(\alpha_\st(\calM))$ pour rétraction).

Pour certaines applications même, le théorème \ref{theo:reciprocite} 
peut s'appliquer alors que la pleine fidélité ne sera \emph{a priori} 
d'aucun secours. C'est typiquement ce qui se passe lorsque l'on souhaite 
relever des applications qui ne sont linéaires (mais par exemple 
semi-linéaires), ou que l'on s'intéresse à des représentations dans des
espaces vectoriels munis de structures supplémentaires (par exemple une
forme quadratique ou symplectique). Nous verrons une application de cela 
en \ref{subsec:descente}.

\section{Compléments}
\subsection{Variante avec coefficients}

Dans la pratique, il arrive souvent que l'on ait besoin d'étudier les
représentations semi-stables, non pas à coefficients dans $\F_p$ mais
dans une extension finie\footnote{Étant donné que l'on ne s'intéresse
qu'à des représentations de dimension finie, il est toujours possible de
faire cette hypothèse} $E$ de $\F_p$. Une telle représentation $V$ peut
également être vue comme une $\F_p$-représentation munie d'un morphisme
d'anneaux $E \to \End_{\F_p[G_K]} (V)$. On est donc amené à considérer
la catégorie $\ModpN {\tilde S \otimes E}$ dont les objets sont les
couples $(\calM, \nu)$ où $\calM \in \ModpN {\tilde S}$ et $\nu : E \to
\End_{\ModpN {\tilde S}}(\calM)$ est un morphisme d'anneaux. Bien sûr,
la notation est justifiée par le fait que la donnée de $\nu$ équivaut à
une structure de $E$-espace vectoriel --- et donc de $(\tilde S \otimes_
{\F_p} E)$-module\footnote{On prendra garde au fait que le Frobenius sur
$\tilde S \otimes_{\F_p} E$ est bien l'élévation à la puissance $p$ sur
$\tilde S$, mais l'identité sur $E$ !} --- sur $\calM$. Toutefois, pour
ce que nous voulons faire ici, la première description que nous avons
donnée sera plus adaptée.

On dispose bien entendu d'un foncteur oubli $\Oub_E : \ModpN {\tilde S 
\otimes E} \to \ModpN {\tilde S}$, $(\calM, \nu) \mapsto \calM$. Il est 
fidèle et conservatif. Par ailleurs, le foncteur $\Tst$ se prolonge en 
$\TstE : \ModpN {\tilde S \otimes E} \to \Rep_E(G_K)$ obtenu simplement 
en faisant agir $E$ sur $\Tst(\calM)$ \emph{via} $\lambda \cdot x = 
\Tst( \nu(\lambda)) (x)$.

\begin{theo}
\label{theo:pylonetE}
La fibration $\TstE$ est un pylonet (contravariant) additif et autodual.
En outre, si $(\calM, \nu) \in \ModpN {\tilde S \otimes E}$, on a :
$$\Max(\calM, \nu) = (\Max(\calM), \Max(\nu))
\quad \text{et} \quad
\Min(\calM, \nu) = (\Min(\calM), \Min(\nu))$$
où $\Max(\nu) : E \to \End_{\ModpN {\tilde S}}(\Max(\calM))$, 
$\lambda \mapsto \Max(\nu(\lambda))$ (et de même pour $\Min(\nu)$).
\end{theo}

\begin{proof}
Le premier point ne pose aucune difficulté particulière : on peut par 
exemple reprendre la démonstration de la section \ref{sec:application} 
en ajoutant l'action de $E$ à chaque étape, ce que nous laissons au 
lecteur.
Pour la seconde assertion, on remarque d'abord que $\TstE(\Max(\calM), 
\Max(\nu)) \simeq \TstE(\calM, \nu)$. Ainsi, par définition des objets
maximaux, on a un morphisme canonique $f : (\Max(\calM), \Max(\nu)) \to 
\Max(\calM, \nu)$ dans la catégorie $\ModpN {\tilde S \otimes E}$. Le
morphisme $\Oub_E(f)$ s'envoie sur un isomorphisme par $\Tst$ et a pour 
source un objet maximal (de $\ModpN {\tilde S}$). On en déduit que c'est
un isomorphisme, et puis que c'est aussi le cas de $f$ en utilisant la 
conservativité de $\Oub_E$. Le cas des objets minimaux se traite 
pareillement.
\end{proof}

Comme pour $\ModpN {\tilde S}$, on note $\MaxpN {\tilde S \otimes E}
= \Max(\ModpN {\tilde S \otimes E})$ et $\MinpN {\tilde S \otimes E}
= \Min(\ModpN {\tilde S \otimes E})$.

\begin{theo}
\label{theo:pleinfidE}
La restriction de $\TstE$ à $\MaxpN {\tilde S \otimes E}$ (resp. 
$\MinpN {\tilde S \otimes E}$) est pleinement fidèle et son image
essentielle est stable par sous-objets et quotients.
\end{theo}

\begin{proof}
La pleine fidélité est une conséquence directe du théorème
\ref{theo:pleinfid} et de la formule $\Max(\calM, \nu) = (\Max(\calM),
\Max(\nu))$ (resp. $\Min(\calM, \nu) = (\Min(\calM), \Min(\nu))$).
La stabilité découle du résultat analogue pour $\ModpN {\tilde S}$ et de
la pleine fidélité puisque cette dernière permet de relever l'action de
$E$.
\end{proof}

\subsection{Passage à une extension finie, donnée de descente}
\label{subsec:descente}

Dans ce paragraphe, on cherche à comprendre comment les catégories
précédentes (et les représentations qu'elles produisent) se comportent
lorsque l'on change le corps $K$. Pour cela, on fixe $L$ une extension
finie de $K$ dont un note $\O_L$ l'anneau des entiers, $\ell$ le corps
résiduel et $G_L$ le groupe de Galois absolu. On note $L_0$ la plus
grande extension non ramifiée (sur $\Q_p$) contenue dans $L$ ; elle
s'identifie à $W(\ell)[1/p]$. Soit $e_L = [L:L_0]$.

À cette situation, il est attaché de nouvelles catégories de modules
définis sur l'anneau $\tilde S_L = \ell[u]/u^{e_L p}$. Afin d'éviter les
confusions, nous indicerons dans la suite, les constantes, les
catégories et les foncteurs par les lettres $K$ ou $L$ selon qu'elles se
réfèrent au corps $K$ ou $L$ ; par exemple, nous noterons $e_K$ et
$e_L$, $\pi_K$ et $\pi_L$ pour les uniformisantes choisies, $\ModpN
{\tilde S_K}$ et $\ModpN {\tilde S_L}$, ou encore $\Tkst$ et $\Tlst$
(pour ne pas confondre avec $\TstE$).

\subsubsection*{Cas d'une extension non ramifiée}

On suppose d'abord que $L/K$ est non ramifiée. L'uniformisante $\pi \in 
\O_K$ reste une uniformisante de $L$ dont le polynôme minimal sur $L_0$ 
est encore $E(u)$. On a donc $e_K = e_L$ et on note à nouveau $e$ cette 
valeur commune. L'extension des scalaires de $k$ à $\ell$ définit un 
foncteur fidèle $\ModpN {\tilde S_K} \to \ModpN {\tilde S_L}$ et on 
vérifie directement que pour $\calM \in \ModpN {\tilde S_K}$, la flèche 
naturelle $\Tkst(\calM) \to \Tlst(\calM \otimes_k \ell)$ est un 
isomorphisme $G_L$-équivariant.

\begin{prop}
\label{prop:descente}
On se donne $L$ une extension \emph{non ramifiée} de $K$.
Soit $E$ une extension finie de $\F_p$. Soit $T$ une $E$-représentation
de $G_K$. On suppose que la restriction de $T$ à $G_L$ est dans
$\Tlst(\ModpN{\tilde S_L \otimes E})$. Alors $T$ est dans
$\Tkst(\ModpN{\tilde S_K \otimes E})$.
\end{prop}

\begin{proof}
Par les propriétés de pleine fidélité, on se ramène facilement au cas 
\og sans coefficients \fg. La clôture galoisienne $M$ de $L$ est encore 
une extension non ramifiée de $K$ et la restriction de $T$ à $G_M$ 
provient d'un objet de $\ModpN{\tilde S_M \otimes E}$ ; on peut donc 
supposer que $L/K$ est galoisienne. Soit $\calM_L$ l'objet de $\MaxpN 
{\tilde S_L}$ associé à $T_{|G_L}$. On a une action naturelle de $G_K$ 
sur $\Mlst(T)$ donnée par $(\sigma f)(x) = \sigma f(\sigma^{-1} x)$, qui 
se factorise à travers $\Gal(L/K) \simeq \Gal(\ell/k)$ puisque les $f 
\in \Mlst(T)$ sont par définition $G_L$-équivariants. Par ailleurs, la
combinaison des théorèmes \ref{theo:equivs} et \ref{theo:reciprocite} 
assure que $\calM_L = \ModN(\Mlst(T))$. Ceci permet de remonter l'action 
de $\Gal(\ell/k)$ à $\calM_L$. De la nullité de $H^1(\Gal(\ell/k), 
\text{GL}_{erd}(\ell))$ (où $d$ est la dimension de $T$), on déduit  
$\calM_L^{\Gal(\ell/k)} \otimes_k \ell \simeq \calM_L$. On pose alors 
$\calM_K = \calM_L^{\Gal(\ell/k)}$ et on vérifie à la main que  
l'isomorphisme $\Tkst(\calM_K) \simeq \Tlst(\calM_L) = T$ est 
$G_K$-équivariant.
\end{proof}

\subsubsection*{Cas d'une extension modérément ramifiée}

On suppose maintenant que l'extension $L/K$ est totalement et modérément
ramifiée. Notons $n$ son degré ; il est premier avec $p$, et on fixe un
entier $m$ tel que $mn \equiv 1 \pmod p$. On suppose de surcroît que $K$
contient toutes les racines $n$-ièmes de l'unité\footnote{Quitte à
remplacer $K$ par une extension non ramifiée, cette hypothèse est
évidemment toujours satisfaite. Par ailleurs, comme cela a été expliqué
précédemment, le passage à une extension non ramifiée ne pose pas
réellement problème.}. Si $\pi_K$ une uniformisante de $\O_K$, le lemme
de Hensel assure que $L$ s'obtient en ajoutant à $K$ une racine $n$-ième
de $\pi_K$. Cette racine $n$-ième est en outre une uniformisante de $L$,
et c'est elle que nous choisissons pour $\pi_L$. L'extension $L/K$ est
galoisienne et son groupe de Galois $\Gal(L/K)$ s'identifie au groupe
des racines $n$-ièmes de l'unité par l'application $\sigma \mapsto
\frac{\sigma(\pi_L)}{\pi_L}$. Soit encore $\pi_{1,L}$ une racine
$p$-ième de $\pi_L$. On pose $\pi_{1,K} = \pi_{1,L}^n$ ; c'est bien une
racine $p$-ième de $\pi_K$.

Nous notons $u_K$ (resp. $u_L$) la variable intervenant dans les
polynômes éléments de $\tilde S_K$ (resp $\tilde S_L$) et $\hat A_K =
\hat A_0 \brac{X_K}$ (resp. $\hat A_L = \hat A_0 \brac{X_L}$) l'anneau
de périodes associé. On dispose d'une inclusion $\tilde S_K \toinj
\tilde S_L$, $u_K \mapsto u_L^n$ qui fait de $\tilde S_L$ un $\tilde
S_K$-module libre de rang $n$. On a également une flèche $\psi_{K,L} :
\hat A_K \to \hat A_L$ défini comme l'unique application $\hat
A_0$-linéaire envoyant $\gamma_i(X_K)$ sur $\gamma_i((1+X_L)^n - 1)$
pour tout $i \geq 0$ ; c'est un isomorphisme d'anneaux $G_L$-équivariant
d'inverse $\psi_{L,K}$ défini comme l'unique application $\hat
A_0$-linéaire envoyant $\gamma_i(X_L)$ sur $\gamma_i((1+X_K)^m - 1)$
pour tout $i \geq 0$. Le diagramme suivant est commutatif :
$$\xymatrix @C=70pt {
\tilde S_K \ar[d]_-{u_K \mapsto u_L^n} \ar[r]^-{u_K \mapsto 
\frac{\pi_{1,K}}{1+X_K}} & \hat A_K \ar@<2pt>[d]^-{\psi_{K,L}} \\
\tilde S_L \ar[r]_-{u_L \mapsto \frac{\pi_{1,L}}{1+X_L}} & \hat A_L
\ar@<2pt>[u]^-{\psi_{L,K}} }$$
L'extension des scalaires de $\tilde S_K$ à $\tilde S_L$ définit de 
façon évidente un foncteur exact et fidèle $\ModpN {\tilde S_K} \to 
\ModpN {\tilde S_L}$.

\begin{prop}
Soit $\calM_K \in \ModpN {\tilde S_K}$. Alors le morphisme 
$$\Tkst(\calM_K) \to \Tlst(\calM_K \otimes_{\tilde S_K} \tilde S_L),
\quad f \mapsto (\psi_{K,L} \circ f) \otimes_{\tilde S_K} \tilde S_L$$
est un isomorphisme $G_L$-équivariant.
\end{prop}

\begin{proof}
On vérifie directement la $G_L$-équivariance et l'injectivité du
morphisme de la proposition. La surjectivité résulte alors de ce
que les espaces de départ et d'arrivée sont des $\F_p$-espaces 
vectoriels de même dimension (en l'occurrence le rang de $\calM_K$
sur $\tilde S_K$).
\end{proof}

Nous souhaitons à présent décrire les représentations de $G_K$ dont la 
restriction à $G_L$ provient d'un objet de $\ModpN {\tilde S_L}$. Pour 
cela, on a besoin au préalable d'étendre l'action galoisienne sur $\hat 
A_L$ à tout $G_K$. Ceci se fait tout simplement en utilisant 
l'isomorphisme $\psi_{L,K}$ (qui est déjà, rappelons-le, 
$G_L$-équivariant). De façon concrète, $G_K$ agit de façon habituelle 
sur $\hat A_0$ et sur $X_L$ par la formule 
$$\sigma X_L = \pa{\frac{\sigma \pi_{1,K}}{\pi_{1,K}}}^m (1+X_L) - 1$$ 
valable pour tout $\sigma \in G_K$. En particulier, l'action de $G_K$ 
sur $\frac{\pi_{1,L}}{1+X_L}$ n'est pas triviale, mais se fait 
\emph{via} le caractère $\omega : \Gal(L/K) \to k^\star$ défini par
$$\omega(\bar \sigma) =  \pa{\frac{\sigma \pi_{1,L}}{\pi_{1,L}}} \cdot 
\pa{\frac{\sigma \pi_{1,K}}{\pi_{1,K}}}^{-m} = \pa{\frac{\sigma 
\pi_{1,L}}{\pi_{1,L}}}^{1-nm}$$
où $\sigma \in G_L$ relève $\bar \sigma$.
La formule précédente a bien un sens car, d'une part, la valeur du
membre de droite ne dépend pas du relevé choisi, et d'autre part, par
définition de $m$, l'exposant $1-nm$ est multiple de $p$, ce qui assure
que $\omega$ prend ses valeurs dans le groupe des racines $n$-ièmes de
l'unité de $\bar K$ qui sont par hypothèse toutes dans $\O_K^\star$ (et
que l'on identifie ensuite aux racines $n$-ièmes de l'unité de $k^\star$
grâce au lemme de Hensel). Ceci nous conduit à définir une action de
$\Gal(L/K)$ sur $\tilde S_L$ en décrétant qu'il agit trivialement sur
$k$ et par l'intermédiaire de $\omega$ sur $u_L$. Le morphisme habituel
$\tilde S_L \to \hat A_L$, $u_L \mapsto \frac{\pi_{1,L}}{1+X_L}$ est
alors $G_K$-équivariant.

\begin{deftn}
Soit $\calM_L$ un objet de $\ModpN {\tilde S_L}$. Une \emph{donnée de
descente} (de $L$ à $K$) sur $\calM_L$ est une action semi-linéaire de 
$\Gal(L/K)$ sur $\calM_L$ respectant $\Fil^r \calM_L$ et commutant à 
$\phi_r$ et $N$.

On note $\ModpNdd {\tilde S_L}$ la catégorie dont les objets sont la
donnée de $\calM_L \in \ModpN {\tilde S_L}$ et d'une donnée de descente
sur $\calM_L$.
\end{deftn}

\noindent
La catégorie $\ModpNdd {\tilde S_L}$ est additive et équipée d'une 
dualité obtenue en définissant sur $\calM_L^\vee = \hom_{\tilde S_L} 
(\calM_L, \tilde S_L)$ une action de $\Gal(L/K)$ par la formule $(\sigma 
f)(x) = \sigma f (\sigma^{-1} x)$ (pour $\sigma \in \Gal(L/K)$, $f \in 
\calM_L^\vee$ et $x \in \calM_L$). En outre, si $\calM_L \in \ModpNdd 
{\tilde S_L}$, la $G_L$-représentation $\Tlst(\calM_L) = \hom_{\pFilpN 
{\tilde S_L}}(\calM_L, \hat A_L)$ se prolonge naturellement à $G_K$ par 
la même formule que précédemment : $\sigma f(x) = \sigma f (\bar 
\sigma^{-1} x)$ où $\bar \sigma$ est l'image de $\sigma$ dans 
$\Gal(L/K)$. On définit comme ceci un foncteur exact et fidèle 
$\ModpNdd {\tilde S_L} \to \Rep_{\F_p} (G_K)$ noté encore $\Tlst$.

\begin{theo}
\label{theo:pylonetdd}
Le foncteur $\Tlst : \ModpNdd {\tilde S_L} \to \Rep_{\F_p} (G_K)$ est un 
pylonet additif et autodual. En outre, si $\calM_L$ est un objet de 
$\ModpNdd {\tilde S_L}$, l'action de $\Gal(L/K)$ s'étend à 
$\Max(\calM_L)$ (resp. $\Min(\calM_L)$) calculé dans $\ModpN {\tilde 
S_L}$ et en fait un objet de $\ModpNdd {\tilde S_L}$ qui s'identifie
à $\Max(\calM_L)$ (resp. $\Min(\calM_L)$) calculé dans $\ModpNdd {\tilde
S_L}$.

De plus, les restrictions de $\Tlst$ à $\Max(\ModpNdd {\tilde S_L})$ 
et $\Min(\ModpNdd {\tilde S_L})$ sont exactes et pleinement fidèles, et 
leur image essentielle est stable par sous-objets et quotients.
\end{theo}

\begin{proof}
Elle est semblable à celle des théorèmes \ref{theo:pylonetE} et 
\ref{theo:pleinfidE}.
\end{proof}

\noindent
{\it Remarque.} Bien entendu, on peut aussi fabriquer des catégories en 
administrant simultanément des données de descente et l'action de 
coefficients. Le théorème précédent se généralise directement à cette
situation composite.

\subsubsection*{Quelques mots sur le cas général}

Lorsque l'extension $L/K$ est une extension galoisienne quelconque, les 
données de descente sur les objets de $\ModpN {\tilde S_L}$ ont été 
définies dans \cite{bcdt}, \S 5.6. Hélas, dans cette situation plus 
générale, on ne peut en général pas relever de façon canonique l'action 
de $\Gal(L/K)$ au niveau de $\hat A_L$ --- ni même au niveau de $\tilde 
S_L$ --- car on ne dispose plus de l'isomorphisme $\psi_{K,L}$. Il 
est alors nécessaire de faire des choix arbitraires, ce qui impose de 
manipuler toute une flopée de conditions de compatibilités pas vraiment 
agréables. Malgré tout, il est probable qu'il subsiste un énoncé 
analogue à celui du théorème \ref{theo:pylonetdd} dans ce contexte plus 
général.

\subsection{Quotients de réseaux}
\label{subsec:reseaux}

Nous nous intéressons ici aux $E$-représentations qui peuvent s'écrire 
comme un quotient (annulé par $p$) de deux réseaux dans une 
$W(E)[1/p]$-représentation semi-stable dont les poids de Hodge-Tate sont 
dans $\{0, \ldots, r\}$. Pour expliquer le lien avec la théorie que nous 
avons développée dans les pages précédentes, nous avons besoin dans un 
premier temps d'introduire la notion de \emph{module fortement 
divisible} dûe à Breuil.

Soit $S$ le complété $p$-adique de l'enveloppe à puissances divisées
(compatibles aux puissances divisées canoniques sur $p$) de $W[u]$ par
rapport au noyau de $s : W[u] \to \O_K$, $u \mapsto \pi$. Il est 
muni :
\begin{itemize}
\item d'une filtration $\Fil^i S$ définie comme le complété $p$-adique
de la filtration donnée par les puissances divisées ;
\item d'un Frobenius $\phi : S \to S$ défini comme l'unique morphisme
d'anneaux, continu pour la topologie $p$-adique, qui agit sur $K_0$ 
comme le Frobenius et qui envoie $u$ sur $u^p$ ;
\item d'un opérateur de monodromie $N : S \to S$ défini comme l'unique
application continue $W$-linéaire qui envoie $u^n$ sur $-n u^n$.
\end{itemize}
Pour $i < p-1$ (et, donc, en particulier pour $i=r$), on a $\phi(\Fil^i
\calM) \subset p^i \calM$, ce qui permet de définir l'application
$\phi_i = \frac{\phi}{p^i} : \Fil^i \calM \to \calM$. On remarque que
le polynôme minimal de $\pi$ sur $K_0$, traditionnellement noté $E(u)$
est un polynôme d'Eisenstein d'où on déduit que $\phi_1(E(u))$ est une 
unité de $S$. On dispose en outre d'un morphisme évident $S \to \tilde 
S$, $u \mapsto u$, $\gamma_i(u^e) \mapsto 0$ pour $i \geq p$. Il permet 
de voir $\tilde S$ comme une $S$-algèbre et se factorise par $S_1 = 
S/pS$. Un module fortement divisible est alors la donnée des points 
suivants :
\begin{enumerate}
\item un $S$-module libre de rang fini $\calM$ ;
\item un sous-module $\Fil^r \calM \subset \calM$ contenant $\Fil^r S
\: \calM$ ;
\item un opérateur $\phi$-semi-linéaire $\phi_r : \Fil^r \calM \to 
\calM$ vérfiant
$$(\forall s \in S) \, (\forall x \in \calM) \quad
\phi_r(sx) = \phi_r(s) \cdot \frac{\phi_r(E(u)^r x)}{\phi_1(E(u))^r}$$
et dont l'image engendre exactement $p^r \calM$ ;
\item un opérateur $N : \calM \to \calM$ vérifiant :
\begin{itemize}
\item (condition de Leibniz) $N(sx) = s N(x) + N(s) x $ pour tout $x \in
\calM$ et $s \in S$ ;
\item (transversalité de Griffith) $E(u) N(\Fil^r \calM) \subset \Fil^r
\calM$ ;
\item le diagramme suivant est commutatif :
\begin{equation}
\label{eq:deffortdiv}
\raisebox{0.5\depth}{\xymatrix @C=50pt {
\Fil^r \calM \ar[r]^-{\phi_r} \ar[d]_-{E(u) N} & \calM 
\ar[d]^-{\phi_1(E(u)) N} \\
\Fil^r \calM \ar[r]^-{\phi_r} & \calM }}
\end{equation}
\end{itemize}
\end{enumerate}
On note $\ModpN S$ la catégorie des modules fortement divisibles, les
morphismes étant naturellement les applications $S$-linéaires commutant
aux structures supplémentaires. On définit de même la catégorie $\ModpN 
{S_1}$ en remplaçant partout $S$ par $S_1 = S/pS$. Une adaptation 
immédiate de la proposition 2.2.2.1 de \cite{breuil-ens} montre que le 
fonction $T : \calM \mapsto \calM \otimes_{S_1} \tilde S$ donne 
naissance à une équivalence de catégories entre $\ModpN {S_1}$ et 
$\ModpN {\tilde S}$ dont un quasi-inverse est donné par la formule 
$$T^{-1}(\calM) = S_1 \otimes_{(\phi),k[u]/u^e} \frac{\Fil^r \calM} 
{u^e \Fil^r \calM}.$$
D'autre part, on dispose d'un foncteur $\Tsthat : \ModpN S \to 
\Rep_{\Z_p} (G_K)$ dont la définition est analogue à celle de $\Tst$ 
mais fait intervenir un anneau de période plus compliqué que nous ne 
souhaitons pas décrire ici. Quoi qu'il en soit, dans \cite{liu-compo}, 
Liu a montré que $\Tsthat$ induit une anti-équivalente entre $\ModpN S$ 
et la catégorie des réseaux dans les représentations semi-stables à 
poids de Hodge-Tate compris entre $0$ et $r$. Finalement, on montre 
qu'un morphisme surjectif $\hat \calM \to \calM$ (avec $\hat \calM \in 
\ModpN S$ et $\calM \in \ModpN {S_1}$) induit 
une surjection $\Tsthat(\hat \calM) \to \Tst \circ T(\calM)$ et donc 
fait apparaître $\Tst \circ T (\calM)$ comme un quotient d'un réseau 
dans une représentation semi-stable.

\begin{lemme}
\label{lem:relevT}
Soient $\hat \calM \in \ModpN S$, $\calM \in \ModpN {S_1}$. Notons
$\pr : \calM \to T(\calM)$ la projection canonique. Soit $f$
un morphisme $S$-linéaire $\hat \calM \to T(\calM)$ compatible
aux structures additionnelles. Alors, il existe un unique morphisme
$S$-linéaire et compatible aux structures addtionnelles $g : \hat 
\calM \to \calM$ tel que $f = \pr \circ g$.
De plus, $g$ est surjectif si, et seulement si $f$ l'est.
\end{lemme}

\begin{proof}
L'unicité de $g$ est facile et laissée au lecteur. Pour l'existence, on 
remarque d'abord que $f$ passe au quotient pour définir un morphisme 
$\bar f : \hat \calM \otimes_S \tilde S \to T(\calM)$ dans la catégorie 
$\ModpN {\tilde S}$. L'image de $\bar f$ par $T^{-1}$ est alors un 
morphisme de $\hat \calM / p\hat \calM$ dans $\calM$, qui composé avec 
la projection $\hat \calM \to \hat \calM / p \hat\calM$ fournit un $g$ 
adéquat.
Évidemment si $g$ est surjectif, $f$ l'est aussi. Réciproquement si $f$
est surjectif, $g \otimes_S W = f \otimes_S W$ l'est aussi, ce qui suffit
à assurer la surjectivitié de $g$ lui-même.
\end{proof}

\noindent
{\it Remarque.} On montre de même qu'un morphisme $f : \hat \calM \to
\RedN(T(\calM))$ se relève en une unique flèche $g : \hat \calM \to 
\calM$.

\bigskip

Tout cela nous conduit à poser la définition suivante :

\begin{deftn}
\label{def:modst}
On note $\Modst {\tilde S \otimes E}$ la sous-catégorie pleine de
$\ModpN {\tilde S \otimes E}$ formée des objets $(\calM,\nu)$ pour
lesquels il existe un module fortement divisible $\hat \calM$, 
un morphisme de $\Z_p$-algèbres $\hat \nu : W(E) \to \End_{\ModpN S} 
(\hat \calM)$ et un morphisme surjectif $S$-linéaire compatible à toutes 
les structures $f : \hat \calM \to \calM$ tels que pour tout $\hat 
\lambda \in W(E)$, $\hat \nu(\hat \lambda)$ stabilise $\ker f$ et 
induise sur $\calM$ l'application $\nu (\lambda)$ où $\lambda$ est la 
réduction de $\hat \lambda$ modulo $p$.
\end{deftn}

\begin{lemme}
\label{lem:filsurj}
Soient $\hat \calM \in \ModpN S$, $\calM \in \ModpN {S_1}$ et $f : \hat 
\calM \to \calM$ un morphisme surjectif compatible à $\Fil^r$ et 
$\phi_r$. Alors le morphisme $\Fil^r \hat \calM \to \Fil^r \calM$ induit 
par $f$ est surjectif.
\end{lemme}

\begin{proof}
Soit $\calF = f(\Fil^r \hat \calM)$. De la surjectivité de $f$, on 
déduit que le module engendré par $\phi_r(\calF)$ est $\calM$ tout 
entier. L'isomorphisme
$$S_1 \otimes_{(\phi),k[u]/u^e} \frac{\Fil^r 
\calM} {u^e \Fil^r \calM + \Fil^p S_1 \calM} \simeq
S_1 \otimes_{(\phi),k[u]/u^e} \frac{\Fil^r T(\calM)} {u^e \Fil^r 
T(\calM)} \stackrel{\sim}{\longrightarrow} \calM$$
montre alors que $\Fil^r \calM = \calF + u^e \Fil^r \calM + \Fil^p S_1
\calM$. Or on a $\Fil^p S_1 \calM \subset \Fil^r S_1 \calM =
\phi_r(\Fil^r S \hat \calM) \subset \calF$, ce qui donne $\Fil^r \calM =
\calF + u^e \Fil^r \calM$. La conclusion s'ensuit facilement en
remarquant que la suite des puissances de $u^e$ (dans $S_1$) s'annule à 
partir d'un certain rang (en l'occurrence $u^{ep}$).
\end{proof}

\begin{theo}
La restriction de $\Tst$ à $\Modst {\tilde S \otimes E}$ est un pylonet
additif et autodual.
\end{theo}

\begin{proof}
Il faut vérifier les axiomes \Ax 1, \Ax 2, \Ax {3a}, \Ax {3b}, \Ax 4 et
\Ax 5. Éventuel\-lement en utilisant les énoncés analogues pour le
foncteur $\Tst$ défini sur la catégorie $\ModpN {\tilde S \otimes E}$
tout entière, on établit facilement \Ax 1, \Ax {3b} et \Ax 5. L'axiome
\Ax {3a} ne pose pas non plus véritablement problème : en reprenant les
notations de la démonstration de la proposition \ref{prop:axiome3a}, il
suffit de montrer que si $\calM_1$ et $\calM_2$ sont dans $\Modst
{\tilde S \otimes E}$, alors il en est de même de $\calM'$, ce qui est
clair puisque $\calM'$ est défini comme un quotient de $\calM_1 \oplus
\calM_2$. La vérification de \Ax 2 est, elle aussi, très simple : il
suffit de justifier que si $f : \calM \to \calM'$ est un morphisme dans
$\Modst {\tilde S \otimes E}$, alors $\calM'' = \coker f$ (calculé dans
$\ModpN {\tilde S \otimes E}$) est dans $\Modst {\tilde S \otimes E}$ et
pour cela, d'après le lemme \ref{lem:relevT}, il suffit d'établir la 
surjectivité de $g : \calM' \to \calM''$... qui résulte directement de 
la construction.

Il ne reste finalement qu'à vérifier \Ax 4. On suppose pour simplifier
que $E = \F_p$, le cas général s'obtenant de la même façon en ajoutant
l'action de $E$ ou de $W(E)$ à chaque étape. Soit $\calM \in \Modst
{\tilde S}$. Par hypothèse, il existe $\hat \calM \in \ModpN S$ muni
d'un morphisme surjectif $f : \hat \calM \to \calM$. Par le lemme
\ref{lem:relevT}, celui-ci se relève en un morphisme surjectif $g : \hat
\calM \to T^{-1}(\calM)$ et par le lemme \ref{lem:filsurj}, $g$ induit
aussi une surjection au niveau des $\Fil^r$. En utilisant les
équivalences de catégories données par les théorèmes 2.2.1 et 2.3.1 de
\cite{carliu}, on montre aisément que $\hat \calM' = \ker g$ est encore
un objet de $\ModpN S$. Ainsi on obtient une suite exacte $0 \to \hat
\calM' \to \hat \calM \to T^{-1}(\calM) \to 0$ qui induit également une
suite exactement au niveau des $\Fil^r$. Le lemme V.3.4.1 de
\cite{caruso-thesis} montre alors l'existence d'une nouvelle suite
exacte\footnote{Dans \emph{loc. cit.}, on utilise la notation \og $\vee$
\fg\ à la place de \og $\star$ \fg.} :
$$0 \to \hat \calM^\star \to \hat \calM'^\star \to T^{-1}(\calM^\star) \to 
0$$
à partir de laquelle on obtient le morphisme surjectif que l'on 
cherchait.
\end{proof}

Notons $\Maxst$ et $\Minst$ les endofoncteurs de $\Modst {\tilde S
\otimes E}$ qui découlent du théorème précédent. Notez bien que si
$\calM$ est dans $\Modst {\tilde S \otimes E}$, c'est aussi un objet de
$\ModpN {\tilde S \otimes E}$ auquel on peut donc appliquer les deux
foncteurs $\Max$ et $\Maxst$. Je ne sais pas si ces foncteurs
coïncident en général, et c'est la raison pour laquelle je préfère
introduire deux notations distinctes. Définissons également 
$\Maxst_{/\tilde S \otimes E}$ (resp. $\Minst_{/\tilde S \otimes E}$)
comme l'image essentielle de $\Maxst$ (resp. $\Minst$).

\begin{theo}
La restriction de $\Tst$ à $\Maxst_{/\tilde S \otimes E}$ (resp.
$\Maxst_{/\tilde S \otimes E}$) est pleinement fidèle et son image
essentielle est stable par sous-objets et quotients.
\end{theo}

\begin{proof}
On ne traite que le cas des objets maximaux, l'autre s'obtenant par
dualtité. La fidélité ne pose aucun problème. Soient $\calM$ et $\calM'$
des objets de $\Maxst_{/\tilde S \otimes E}$ et $f : \TstE(\calM') \to
\TstE(\calM)$ une application $G_K$-équivariante. Par le théorème
\ref{theo:pleinfidE}, $f$ provient d'un morphisme $g : \Max(\calM) \to
\Max(\calM')$. Par ailleurs, par le lemme \ref{lem:redinjectif},
$\RedN(\Maxst(\calM))$ et $\RedN(\Maxst(\calM'))$ apparaissent
respectivement comme des sous-modules de $\RedN(\Max(\calM))$ et
$\RedN(\Max(\calM'))$. Pour établir la pleine fidélité, il suffit de
montrer que $\RedN(g)$ envoie $\RedN \circ \Maxst(\calM)$ sur $\RedN 
\circ \Maxst(\calM')$. Considérons $\hat \calM$ et $\hat \calM'$ des 
modules fortement divisibles munis de surjections $h : \hat \calM \to 
\RedN \circ \Maxst (\calM)$ et $h' : \hat \calM' \to \RedN \circ 
\Maxst(\calM')$ et attardons-nous sur le morphisme
$$(\RedN(g) \circ h) \oplus h' : \hat \calM \oplus \hat \calM'
\to \RedN \circ \Max(\calM').$$
Soit $\calM''$ son image. Par la proposition \ref{prop:redinterm},
$\TstE(\calM'') \simeq \TstE \circ \Maxst(\calM)$ (la compatibilité à
l'action de $E$ venant de la fonctorialité), et donc par maximalité de
$\Maxst(\calM')$, on obtient $\calM'' \subset \Red \circ
\Maxst(\calM')$. Ceci entraîne $\Red(g) \circ h (\hat \calM) \subset
\Red \circ \Maxst(\calM')$, \emph{i.e.} $\Red(g)(\Red \circ
\Maxst(\calM)) \subset \Red \circ \Maxst(\calM')$ comme voulu.

Lorsque $E = \F_p$, la stabilité par sous-objets découle directement de
la proposition \ref{prop:decoup}. La stabilité par quotients s'obtient
par dualité, tandis que le cas des coefficients quelconques se fait en
relevant l'action de $E$ grâce à la pleine fidélité.
\end{proof}

\subsubsection*{Représentations cristallines} 

On peut également s'intéresser aux réseaux à l'intérieur de 
repré\-sentations cristallines plutôt que semi-stables ; au niveau des
modules fortement divisibles, ceci correspond à $N \equiv 0 \pmod {uS + 
\Fil^p S}$, c'est-à-dire $N(\hat \calM) \subset (uS + \Fil^p S) \hat
\calM$. En réalité, on peut légèrement simplifier cette condition comme 
l'affirme le lemme suivant.

\begin{lemme}
Soit $\hat \calM$ un module fortement divisible tel que $N(\hat \calM) 
\subset (uS + \Fil^p S) \hat \calM$. Alors $N(\hat \calM) \subset u \hat 
\calM$.
\end{lemme}

\begin{proof}
Le diagramme \eqref{eq:deffortdiv} montre que $N \circ \phi_r (\Fil^r 
\hat \calM)$ est inclus dans $\phi_r(\Fil^r \hat \calM)$. Pour estimer ce 
dernier, on utilise la proposition 4.1.2 de \cite{liu-compo} qui assure 
l'existence de $x_1, \ldots, x_d \in \Fil^r \hat \calM$ tels que
\begin{itemize}
\item[$\bullet$] $Fil^r \hat \calM$ est engendré par les $x_i$ et 
$\Fil^p S \: \hat\calM$ ;
\item[$\bullet$] les $e_i = \phi_r(x_i)$ forment une base de $\hat 
\calM$ ;
\item[$\bullet$] les $E(u)^r e_i$ s'expriment comme une combinaison 
linéaire à coefficients dans $S$ des $x_i$.
\end{itemize}
Ainsi, en définissant 
$$T = \Bigg\{ \sum_{i \geq 0} a_i \frac{E(u)^i}{(i+r)!} , \quad 
a_i \in W[u], \, \lim_{i \to \infty} a_i = 0 \Bigg\}$$
on a l'inclusion $\Fil^r \hat \calM \subset T x_1 + \cdots + T x_r$ (où
tout est vu par exemple dans $\hat \calM \otimes_S K_0[[u]]$). Par
suite, $\phi_r(\Fil^r \hat \calM)$ est contenu dans le $\phi(T)$-module
engendré par les $e_i$.

Montrons que $\phi(T) \cap (uS + \Fil^p S) \subset u S$. On considère
pour cela un élément $x$ dans l'intersection précédente, et on souhaite
montrer qu'il est dans $uS$. On peut écrire $x = \phi(y)$ avec $y =
\sum_{i \geq 0} a_i \frac{E(u)^i}{(i+r)!}$ où $a_i \in W[u]$ converge
vers $0$. En regardant modulo $uS + \Fil^p S$, on obtient $\sum_{i \geq
0} a_i(0) \frac{E(0)^i}{(i+r)!} = 0$, puis
\begin{eqnarray*}
x & = &
\sum_{i \geq 0} \phi(a_i) \frac{\phi(E(u))^i - \phi(E(0))^i}{(i+r)!} +
\sum_{i \geq 0} [\phi(a_i) - \phi(a_i(0))] \frac{\phi(E(0))^i}{(i+r)!} \\
& = & [\phi(E(u)) - \phi(E(0))] \sum_{s,t \geq 0} \phi(a_{s+t+1})
\frac{\phi(E(u))^s \phi(E(0))^t}{(s+t+r+1)!} +
\sum_{i \geq 0} [\phi(a_i) - \phi(a_i(0))] \frac{\phi(E(0))^i}{(i+r)!}
\end{eqnarray*}
Étant donné que $u$ --- et même à vrai dire $u^p$ --- divise $\phi(E(u))
- \phi(E(0))$ et $\phi(a_i) - \phi(a_i(0))$, il suffit de justifier que
$\frac{\phi(E(u))^s \phi(E(0))^t}{(s+t+r+1)!} \in S$ pour tous entiers
$s$ et $t$. Le numérateur de cette dernière fraction est manifestement
divisible par $p^{s+t}$. De plus,
$$v_p\pa{\frac{p^{s+t}}{(s+t+r+1)!}} > s+t - \frac{s+t+r+1}{p-1} = 
\frac{p-2}{p-1}(s+t) - \frac{r+1}{p-1} \geq - \frac{r+1}{p-1}
\geq -1$$
d'où on déduit que $\frac{p^{s+t}}{(s+t+r+1)!} \in W \subset S$. Au
final, $\frac{\phi(E(u))^s \phi(E(0))^t}{(s+t+r+1)!} \in S$ comme voulu.

On conclut maintenant la preuve du lemme comme suit. Par ce qui précède,
on a $N \circ \phi_r(\Fil^r \hat \calM) \subset u \hat \calM$, et donc
en particulier $N(e_i) \in u \hat \calM$ pour tout $i \in \{1, \ldots,
x_d\}$. Par ailleurs, on vérifie tout de suite que pour tout $s \in S$,
$N(s)$ est divisible par $u$. Ainsi $N(s e_i) = N(s) e_i + s N(e_i)$ 
est lui aussi multiple de $u$. La valeur de $N$ sur n'importe quelle
combinaison linéaire des $e_i$ est donc multiple de $u$. Comme les $e_i$
forment une base de $S$, on a bien démontré que $N(\hat \calM) \subset
\hat \calM$.
\end{proof}

On peut alors adapter la définition \ref{def:modst} dans ce contexte :

\begin{deftn}
\label{def:modcris}
On note $\Modcris {\tilde S \otimes E}$ la sous-catégorie pleine de
$\ModpN {\tilde S \otimes E}$ formée des objets $(\calM,\nu)$ pour
lesquels il existe un module fortement divisible $\hat \calM$ avec 
$N(\hat \calM) \subset u \hat \calM$, un morphisme de $\Z_p$-algèbres 
$\hat \nu : W(E) \to \End_{\ModpN S} (\hat \calM)$ et un morphisme 
surjectif $S$-linéaire compatible à toutes les structures $f : \hat 
\calM \to \calM$ tels que pour tout $\hat \lambda \in W(E)$, $\hat 
\nu(\hat \lambda)$ stabilise $\ker f$ et induise sur $\calM$ 
l'application $\nu (\lambda)$ où $\lambda$ est la réduction de $\hat 
\lambda$ modulo $p$.
\end{deftn}

\begin{theo}
La restriction de $\Tst$ à $\Modcris {\tilde S \otimes E}$ est un
pylonet additif et autodual. La restriction de $\Tst$ à la catégorie
des objets maximaux (resp. minimaux) correspondants est pleinement 
fidèle et son image essentielle est stable par sous-objets et
quotients.
\end{theo}

\begin{proof}
C'est la même que dans le cas semi-stable.
\end{proof}

De façon similaire, on peut considérer la sous-catégorie pleine de 
$\ModpN {\tilde S \otimes E}$ comprenant les objets $\calM$ pour 
lesquels $N(\calM) \subset u \calM$ (sans demander, donc, qu'il existe 
un relèvement sous forme de module fortement divisible). Par les mêmes 
méthodes, on a encore un théorème analogue dans cette dernière 
situation.

\subsection{Objets simples}

On suppose dans cette sous-section $er \geq p-1$ (le cas $er < p-1$ a
déjà été étudié dans \cite{caruso-crelle}). On note $K^\nr \subset \bar
K$ l'extension maximale non ramifiée de $K$. Son corps résiduel
s'identifie à une clôture algébrique de $k$, notée $\bar k$. Pour tout
entier $d$, on note $\F_{p^d}$ l'unique sous-corps de $\bar k$ de
cardinal $p^d$. On fixe par ailleurs $E$ une extension finie de $\F_p$
de degré $h$, ainsi qu'un isomorphisme $\tau : E \to \F_{p^h}$. Dans la
suite, on supposera toujours que l'image de $\tau$ est incluse dans $k$
et on utilisera cette hypothèse pour identifier $E$ à un sous-corps de
$k$.

\medskip

Soit $\calR_h$ l'ensemble des classes d'équivalence d'éléments de
$\Z_{(p)}$ (le localisé de $\Z$ en $p$) pour la relation d'équivalence
suivante : $a \sim b$ si, et seulement s'il existe un entier $n$ tel que
$a \equiv p^{hn} b \pmod \Z$. \emph{Via} l'écriture en base $p$, les
éléments de $\calR_h$ s'identifient à l'ensemble des suites $(a_i)$
périodiques (depuis le début) d'entiers compris entre $0$ et $p-1$ où on
a identifié la suite $(a_i)$ à la suite $(a_{i+h})$, et où on a ôté la
suite constante égale à $p-1$.
À tout $a \in \calR_h$, on associe un objet $(\calM(a),\nu_a)$ de
$\MaxpN {\tilde S \otimes E}$ défini comme suit. On choisit $(a_i)$ une
suite périodique qui représente $a$, on note $d$ sa plus petite période, 
$d_h = \ppcm(d,h)$ et on définit :
\begin{itemize}
\item[$\bullet$] $\calM(a) = \bigoplus_{i \in \Z/d_h\Z} \tilde S \cdot
e_i$ ;
\item[$\bullet$] $\Fil^r \calM(a) = \sum_{i \in \Z/d_h\Z} \tilde S \cdot
u^{er-a_i} e_i$ ;
\item[$\bullet$] $\phi_r(u^{er-a_i} e_i) = (-1)^r e_{i+1}$ ;
\item[$\bullet$] $N(e_i) = 0$ ;
\item[$\bullet$] $\nu_a(\lambda)(e_i) = \lambda^{p^i} e_i$ ($\lambda
\in E \subset k$).
\end{itemize}
À partir de la proposition 3.6.7 de \cite{carliu}, on montre facilement
que $\calM(a)$ est un objet de $\MaxpN {\tilde S \otimes E}$. De plus, 
on vérifie sans mal qu'il ne dépend pas (à isomorphisme près) du choix 
du représentant $(a_i)$. 

On peut en outre déterminer la restriction au groupe d'inertie, noté
$I_K$, de la représentation galoisienne associée à $\calR_h$. Pour cela,
on a tout d'abord besoin de rappeler la définition des caractères
fondamentaux de Serre. Pour tout entier $d$, on définit $\theta_d : I_K 
\to \mu_{p^d-1} (\bar K) \simeq \F_{p^d}^\star$, $g \mapsto \frac{g 
\pi_d} {\pi_d}$, l'isomorphisme entre $\mu_{p^d-1} (\bar K)$ et 
$\F_{p^d}^\star$ étant induit par la réduction modulo l'idéal maximal.
(On rappelle que $\pi_d$ une racine $p^d$-ième fixée de $\pi$.)

\begin{prop}
\label{prop:tstsimple}
Soient $a \in \calR_h$, $(a_i)$ une suite périodique représentant $a$ et
$d$ sa plus petite période. Alors, en tant que $E$-représentation de
$I_K$ :
$$\TstE(\calM(a),\nu_a) = \F_{p^d}\big(\theta_d^{a_0 + p^{d-1} a_1 + 
\cdots + p a_{d-1}}\big) \otimes_{\F_{p^d} \cap E} E$$
où $\F_{p^d}(\psi)$ désigne la $\F_{p^d}$-représentation de $G_K$ de
dimension $1$ où l'action se fait par le caractère $\psi$.
\end{prop}

\begin{proof}
Elle est semblable à celle du théorème 5.2.2 de \cite{caruso-crelle} ;
nous nous contentons donc de renvoyer à cette référence. On prendra
toutefois garde au twist qui apparaît dans la définition de $\hat A_0$
qui n'est pas discuté avec beaucoup d'attention dans \emph{loc. cit.},
et peut facilement être source d'erreurs dans les calculs (l'erreur se
manifestant le plus souvent par un décalage d'indice dans l'exposant
de $\theta_d$).
\end{proof}

\begin{cor}
On suppose $er \geq p-1$ et $k$ algébriquement clos.
Les objets simples de $\MaxpN {\tilde S \otimes E}$ sont exactement
les $\calM(a)$, $a \in \calR_h$. De plus, ils sont deux à deux non
isomorphes.
\end{cor}

\begin{proof}
Les arguments des paragraphes 1.6 et 1.7 de \cite{serre} montrent que
les $E$-représentations irréduc\-tibles de $G_K = I_K$ sont exactement 
les $\TstE(\calM(a), \nu_a)$ et qu'elles sont deux à deux non 
isomorphes. Le corollaire provient alors de la pleine fidélité de $\TstE 
: \MaxpN {\tilde S \otimes E} \to \Rep_E (G_K)$ (théorème 
\ref{theo:pleinfidE}).
\end{proof}




\end{document}